\documentstyle[11pt]{article}
\input epsf.tex
\font\twelvegoth=eufm10 at 12pt
\font\tengoth=eufm10
\font\sevengoth=eufm7
\font\fivegoth=eufm5

\newfam\gothfam
\textfont\gothfam=\tengoth
\scriptfont\gothfam=\sevengoth
\scriptscriptfont\gothfam=\fivegoth
\def\goth{\fam\gothfam\tengoth}

\font\tensym=msbm10
\font\sevensym=msbm7
\font\fivesym=msbm5

\newfam\ssymfam
\textfont\ssymfam=\tensym
\scriptfont\ssymfam=\sevensym
\scriptscriptfont\ssymfam=\fivesym
\def\ssym{\fam\ssymfam\tensym}

\newtheorem{example}{Example}[section]
\newtheorem{theorem}[example]{Theorem}
\newtheorem{corollary}[example]{Corollary}
\newtheorem{proposition}[example]{Proposition}
\newtheorem{lemma}[example]{Lemma}
\newtheorem{remark}[example]{Remark}

\def\boxit#1#2{\setbox1=\hbox{\kern#1{#2}\kern#1}%
\dimen1=\ht1 \advance\dimen1 by #1 \dimen2=\dp1 \advance\dimen2 by #1
\setbox1=\hbox{\vrule height\dimen1 depth\dimen2\box1\vrule}%
\setbox1=\vbox{\hrule\box1\hrule}%
\advance\dimen1 by .4pt \ht1=\dimen1
\advance\dimen2 by .4pt \dp1=\dimen2 \box1\relax}

\def\resp{{\em resp.$\ $}}
\def\proof{\medskip\noindent {\it Proof --- \ }}

\def\cqfd{\hfill $\Box$ \bigskip}
\def\adots{\mathinner{\mkern2mu\raise1pt\hbox{.}
\mkern3mu\raise4pt\hbox{.}\mkern1mu\raise7pt\hbox{.}}}
\def\<{\lngle\,}
\def\>{\,\rangle}

\def\cf{{\it cf.$\ $}}

\def\ie{{\it i.e. }}
\def\eg{{\it e.g. }}
\def\tab{{\rm Tab\, }}

\def\SG{\hbox{\goth S}}

\def\N{{\bf N}}
\def\Z{{\bf Z}}
\def\C{{\bf C}}
\def\R{{\bf R}}
\def\Q{{\bf Q}\, }
\def\T{{\cal T}}
\def\F{{\cal F}}
\def\P{{\cal P}}
\def\PB{{\bf P}}
\def\IB{{\bf I}}
\def\FB{{\bf F}}
\def\EAHB{\widehat{\bf H}}

\def\U{\widetilde{\cal U}}
\def\V{\widetilde{\cal V}}
\def\UU{{\cal U}}
\def\VV{{\cal V}}
\def\I{{\cal J}}
\def\S2{{\cal S}}

\def\sgn{{\rm sgn\, }}
\def\tr{{\rm Tr}}
\def\mod{{\rm \ mod\ }}
\def\ch{{\rm ch\, }}

\def\h{\hbox{\goth h}}
\def\gl{\hbox{\goth gl\hskip 1pt}}
\def\slchap{\widehat{\hbox{\goth sl}}}

\def\ASG{\widetilde{\hbox{\goth S}}}
\def\EASG{\widehat{\hbox{\goth S}}}
\def\EAH{\widehat H}

\def\A{{\cal A}}
\def\G{{\cal G}}
\def\L{{\cal L}}
\def\desc{{\rm desc}}

\def\bar{\overline}

\def\Fr{{\rm Fr}}
\def\sp{{\rm spin}}
\def\<{\langle}
\def\>{\rangle}
\def\finex{\hfill $\Diamond$}
\def\op{{\rm op}}
\def\Im{{\rm im\,}}
\def\pr{{\rm pr\,}}
\def\PP{{\ssym P}}
\def\E{{\cal E}}
\def\FI{{\bf F}_\infty}
\def\HI{{\bf H}_\infty}

\def\deg{{\rm deg\,}}
\def\today{\number\day \space\ifcase\month\or
 Janvier\or F\'evrier\or Mars\or Avril\or Mai\or Juin\or
 Juillet\or Ao\^ut\or Septembre\or Octobre\or Novembre\or D\'ecembre\fi
 \space\number\year}

\title{\bf Littlewood-Richardson coefficients and Kazhdan-Lusztig polynomials}
\author{\rm Bernard {\sc Leclerc}\thanks{
D\'epartement de Math\'ematiques,
Universit\'e de Caen, Campus II,
Boulevard Mar\'echal Juin,
BP 5186, 14032 Caen cedex, France.}
\rm \ and Jean-Yves {\sc Thibon}\thanks{
Institut Gaspard Monge,
Universit\'e de Marne-la-Vall\'ee,
Champs-sur-Marne,
77454 Marne-la-Vall\'ee cedex 2, France.}
}

\date{}

\begin{document}
\maketitle

\vskip 1cm

\begin{abstract}
We show that the Littlewood-Richardson coefficients are values at 1
of certain parabolic Kazhdan-Lusztig polynomials for affine symmetric groups.
These $q$-analo\-gues of Littlewood-Richardson multiplicities 
coincide with those previously introduced in \cite{LLT2} in terms
of ribbon tableaux.
\end{abstract}

\vskip 0.6cm

\section{Introduction} \label{SECT1}
Let $\lambda = (\lambda_1 \ge \ldots \ge \lambda_r\ge 0)$ and 
$\mu=(\mu_1 \ge \ldots \ge \mu_r\ge 0)$ denote two partitions 
of length $\le r$, identified in the usual way with dominant
integral weights of the complex Lie algebra $\gl_r$.
It was shown by Lusztig \cite{Lu2} that the multiplicity 
$K_{\lambda,\mu}$ of the weight $\mu$ in the finite-dimensional
irreducible representation $W(\lambda)$ of $\gl_r$ with highest 
weight $\lambda$ is the value at 1 of a certain Kazhdan-Lusztig 
polynomial $P_{n_\mu, n_\lambda}$ for the affine symmetric group $\EASG_r$. 
(For the definition of $n_\lambda$, see below Section~\ref{SUBSECT21}).
Moreover, Lusztig proved \cite{Lu1} that the polynomial 
$P_{n_\mu, n_\lambda}(q)$ is equal to the Kostka-Foulkes polynomial 
$K_{\lambda,\mu}(q)$ defined as the coefficient of the Schur function 
$s_\lambda$ on the basis of Hall-Littlewood function $P_\mu(q)$ \cite{Mcd}.
A combinatorial expression of $K_{\lambda,\mu}(q)$ had previously
been given by Lascoux and Sch\"utzenberger in terms of semi-standard
Young tableaux \cite{Sc,Mcd}.

It is well known that $K_{\lambda,\mu}$ is also equal to the
multiplicity of $W(\lambda)$ as an irreducible component
of the tensor product 
\[
W(\mu_1)\otimes \cdots \otimes W(\mu_r)
\]
of symmetric powers of the vector representation of $\gl_r$.
Let now $\nu^{(1)},\ldots ,\nu^{(s)}$ be arbitrary dominant
weights and let $c_{\nu^{(1)},\ldots ,\nu^{(s)}}^\lambda$
denote the multiplicity of $W(\lambda)$ in
\[
W(\nu^{(1)})\otimes \cdots \otimes W(\nu^{(s)})\,.
\]
A $q$-analogue $c_{\nu^{(1)},\ldots ,\nu^{(s)}}^\lambda(q)$
of this multiplicity
has been introduced in \cite{LLT2} by means of certain generalizations 
of semi-standard Young tableaux called ribbon tableaux, and it has been 
proved that when the partitions $\nu^{(j)}$ have only one part $\mu_j$ 
\[
c_{\mu_1,\ldots ,\mu_r}^\lambda(q) = K_{\lambda,\mu}(q).
\]
The purpose of this paper is to establish that for all
$\nu^{(1)},\ldots ,\nu^{(s)},\lambda$ the
$c_{\nu^{(1)},\ldots ,\nu^{(s)}}^\lambda(q)$ are Kazhdan-Lusztig
polynomials for the group $\EASG_r$.

Let us outline how this result is obtained.
As mentioned by Lusztig in \cite{Lu2}, the expression of the
weight multiplicity $K_{\lambda,\mu}$ as a value at 1 of a Kazhdan-Lusztig
polynomial might be deduced from the conjecture of \cite{Lu0}
for the characters of irreducible representations of 
$GL_r$ over an algebraically closed field of characteristic $n\ge r$ 
together with the Steinberg tensor product theorem.  
In \cite{Lu4,Lu5} a similar conjecture was formulated for the
characters of irreducible representations of $U_q(\gl_r)$
when $q^2$ is a primitive $n$th root of 1. 
A remarkable feature of the quantum conjecture is that the 
restriction $n\ge r$ is no longer necessary.
This conjecture is
now proved due to work of Kazhdan-Lusztig and Kashiwara-Tanisaki. 
On the other hand
Lusztig has derived in \cite{Lu4} an analogue of the Steinberg
tensor product theorem for the quantum case. 
From these two facts, it is easy to deduce that the Littlewood-Richardson
multiplicities are value at 1 of Kazhdan-Lusztig polynomials
(see below, Section~\ref{SECT3}).

However this would not provide the link with the $q$-analogues
defined by means of ribbon tableaux. 
We shall therefore follow a different approach and
rely on the construction given in \cite{LT} of a canonical basis 
in the level 1 Fock space representation of the quantum affine 
algebra $U_q(\slchap_n)$. This canonical basis satisfies a 
formal $q$-analogue of Steinberg's tensor product theorem which
may be formulated in terms of the combinatorics of ribbon tableaux.
On the other hand, Varagnolo and
Vasserot \cite{VV} have recently verified a conjecture of 
\cite{LT}. They proved that the coefficients of the expansion 
of this canonical basis on the standard basis of $q$-wedge products
coincide with the Kazhdan-Lusztig polynomials occuring in
Lusztig's conjecture. 
Using these two results, we are able to express the 
$c_{\nu^{(1)},\ldots ,\nu^{(s)}}^\lambda(q)$ as Kazhdan-Lusztig
polynomials.

More precisely, they belong to a family of parabolic analogues of
Kazhdan-Lusztig polynomials introduced by Deodhar \cite{De1,De2}.
There are two types of such polynomials associated with the
Hecke algebra modules obtained by inducing respectively the characters
$T_i \mapsto -q$ and $T_i \mapsto q^{-1}$ of a parabolic
subalgebra.
The $c_{\nu^{(1)},\ldots ,\nu^{(s)}}^\lambda(q)$ turn out
to belong to the  family denoted by $\widetilde P_{x,y}^{J}$ in \cite{De2}
and by $n_{x,y}$ in \cite{So1}, which is less well understood.
In particular $\widetilde P_{x,y}^{J}$ may be 0 even if $x<y$ in the
Bruhat ordering. Also, since the $\widetilde P_{x,y}^{J}$ are equal to
alternating sums of ordinary Kazhdan-Lusztig polynomials, it is
not a priori clear whether these polynomials have non-negative coefficients.
However, according to experts, it seems probable that
they admit a geometrical interpretation in terms of Schubert
varieties of finite codimension in an affine flag manifold.
This would settle the positivity conjecture VI.3 of \cite{LLT2}.
Note that in the case of two factors the polynomials
$c_{\nu^{(1)},\nu^{(2)}}^\lambda(q)$ are known to have 
non-negative coefficients because of their combinatorial
interpretation in terms of Yamanouchi domino tableaux given
in \cite{CL}.

That the non-vanishing 
of the polynomials $\widetilde P_{x,y}^{J}$ is a difficult
problem should not be too surprising. 
Indeed, our result shows that this  
contains as a special case the 
non-vanishing of the Littlewood-Richardson coefficients.
There has been some important recent progress by Klyachko
on this classical subject \cite{K} using toric vector bundles on the
projective plane (see the reviews of Zelevinsky \cite{Z}
and Fulton \cite{Fu}).
Maybe some new understanding will arise from 
the connection with affine Schubert varieties.

A few comments concerning the growing literature on $q$-analogues
of Littlewood-Richardson coefficients are in order.
In \cite{ShW} Shimozono and Weyman have studied the Poincar\'e 
polynomials of isotypic components of some virtual graded $GL_r$-modules 
supported in the closure of a nilpotent conjugacy class. 
These are $q$-analogues of Littlewood-Richardson multiplicities
$c_{\nu^{(1)},\ldots ,\nu^{(s)}}^\lambda$ satisfying a $q$-Kostant 
formula and a Morris-like recurrence.
In the case where all
partitions $\nu^{(j)}$ are rectangular (\ie the corresponding
weights are multiples of a single fundamental weight)
and are arranged in non-increasing order of width,
these polynomials have non-negative coefficients.
(This is not true in general, but see \cite{ShW}, Conjecture 4.)
In this case, a combinatorial interpretation 
in terms of semi-standard Young tableaux
was given by Shimozono \cite{Sh1,Sh2}, which shows that they coincide
with the generalized Kostka-Foulkes polynomials studied by
Schilling and Warnaar \cite{SchWa} in relation with exactly 
solvable lattice models and Rogers-Ramanujan type identities. 
A different combinatorial interpretation
using rigged configurations has been conjectured by Kirillov
and Shimozono \cite{KS} and recently verified \cite{KSS}. 
It is believed that for rectangular
shapes in non-increasing order these Poincar\'e polynomials 
are equal to the corresponding 
$c_{\nu^{(1)},\ldots ,\nu^{(s)}}^\lambda(q)$ but the reason
for that is still unclear. 

Let us describe more precisely the contents of this paper.
The results rely mainly on four sources, namely the parabolic 
analogue of Kazhdan-Lusztig polynomials developed by Deodhar
in \cite{De1,De2}, our joint paper with Lascoux on ribbon 
tableaux and generalizations of Kostka-Foulkes
polynomials \cite{LLT2}, our previous note \cite{LT}, and 
the paper of Varagnolo and Vasserot \cite{VV}.
Since \cite{LT} contains no proofs, and since only a small part of
\cite{LLT2} and \cite{VV} is needed to obtain our results, we thought it would
be appropriate to provide a self-contained exposition
of this material. 
Thus the style of the paper is openly expository and we hope
it can be read without a previous knowledge of these    
four sources.
However for what concerns parabolic Kazhdan-Lusztig polynomials,
we decided to omit the proofs because they can be found
in the optimum exposition by Soergel of Kazhdan-Lusztig theory from 
scratch \cite{So1}.

So in Section~\ref{SECT2} we explain all the necessary background
on (extended) affine symmetric groups $\EASG_r$ and their Hecke algebras 
$\EAH_r$.
In particular we introduce the two presentations (Coxeter-type
and Bernstein-type) and give the relations between them.
Following \cite{VV} we construct a representation of $\EAH_r$
on the weight lattice $\P_r$ of $\gl_r$ and introduce its two
Kazhdan-Lusztig bases.
The coefficients of these bases on the basis of weights
are the parabolic Kazhdan-Lusztig polynomials (for various 
parabolic subgroups).

In Section~\ref{SECT3}, we recall the Lusztig conjecture for 
quantum $\gl_r$ at an $n$th root of 1, the tensor product theorem,
and using a formula of Littlewood we deduce from this 
that the Littlewood-Richardson coefficients are value at 1 of
parabolic Kazhdan-Lusztig polynomials (Theorem~\ref{TH3}).

In Section~\ref{SECT4} we recall following \cite{LLT2} the
definitions of ribbon tableaux and their spin, we 
introduce the $q$-analogues $c_{\nu^{(1)},\ldots ,\nu^{(s)}}^\lambda(q)$,
and we state our main result (Theorem~\ref{TH4}).

In Section~\ref{SECT5} we explain the construction
of \cite{VV} and consider a quotient $\F_r$ of $\P_r$
whose bases are naturally labelled by dominant integral 
$\gl_r$-weights. 
This space can be identified in a natural way with the 
(finitized) $q$-deformed Fock space of Kashiwara, Miwa and
Stern \cite{KMS} considered in \cite{LT}.
Projecting on $\F_r$ the Kazhdan-Lusztig involution of $\P_r$
one gets the involution defined in \cite{LT} in terms
of $q$-wedge products. This implies that the canonical bases
of \cite{LT} have coefficients given by some parabolic
Kazhdan-Lusztig polynomials (Theorem~\ref{THVV}).

In Section~\ref{SECT6} we study the action of the center 
$Z(\EAH_r)$ of $\EAH_r$ on $\F_r$ and show that it can be expressed via
the combinatorics of ribbon tableaux. We then prove 
that the vectors $G^-_{\lambda + \rho}$ of the canonical basis  
indexed by non-restricted weights are obtained from the 
restricted ones by acting with an element of $Z(\EAH_r)$.
This should be regarded as an analogue in this setting of
the Steinberg-Lusztig tensor product theorem. 
Then we give the proof of Theorem~\ref{TH4}.

In Section~\ref{SECT7} we review the construction of
Kashiwara, Miwa and Stern of the Fock space $\FI$ 
obtained by taking the limit $r\rightarrow \infty$ in $\F_r$.
It affords a level 1 integrable representation of 
the quantum affine algebra $U_q(\slchap_n)$.
We investigate the behaviour of the canonical bases of 
$\FI$ with respect to the semi-linear involution
induced by the conjugation of partitions, and derive from
this a symmetry of the polynomials 
$c_{\nu^{(1)},\ldots ,\nu^{(s)}}^\lambda(q)$
(Theorem~\ref{THCONJ}) and 
an inversion formula for parabolic Kazhdan-Lusztig
polynomials (Corollary~\ref{CORINV}).
This formula, together with a result of Du, Parshall
and Scott \cite{DPS}, provides an alternative proof
of Soergel's character formula for tilting modules 
in type $A$ (Remark~\ref{REMTILT}).

Finally Section~\ref{SECT8} provides some numerical tables 
of $q$-Littlewood-Richardson
multiplicities and Kazhdan-Lusztig polynomials, which may serve as
examples of the results discussed in the text.


\section{Affine symmetric groups and their Hecke algebras} \label{SECT2}

\subsection{Affine symmetric groups} \label{SUBSECT21}

Let $\ASG_r$ denote the Coxeter group of type $\widetilde A_{r-1}$.
For $r=2$, this is the group generated by $s_0, s_1$ subject
to the relations $s_0^2=s_1^2=1$.
For $r>2$, $\ASG_r$ is generated by $s_0, s_1, \ldots , s_{r-1}$ subject
to 
\begin{eqnarray}
&&s_i s_{i+1} s_i=s_{i+1}s_is_{i+1},\label{EQ_s1}\\
&&s_is_j=s_js_i,\qquad\qquad\qquad\qquad (i-j \not = \pm1),\label{EQ_s2}\\
&&s_i^2=1,\label{EQ_s3}
\end{eqnarray}
where the subscripts are understood modulo $r$.
The subgroup generated by $s_1,\ldots ,s_{r-1}$ is isomorphic
to the symmetric group $\SG_r$.
The group $\ASG_r$ has a concrete realization as an affine reflection group.
Let $(\epsilon_1,\ldots ,\epsilon_r)$ denote the standard
basis of $\R^r$, and define a scalar product by putting
$(\epsilon_i,\epsilon_j)=\delta_{ij}$. 
Set $\alpha_i = \epsilon_i-\epsilon_{i+1}\ (1\le i\le r-1)$ and
$\alpha_0 = \epsilon_r - \epsilon_1$.
Let $\h_r$ denote a Cartan subalgebra of $\gl_r$.
We identify $\R^r$ with (the real part of) $\h^*_r$ in the
usual way, so that $P=P_r:=\bigoplus_{i=1}^r \Z \epsilon_i$
becomes the weight lattice, $Q=Q_r:=\bigoplus_{i=1}^{r-1} \Z \alpha_i$
the root lattice, $\alpha_i\ (1\le i \le r-1)$ the simple roots,
$-\alpha_0$ the highest root, {\it etc.} 
For $\alpha \in \R^r$ and $m \in \Z$, denote by $S_{\alpha,m}$
the affine reflection defined by
\[
S_{\alpha,m}(\lambda) = \lambda - 
2{(\lambda,\alpha)+m\over (\alpha,\alpha)} \, \alpha .
\]
Then for any $m$, the assignment 
\[
s_0 \mapsto S_{\alpha_0,m} \,, \quad
s_i \mapsto S_{\alpha_i,0} \quad (1\le i \le r-1)
\]
defines a faithful representation $\pi_m$ of $\ASG_r$ as a discrete subgroup
of the group of affine transformations of $\R^r$.
In coordinates, we have 
\begin{eqnarray*}
&&\pi_m(s_i)(\lambda) = 
(\lambda_1,\ldots ,\lambda_{i+1},\lambda_i , \ldots ,\lambda_r),
\qquad \qquad (1\le i \le r-1), \\
&&\pi_m(s_0)(\lambda) = 
(\lambda_r+m,\lambda_2,\ldots,\lambda_{r-1},\lambda_1-m).
\end{eqnarray*}
Note that for $s\in\SG_r$, $\pi_m(s)$ does not depend on $m$.
We shall therefore simplify the notation and write
$s\lambda$ in place of $\pi_m(s)(\lambda)$.

This realization shows that $\ASG_r$ contains a large commutative subgroup,
namely the image under $\pi_m^{-1}$ of the group of translations by the 
vectors of the lattice $mQ$.
Write $\tr(\lambda)$ for the translation by $\lambda \in \R^r$,
and let $t_i$ denote the element of $\ASG_r$ corresponding 
to $\tr(m\alpha_i)$ under $\pi_m$. Then one can check that
\begin{eqnarray*}
&&t_1 = (s_0 s_{r-1} s_{r-2} \cdots s_3 s_2)(s_3 s_4 \cdots s_{r-1} s_0 s_1), \\
&&t_2 = (s_1 s_0 s_{r-1} \cdots s_4 s_3)(s_4 s_5 \cdots s_0 s_1 s_2), \\
&&\vdots \qquad \qquad \qquad  \vdots \\
&&t_{r-1} = (s_{r-2} s_{r-3} s_{r-4} \cdots s_1 s_0)
(s_1 s_2 \cdots s_{r-3} s_{r-2} s_{r-1}), \\
&&t_0 = (s_{r-1} s_{r-2} s_{r-3} \cdots s_2 s_1)
       (s_1 s_2 \cdots s_{r-2} s_{r-1} s_0).
\end{eqnarray*}

It will be convenient to enlarge $\ASG_r$ by adding the translations
by vectors of the lattice $mP$.
Abstractly, this extended affine symmetric group that we shall
denote by $\EASG_r$ may be defined as the group generated by
$s_0, s_1, \ldots , s_{r-1},\tau$ subject to relations
(\ref{EQ_s1}), (\ref{EQ_s2}), (\ref{EQ_s3}) together with
\begin{equation}
\tau s_i = s_{i+1} \tau ,
\end{equation}
where again subscripts are understood modulo $r$.
It is clear that each $w\in\EASG_r$ can be written in a unique
way as 
\begin{equation}
w = \tau^k \sigma, \qquad (k \in \Z,\ \sigma \in \ASG_r).
\end{equation}
An alternative useful presentation is as follows. The group
$\EASG_r$ is generated by the elements $s_1, \ldots , s_{r-1}, y_1, \ldots ,y_r$
subject to relations (\ref{EQ_s1}), (\ref{EQ_s2}), (\ref{EQ_s3})
with all indices between $1$ and $r-1$ together with
\begin{eqnarray}
&&y_iy_j = y_jy_i, \\
&&s_iy_j = y_js_i \mbox{ for } j\not = i,i+1, \\
&&s_iy_is_i = y_{i+1}. \label{EQ8}
\end{eqnarray}
The homomorphism $\pi_m$ can then be extended to $\EASG_r$ by setting
\[
\pi_m(y_i) := \tr(m\epsilon_i), \qquad 
\pi_m(\tau) := 
S_{\alpha_1,0}S_{\alpha_2,0} \cdots S_{\alpha_{r-1},0}\,\tr(m\epsilon_r),
\]
or in coordinates
\begin{eqnarray*}
&&\pi_m(y_i)(\lambda) =
(\lambda_1,\ldots ,\lambda_i+m, \ldots ,\lambda_r),
\qquad \qquad (1\le i \le r), \\
&&\pi_m(\tau)(\lambda) =
(\lambda_r+m,\lambda_1,\ldots,\lambda_{r-2},\lambda_{r-1}).
\end{eqnarray*}
The following equations relate the two above presentations
of $\EASG_r$:
\begin{eqnarray}
&&y_i = s_{i-1}s_{i-2}\cdots s_1 s_0 s_{r-1}s_{r-2}\cdots s_{i+1}\tau, 
\qquad (1\le i \le r)\label{EQ9} \\
&&\tau = s_1 s_2 \cdots s_{r-1} y_r , \\
&&s_0 = s_{r-1}s_{r-2}\cdots s_2 s_1 s_2 \cdots s_{r-1} y_1^{-1}y_r.
\end{eqnarray}
(In (\ref{EQ9}) the subscripts are understood modulo $r$.)

\begin{figure}[t]
\begin{center}
\leavevmode
\epsfxsize =10cm
\epsffile{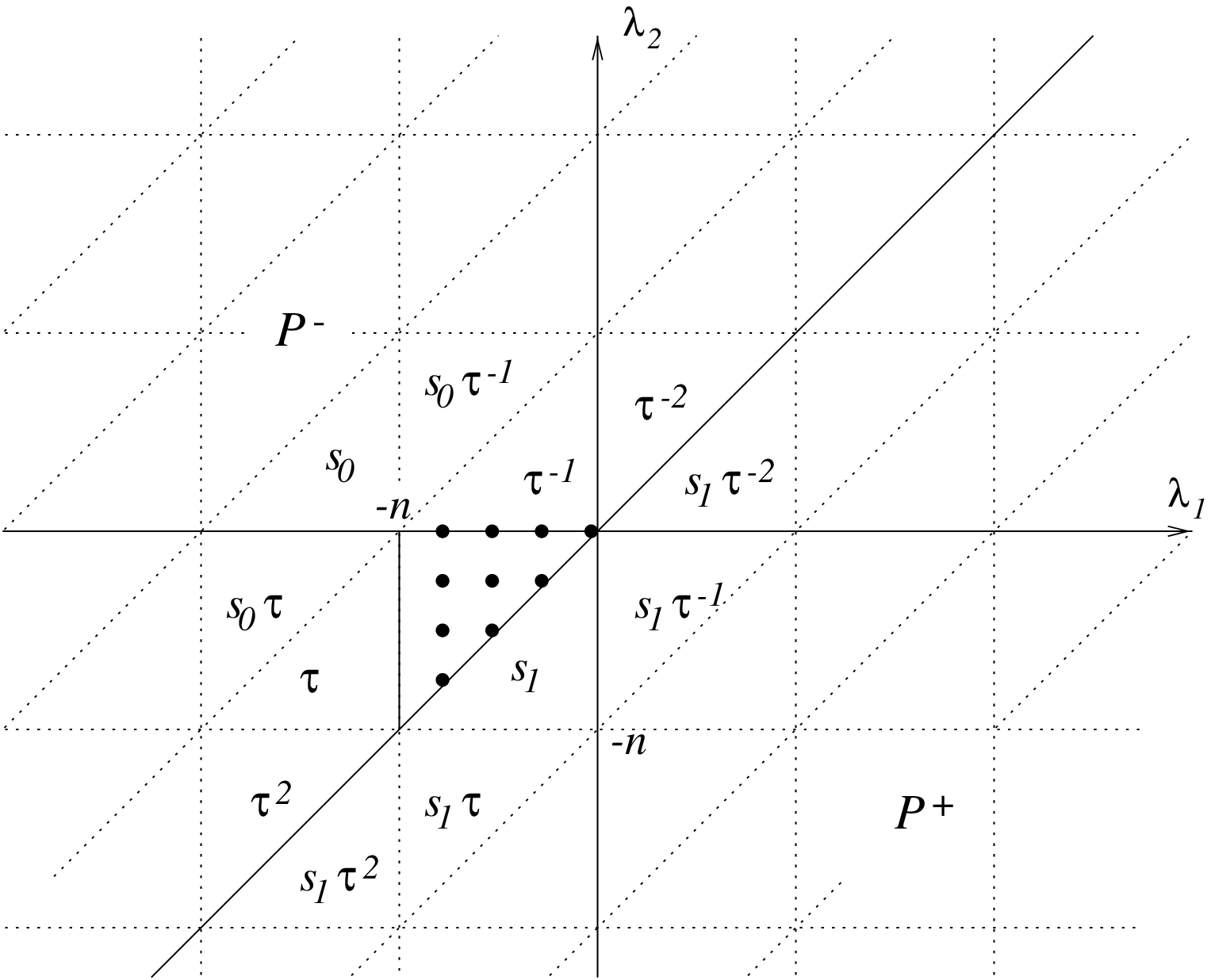}
\end{center}
\caption{\label{FIG0} The action of $\EASG_2$ on $P_2$ via $\pi_{-n}$}
\end{figure}

Note that $\EASG_r$ is not a Coxeter group. However, one can still
define a Bruhat order and a length function. 
Let $w = \tau^k \sigma, w' = \tau^m \sigma'$ with 
$k, m \in \Z,\ \sigma, \sigma' \in \ASG_r$. We say that 
$w<w'$ if and only if $k=m$ and $\sigma < \sigma'$, and we
put $\ell(w) := \ell(\sigma)$.
Define 
\[
A_{r,m} := \left|
\matrix{ 
\{\lambda \in \R^r \ \vert \
m > \lambda_1 \ge \lambda_2 \ge \cdots \ge\lambda_r\ge 0 \}
\mbox{\ \ if \ \ } m>0, \cr
\{\lambda \in \R^r \ \vert \ 
 m < \lambda_1 \le \lambda_2 \le \cdots \le\lambda_r \le 0\}
\mbox{\ \ if\ \ } m<0, 
}\right.
\]
and $\A_{r,m}:=A_{r,m} \cap P$ (see Figure~\ref{FIG0}).
Then $A_{r,m}$ is a fundamental domain for the action of $\EASG_r$
on $\R^r$ via $\pi_m$, that is, each orbit intersects it
in a unique point. Let $\lambda \in P$, and let 
$\nu$ be the intersection of $\A_{r,m}$ with the orbit of $\lambda$.
Then there is a unique $w(\lambda,m)\in \EASG_r$ of minimal length
such that $\pi_m(w(\lambda,m))(\nu) = \lambda$.
Let $\SG_{\nu,m}$ be the parabolic subgroup consisting of
the $w$ such that $\pi_m(w)(\nu)=\nu$. 
(Since $|\nu_1-\nu_r|<m$, $\SG_{\nu,m} \subset \SG_r$.)
Then $w(\lambda,m)$ is the minimal length representative of the
coset
\[
w(\lambda,m)\SG_{\nu,m} = \{w\in \EASG_r \ | \ \pi_m(w)(\nu)=\lambda \}.
\]
In this way, we can associate to the data $(\lambda,m)$ a certain
element $w(\lambda,m)$ of $\EASG_r$.
This will allow us to pass from the indexation by weights of the 
Littlewood-Richardson coefficients to the indexation by elements
of $\EASG_r$ of the Kazhdan-Lusztig polynomials.

\begin{example}{\rm
Take $r=3$ and $\lambda = (5,3,0)$. Then
\[
w(\lambda,3)=\tau^2 s_1, \qquad
w(\lambda,2)=\tau^3 s_0 s_1 s_2,\qquad
w(\lambda,-2)=\tau^{-5}s_1s_0s_2s_1s_2s_0.
\]\finex
}
\end{example}
For $\lambda=(\lambda_1,\ldots,\lambda_r) \in P$ set $\lambda_0:=\lambda_r+m,
\ \lambda_{r+1}:=\lambda_1-m$
and define the descent function
\[
\desc(\lambda,i,m) := \left\{
\matrix{
1& \mbox{ if } \lambda_i>\lambda_{i+1}, \cr
0& \mbox{ if } \lambda_i=\lambda_{i+1}, \cr
-1& \mbox{ if } \lambda_i<\lambda_{i+1},
}\right.
\qquad \qquad (0\le i \le r).
\]
Note that $\desc(\lambda,0,m) = \desc(\lambda,r,m)$ and
$\desc(\lambda,i,m) = \desc(\pi_m(\tau)(\lambda),i+1,m)$.
Geometrically, $\desc(\lambda,i,m) = 0$ means that $\lambda$ lies
on the reflecting hyperplane ${\cal H}_m$ of $\pi_m(s_i)$, \ie 
$\pi_m(s_i)(\lambda) = \lambda$, and $\desc(\lambda,i,m) = \sgn(m)$
means that $\lambda$ belongs to the 1/2-space defined by ${\cal H}_m$
which contains the fundamental domain $\A_{r,m}$, 
\ie $s_i w(\lambda,m) > w(\lambda,m)$.

\begin{lemma} \label{LEM2.1}
Let $\lambda \in P$, $i\in\{0,\ldots ,r-1\}$ and $m\in\Z^*$.
Let $\nu = w(\lambda,m)^{-1}(\lambda)$ be the point of $\A_{r,m}$
congruent to $\lambda$ under $\pi_m$.
One has the three following alternatives:
\[
\matrix{
\mbox{\rm (i) \ \ \ }\desc(\lambda,i,m) = \sgn(m)\hfill & \Longleftrightarrow &
s_iw(\lambda,m)= w(s_i\lambda,m)
> w(\lambda,m),\hfill \cr
\mbox{\rm (ii) \ \ }\desc(\lambda,i,m) = 0 \hfill &\Longleftrightarrow &
s_iw(\lambda,m)  = w(\lambda,m)s_j \mbox{ for some } 
s_j\in\SG_{\nu,m},  \,\hfill \cr
\mbox{\rm (iii) \ }\desc(\lambda,i,m) = -\sgn(m) \hfill &\Longleftrightarrow &
s_iw(\lambda,m)= w(s_i\lambda,m) < w(\lambda,m).\hfill
}
\]
\end{lemma}
\proof This is a reformulation of Lemma~2.1~(iii) of \cite{De1}.
Indeed, $\desc(\lambda,i,m) = -\sgn(m)$ if and only if 
$s_iw(\lambda,m) < w(\lambda,m)$ and in this case  
$s_iw(\lambda,m)= w(s_i\lambda,m)$ by \cite{De1}.
Also, $\desc(\lambda,i,m) = 0$ if and only if 
$s_iw(\lambda,m)\nu = w(\lambda,m)\nu$
which shows that $s_iw(\lambda,m)$ belongs to the same coset
as $w(\lambda,m)$ and is not minimal in this coset.
In this case, by \cite{De1}, there exists $s_j\in\SG_{\nu,m}$ such that
$s_iw(\lambda,m)  = w(\lambda,m)s_j$.
Finally, $\desc(\lambda,i,m) = \sgn(m)$ if and only if
$s_iw(\lambda,m) > w(\lambda,m)$ and 
$s_iw(\lambda,m)\nu \not = w(\lambda,m)\nu$.
In that case, again by \cite{De1}, $s_iw(\lambda,m)$ is minimal
in its coset and thus equal to $w(s_i\lambda,m)$.
\cqfd

\noindent
If $\nu$ is regular, that is, $\SG_{\nu,m} = \{1\}$, then case (ii)
does not occur and we obtain the following criterion
\begin{equation}
s_iw > w \quad \Longleftrightarrow \quad 
\desc(w\nu,i,m) = \sgn(m),\qquad (w\in\EASG_r).
\end{equation}
In particular, taking $m=r$ and $\nu=\rho:=(r-1,r-2,\ldots ,1,0)$
we get that
\begin{equation} \label{EQ13}
s_iw > w \quad \Longleftrightarrow \quad 
\desc(w\rho,i,r) = 1, \qquad (w\in\EASG_r).
\end{equation}

For $\lambda \in P$, set 
$y^\lambda := y_1^{\lambda_1}\cdots y_r^{\lambda_r}$.
Every $w\in\EASG_r$ has a unique decomposition
of the form $w=y^\lambda s$, where $\lambda \in P$
and $s\in\SG_r$. 
Therefore each coset $w\SG_r$ contains a unique 
element $y^\lambda$.
It follows from (\ref{EQ8}) that for $s\in \SG_r$,
$sy^\lambda = y^{s\lambda}s$.
This implies that each double coset $\SG_r w \SG_r$ in $\EASG_r$
contains a unique element $y^\lambda$ with 
$\lambda \in 
P^+:=\{\mu \in P \ |\ \mu_1\ge \mu_2 \ge \ldots \ge \mu_r\}$,
the set of dominant weights.
For $\lambda \in P^+$, we denote by $n_\lambda$
the element of maximal length in $\SG_r y^\lambda \SG_r$.
\begin{lemma} \label{LEM2.2}
Let $\lambda\in P^+$, $\mu\in P^-:= -P^+$ and  $s\in\SG_r$.
We have
\[
\ell(sy^\lambda)=\ell(s) + \ell(y^\lambda),\qquad
\ell(y^\mu s)=\ell(y^\mu)+\ell(s).
\]
In particular $n_\lambda = w_0 y^{\lambda}$, where
$w_0$ denotes the longest element of $\SG_r$.
\end{lemma}
\proof 
If $\lambda\in P^+$ then $\alpha:=y^\lambda\rho$ satisfies 
$\alpha_1>\alpha_2>\cdots >\alpha_r$. Let $s=s_{i_1}\cdots s_{i_k}$
be a reduced decomposition of $s$. By repeated applications
of (\ref{EQ13}) we see that $\ell(sy^\lambda) = \ell(y^\lambda)+k$,
which proves the first statement.
The case of $\mu$ is similar.
Finally, $w_0y^\lambda$ belongs to the double coset of $y^\lambda$
and for $s\in\SG_r$ $(s\not = 1)$, 
$\ell(sw_0y^\lambda)=\ell(sw_0)+\ell(y^\lambda) < \ell(w_0)+\ell(y^\lambda)$
so that $sw_0y^\lambda$ is not maximal. The argument is similar
for right multiplication by $s$, since $w_0y^\lambda = y^{w_0\lambda}w_0$,
and $w_0\lambda\in P^-$.
\cqfd

\begin{example}{\rm
Take $r=3$. Then
\[
n_{(2,1,0)}= s_2s_1s_2y_1^2y_2 =s_2s_1s_2s_0s_2s_1s_2\tau^3, \qquad
n_{(1,1,1)}= s_2s_1s_2y_1y_2y_3 =s_2s_1s_2\tau^3.
\]\finex
}
\end{example}
In fact, Lemma~\ref{LEM2.2} easily follows from a general formula of
Iwahori and Matsumoto (\cite{IM}, Prop. 1.23) which in our case reads
\begin{equation}
\ell(sy^\lambda)=
\sum_{\stackrel{\scriptstyle 
i < j }{s(i)<s(j)}} |\lambda_i -  \lambda_j|
+ 
\sum_{\stackrel{\scriptstyle i < j }{s(i)>s(j)}} |\lambda_i -  \lambda_j+1|,
\qquad (s\in\SG_r,\ \lambda\in P).
\end{equation}
In particular, if $\lambda \in P^+$ then
$
\ell(y^\lambda) = \sum_{i=1}^r (r+1-2i)\lambda_i,$
which shows that
\begin{equation}\label{EQ15}
\ell(y^\lambda)+\ell(y^\mu) = \ell(y^{\lambda+\mu}), \qquad
(\lambda,\mu \in P^+).
\end{equation}
\begin{lemma}\label{LEM1}
Let $\lambda \in P^+$ and set $\lambda^*:=w_0(-\lambda)$. 
Then, for all $n\ge r$ one has
$$
w(n\lambda+\rho,-n) = n_{\lambda^*}\,\tau^{-r+1}.
$$
\end{lemma}

\proof
Since $n\ge r$, the weight 
\[
\nu:=\pi_{-n}(\tau^{r-1} w_0)(\rho) = (1-n,2-n,\ldots ,r-1-n,0)
\]
belongs to $\A_{r,-n}$ and we have
\[
n\lambda + \rho = \pi_{-n}(y^{-\lambda})(\rho)
=\pi_{-n}(y^{-\lambda}w_0\tau^{-r+1})(\nu).
\]
The stabilizer of $\nu$ in $\pi_{-n}(\EASG_r)$ is trivial,
that is, $\nu$ is a regular weight.
Therefore we get 
\[
w(n\lambda + \rho,-n)=y^{-\lambda}w_0\tau^{-r+1}=w_0y^{w_0(-\lambda)}\tau^{-r+1}
=n_{\lambda^*}\,\tau^{-r+1}.
\]
\cqfd

\subsection{Affine Hecke algebras} \label{SUBSECT22}

The Hecke algebra $\EAH_r := H(\EASG_r)$ is the
algebra over $\Z[q,q^{-1}]$ with basis $T_w \ (w\in \EASG_r)$
and multiplication defined by 
\begin{eqnarray}
&&T_w T_{w'} = T_{w w'} \mbox{\ \ if\ \ } \ell(w w') = \ell(w) + \ell(w'), \\
&&(T_{s_i}-q^{-1})(T_{s_i}+q) = 0.
\end{eqnarray}
There is a canonical involution $x\mapsto \overline x$ of $\EAH_r$
defined as the unique ring homomorphism such that $\bar{q} = q^{-1}$
and $\bar{T_w} = (T_{w^{-1}})^{-1}$.

To simplify notation, we put $T_i:=T_{s_i}$ and we write $\tau$
instead of $T_\tau$. 
Then we have the two following presentations of $\EAH_r$
corresponding to the two above presentations of $\EASG_r$
(see \cite{Lu1,LuAMS}).
First, $\EAH_r$ is the algebra generated by 
$T_i \ (0\le i \le r-1)$ and an invertible element $\tau$ subject to 
the relations
\begin{eqnarray}
&&T_i T_{i+1} T_i=T_{i+1}T_iT_{i+1},\label{EQT1}\\
&&T_iT_j=T_jT_i,\qquad\qquad\qquad\qquad (i-j \not = \pm1),\label{EQT2}\\
&&(T_i-q^{-1})(T_i+q) = 0,\label{EQT3}\\
&&\tau T_i = T_{i+1} \tau. \label{EQT4}
\end{eqnarray}
Alternatively, $\EAH_r$ is the algebra generated by
$T_i \ (1\le i \le r-1)$ and invertible elements $Y_i \ (1\le i \le r)$ 
subject to the relations (\ref{EQT1}), (\ref{EQT2}), (\ref{EQT3})
with subscripts between 1 and $r-1$ together with
\begin{eqnarray}
&&Y_iY_j = Y_jY_i, \label{EQB1}\\
&&T_iY_j = Y_jT_i \mbox{ for } j\not = i,i+1, \label{EQB2} \\
&&T_iY_iT_i = Y_{i+1}. \label{EQB3}
\end{eqnarray}
The following equations relate the two above presentations
of $\EAH_r$:
\begin{eqnarray}
&&Y_i = T_{i-1}T_{i-2}\cdots T_1 
T_0^{-1} T_{r-1}^{-1}T_{r-2}^{-1}\cdots T_{i+1}^{-1}\tau,
\qquad (1\le i \le r)\label{EQT9} \\
&&\tau = T_1^{-1} T_2^{-1} \cdots T_{r-1}^{-1} Y_r , \\
&&T_0 = T_{r-1}^{-1}T_{r-2}^{-1}\cdots 
T_2^{-1} T_1^{-1} T_2^{-1} \cdots T_{r-1}^{-1} Y_1^{-1}Y_r.
\end{eqnarray}
(In (\ref{EQT1})(\ref{EQT2})(\ref{EQT4})(\ref{EQT9}) 
the subscripts are understood modulo $r$.)

Note that for $\lambda \in P$, we have two natural elements in 
$\EAH_r$ corresponding to the translation by $\lambda$, namely,
$Y^\lambda := Y_1^{\lambda_1}\cdots Y_r^{\lambda_r}$ and 
$T_\lambda:=T_{y^\lambda}$.
They do not coincide in general. 
(For example if $r=3$, $T_{y_2}=T_1T_0\tau$ and $Y_2=T_1T_0^{-1}\tau$.)
In fact the $T_\lambda$ do
not commute in general. However it follows from (\ref{EQ15}) that 
$T_\lambda T_\mu=T_{\lambda+\mu}=T_\mu T_\lambda$ if $\lambda,\mu \in P^+$.
Let $\lambda\in P$ be written as $\lambda = \lambda' - \lambda''$
with $\lambda', \lambda'' \in P^+$.
Bernstein has introduced an element $\widehat{T}_\lambda \in \EAH_r$ by
\[
\widehat{T}_\lambda := T_{\lambda'\rule{0mm}{2.7mm}} T_{\lambda''}^{-1}.
\]
This element is well-defined, \ie it does not depend on the choice
of $\lambda'$ and $\lambda''$, and 
$\widehat{T}_\lambda\widehat{T}_\mu=
\widehat{T}_{\lambda+\mu}=\widehat{T}_\mu\widehat{T}_\lambda$ 
for all $\lambda,\mu \in P$.
With this notation one can check that 
\begin{equation}\label{EQBERN}
Y^{\lambda}= \bar{\widehat{T}_\lambda} 
= T_{-\lambda'}^{-1}T_{-\lambda''\rule{0mm}{2.7mm}}
= T_{-\lambda''\rule{0mm}{2.7mm}} T_{-\lambda'}^{-1}.
\end{equation}
In particular, if $\lambda\in P^-$ then 
\begin{equation}\label{EQ24}
Y^\lambda = T_\lambda .
\end{equation}
%

\subsection{Action of $\EAH_r$ on the weight lattice} \label{SUBSECT23}

Let $\P = \P_r := \Z[q,q^{-1}]\otimes_\Z P$.  
We shall use the descent function to $q$-deform 
the representation $\pi_m$ of
$\EASG_r$ on $P$ into a representation $\Pi_m$ of $\EAH_r$ on $\P$.
Indeed, it follows from Lemma~\ref{LEM2.1} that 
$\EAH_r$ acts on $\P$ by $\Pi_m(\tau)(\lambda) := \pi_m(\tau)(\lambda)$ and
\[
\Pi_m(T_i)(\lambda) := \left\{
\matrix{
\pi_m(s_i)(\lambda)\hfill &\mbox{if \ \ } \desc(\lambda,i,m) = \sgn(m),\hfill \cr
q^{-1} \lambda \hfill&\mbox{if \ \  } 
\desc(\lambda,i,m) = 0,\hfill \cr
\pi_m(s_i)(\lambda) + (q^{-1}-q)\lambda  &\mbox{if \ \ }  
\desc(\lambda,i,m) = -\sgn(m),\hfill
}\right.
 \quad (0\le i \le r-1).
\]
{\bf Warning}\quad From now on in order to simplify the notation
we shall often omit the dependence on $m$ and write for example
$T_i\lambda$ in place of $\Pi_m(T_i)(\lambda)$, or
$s_i\lambda$ in place of $\pi_m(s_i)(\lambda)$.
We hope that this will not create confusion.

\medskip\noindent
In terms of the Kazhdan-Lusztig elements $C'_i:=T_i+q$ and $C_i:=T_i-q^{-1}$
we have
\[
C'_i\lambda = \left\{
\matrix{
s_i\lambda + q\lambda \hfill &\mbox{if \ \ } \desc(\lambda,i,m) = \sgn(m),
\hfill \cr
(q+q^{-1}) \lambda \hfill&\mbox{if \ \  }
\desc(\lambda,i,m) = 0,\hfill \cr
s_i\lambda + q^{-1}\lambda  &\mbox{if \ \ }
\desc(\lambda,i,m) = -\sgn(m),\hfill
}\right.
 \quad (0\le i \le r-1),
\]
\[
C_i\lambda = \left\{
\matrix{
s_i\lambda - q^{-1}\lambda \hfill &\mbox{if \ \ } \desc(\lambda,i,m) = \sgn(m),
\hfill \cr
0 \hfill&\mbox{if \ \  }
\desc(\lambda,i,m) = 0,\hfill \cr
s_i\lambda - q\lambda \hfill &\mbox{if \ \ }
\desc(\lambda,i,m) = -\sgn(m),\hfill
}\right.
 \quad (0\le i \le r-1).
\]
These formulas show that the $\EAH_r$-module $\P$ decomposes as 
\[
\P = \bigoplus_{\nu \in \A_{r,m}} \EAH_r \nu .
\]
Moreover, each summand of the right-hand side is isomorphic to 
an induced module. Indeed, let $\EAH_{\nu,m}$ be the 
subalgebra of $\EAH_r$ generated by the $T_i$ such that
$s_i\nu = \nu$, and let ${\bf 1}_{q^{-1}}$ denote
the 1-dimensional $\EAH_{\nu,m}$-module in which $T_i$ acts
by multiplication by $q^{-1}$. Then
\[
\EAH_r \nu \simeq \EAH_r \otimes_{\EAH_{\nu,m}} {\bf 1}_{q^{-1}} ,
\]
the isomorphism being given by
\begin{equation}\label{EQISO}
\lambda \in \EASG_r \nu \mapsto T_{w(\lambda,m)}\otimes 1.
\end{equation}
In particular $\lambda = T_{w(\lambda,m)}\nu$.
\subsection{Kazhdan-Lusztig polynomials}\label{SECT24}
	
The module $\P_\nu := \EAH_r \nu$ is a parabolic module of the
type considered by Deodhar in \cite{De1,De2}. 
(Note that if $\nu$ is a regular weight, then $\P_\nu$
is just the regular representation of $\EAH_r$.)
Therefore $\P_\nu$ has two Kazhdan-Lusztig bases constructed as follows
(see \cite{So1}).
Define a semi-linear involution on $\P_\nu$ by
\[
\bar{q} := q^{-1}, \qquad 
\bar{x\nu} := \bar{x}\nu \quad (x\in\EAH_r) ,
\]
and two lattices 
\[
L^+_\nu := \bigoplus_{\lambda \in \EASG_r\nu} \Z[q]\lambda ,
\qquad
L^-_\nu := \bigoplus_{\lambda \in \EASG_r\nu} \Z[q^{-1}]\lambda.
\]
Then there are two bases $C^+_\lambda,\ C^-_\lambda \ (\lambda \in \EASG_r\nu)$
characterized by
\[
\bar{C^+_\lambda}=C^+_\lambda ,\qquad \bar{C^-_\lambda}=C^-_\lambda,
\]
and 
\[
C^+_\lambda \equiv \lambda \mod qL^+_\nu , \qquad
C^-_\lambda \equiv \lambda \mod q^{-1}L^-_\nu . 
\]
When $\nu$ is regular these bases coincide with the Kazhdan-Lusztig
bases $C_w'$ and $C_w$ respectively under the isomorphism (\ref{EQISO}).

These bases can be computed recursively as follows \cite{So1}.
First, by definition, $C^+_\nu=C^-_\nu = \nu$, and more generally
$C^+_{\tau^k\nu}=C^-_{\tau^k\nu} = \tau^k\nu \ (k\in\Z)$.
Let $\lambda \in \EASG_r\nu$ and suppose that $C^+_\mu$
(\resp $C^-_\mu$) has already been calculated for all $\mu<\lambda$,
that is, such that $w(\mu,m)<w(\lambda,m)$.
Then compute $v^+_\lambda = C'_i\,C^+_\mu$
(\resp $v^-_\lambda = C_i\,C^-_\mu$) where $\mu$ and $i$ satisfy
$s_i(\mu) = \lambda$ and $\desc(\mu,i,m)=\sgn(m)$.
Then $v^+_\lambda$ (\resp $v^-_\lambda$) is invariant under the
bar-involution and belongs to $L^+_\nu$ (\resp $L^-_\nu$).
Write 
\[
v^+_\lambda \equiv \lambda + \sum_\alpha a_\alpha \alpha \mod qL^+_\nu,
\qquad
(\mbox{\resp } 
v^-_\lambda \equiv \lambda + \sum_\beta b_\beta \beta \mod q^{-1}L^-_\nu),
\]
where $a_\alpha,\ b_\beta \in \Z$.
The weights $\alpha$, $\beta$ occuring in the right-hand side are
certainly $<\lambda$ and we obtain
\[
C^+_\lambda = v^+_\lambda  - \sum_\alpha a_\alpha C^+_\alpha
\qquad
(\mbox{\resp } 
C^-_\lambda = v^-_\lambda - \sum_\beta b_\beta C^-_\beta)
.
\]
\begin{example}{\rm
Let us take $r=3$, $m=-2$ and compute $C^-_{(0,6,1)}$.
We have $w((0,6,1),-2)=s_2s_0s_1s_2s_0\tau^{-4}$ and
$\nu := w((0,6,1),-2)^{-1} (0,6,1) = (-1,0,0)$.
Clearly 
\[
C^-_{(2,2,3)}=C^-_{\tau^{-4}(-1,0,0)}=(2,2,3).
\]
Then we compute successively ($t=q^{-1}$)
\begin{eqnarray*}
&&v^-_{(1,2,4)}=C^-_{(1,2,4)}= (1, 2, 4)-t(2, 2, 3) , \\
&&v^-_{(1,4,2)}=C^-_{(1,4,2)}= 
(1, 4, 2)-t(1, 2, 4)-t(2, 3, 2)+t^2(2, 2, 3) ,\\
&&v^-_{(4,1,2)}=C^-_{(4,1,2)}= 
(4, 1, 2)-t(1, 4, 2)-t(2, 1, 4)+t^2(1, 2, 4)-t(3, 2, 2)+t^2(2, 3 , 2),\\
&&v^-_{(0,1,6)}=C^-_{(0,1,6)}=
(0, 1, 6)-t(4, 1, 2)-t(0, 4, 3)+t^2(1, 4, 2)+t^2(2, 2, 3)-t(1, 2 , 4)\\
&&\hskip 3.5cm -t(0, 2, 5)+t^2(3, 2, 2)+t^2(0, 3, 4)-t^3(2, 3, 2),\\
&&v^-_{(0,6,1)}=(0, 6, 1)-t(0, 1, 6)-t(4, 2, 1)+t^2(4, 1, 2)-t(0, 3, 4)\\
&&\hskip 3.5cm
+(0, 4, 3)+2t^2(1, 2, 4)-2t(1, 4, 2)+2t^2(2, 3, 2)-2t^3(2, 2, 3)\\
&&\hskip 3.5cm -t(0, 5, 2)+t^2(0 , 2, 5)+t^2(0, 4, 3)-t^3(0, 3, 4).
\end{eqnarray*}	
We see that $v^-_{(0,6,1)}\equiv (0, 6, 1)+(0, 4, 3) \mod tL^-_\nu$.
Thus subtracting the previously calculated element
\[
C^-_{(0,4,3)}=(0, 4, 3)-t(0, 3, 4)-t(1, 4, 2)+t^2(1, 2, 4)+t^2(2, 3, 2)-
t^3(2, 2, 3)
\]
we get
\begin{eqnarray*}
&&
C^-_{(0,6,1)}=(0, 6, 1)-t(0, 1, 6)-t(4, 2, 1)+t^2(4, 1, 2)+t^2(1, 2, 4)
-t(1, 4 , 2)\\
&&\hskip 1.5cm +t^2(2, 3, 2)-t^3(2, 2, 3)-t(0, 5, 2)+t^2(0, 2, 5)+t^2(0, 4, 3)
-t^3(0 , 3, 4).
\end{eqnarray*}
\finex}
\end{example}
Put
\[
C'_w=\sum_{x\in\EASG_r} P_{x,w}(q)\, T_x .
\]
Then (see \cite{So1}, Theorem 2.7)
\[
C_w=\sum_{x\in\EASG_r} P_{x,w}(-q^{-1})\, T_x .
\]
The $P_{x,w}$ are the Kazhdan-Lusztig polynomials (up to 
a factor $q^{\ell(w)-\ell(x)}$ and the change of variable $q\mapsto q^{-2}$).
Similarly for $\lambda\in\EASG_r\nu$ write 
\[
C^+_\lambda=\sum_{\mu \in  \EASG_r\nu} P^+_{\mu,\lambda}(q)\, \mu,
\qquad
C^-_\lambda=\sum_{\mu \in  \EASG_r\nu} P^-_{\mu,\lambda}(-q^{-1})\, \mu .
\]
Then $P^+_{\mu,\lambda}$ and $P^-_{\mu,\lambda}$ are Deodhar's
polynomials $P^J_{w(\mu,m),w(\lambda,m)}$ and 
$\widetilde P^J_{w(\mu,m),w(\lambda,m)}$ respectively 
(again up to a power of $q$ and the change of 
variable $q\mapsto q^{-2}$), where
$J$ is the set of indices $i$ of the Coxeter generators $s_i \in \SG_{\nu,m}$.
Their expression in terms of ordinary Kazhdan-Lusztig polynomials
is given by 

\begin{theorem}[Deodhar \cite{De1,De2}]\label{THDEO}
Let $w_{0,\nu}$ be the longest element of $\SG_{\nu,m}$.
Then
\[
P^+_{\mu,\lambda} = P_{w(\mu,m)w_{0,\nu}\,,\,w(\lambda,m)w_{0,\nu}}, \qquad
P^-_{\mu,\lambda} = \sum_{z\in\SG_{\nu,m}}
(-q)^{\ell(z)}  P_{w(\mu,m)z\,,\,w(\lambda,m)}.
\]
\end{theorem}
We shall also need the following simple observation (see \cite{So1}, Remark 3.2.4).
Suppose that $\desc(\lambda,i,m) = \desc(\mu,i,m) = -\sgn(m)$. Then
\begin{equation}\label{OBSERV}
P^+_{s_i\mu,\lambda} = qP^+_{\mu,\lambda}, \qquad
P^-_{s_i\mu,\lambda} = qP^-_{\mu,\lambda}.
\end{equation}
This follows from the fact that if $\desc(\lambda,i,m)=-\sgn(m)$ then
\[
C'_iC^+_\lambda = (q+q^{-1})\,C^+_\lambda , \qquad
C_iC^-_\lambda = -(q+q^{-1})\,C^-_\lambda.
\]

\section{Littlewood-Richardson coefficients and Kazhdan-Lusztig polynomials}
\label{SECT3}

\subsection{The Lusztig conjecture}

Let $U_q(\gl_r)$ be the quantum enveloping algebra of $\gl_r$.
This is a $\Q(q)$-algebra with generators $E_i,\ F_i, \ q^{\pm\epsilon_j}\ 
(1\le i \le r-1,\ 1\le j \le r)$. The relations are standard \cite{Jim}
and will be omitted. To avoid confusion when $q$ is specialized 
to a complex number, we shall write $K_j^{\pm}$ in place of $q^{\pm\epsilon_j}$.
Let $U_{q,\Z}(\gl_r)$ denote the $\Z[q,q^{-1}]$-subalgebra  
generated by the elements 
\[
E_i^{(k)}:={E_i^k\over [k]!},\quad 
F_i^{(k)}:={F_i^k\over [k]!},\quad  
K_j^{\pm}, \qquad (k\in\N),
\]
where $[k]!:=[k][k-1]\cdots [2][1]$ and $[k]:=(q^k-q^{-k})/(q-q^{-1})$.
Let $\zeta \in \C$ be such that $\zeta^2$ is a primitive $n$th
root of 1. One defines 
$U_{\zeta}(\gl_r) := U_{q,\Z}(\gl_r)\otimes_{\Z[q,q^{-1}]} \C$
where $\Z[q,q^{-1}]$ acts on $\C$ by $q\mapsto \zeta$ \cite{Lu4,Lu5}.

Let $\lambda \in P^+_r$. There is a unique finite-dimensional
$U_q(\gl_r)$-module (of type 1) $W_q(\lambda)$ with highest weight $\lambda$.
Its character is the same as in the classical case and is 
given by Weyl's character formula
\begin{equation}
\ch W_q(\lambda) = s_\lambda(e^{\epsilon_1},\ldots , e^{\epsilon_r}),
\end{equation}
where $s_\lambda$ denotes the Schur function (see \cite{Mcd}).
Fix a highest weight vector $u_\lambda \in W_q(\lambda)$ and 
denote by $W_{q,\Z}(\lambda)$ the $U_{q,\Z}(\gl_r)$-submodule
of $W_q(\lambda)$ generated by acting on $u_\lambda$.
Finally, put 
\[
W_\zeta(\lambda):= W_{q,\Z}(\lambda)\otimes_{\Z[q,q^{-1}]} \C.
\]
This is a $U_{\zeta}(\gl_r)$-module called a Weyl module \cite{Lu4}.
By definition $\ch W_\zeta(\lambda) =\ch W_q(\lambda) $.

There is a unique simple quotient of $W_\zeta(\lambda)$ denoted by $L(\lambda)$.
Its character is given in terms of the characters of the Weyl modules
by the so-called Lusztig conjecture.
Put $m=-n$ (this assumption will be in force for the whole Section~\ref{SECT3})
and consider the action of $\EASG_r$ on $P$ via $\pi_m$.
For $\lambda \in P^+$ write $\nu := w(\lambda+\rho,m)^{-1}(\lambda+\rho)$. Then
\begin{theorem}[Kazhdan-Lusztig, Kashiwara-Tanisaki]
\begin{equation}\label{LC}
\ch L(\lambda)  = \sum_w
(-1)^{\ell(w(\lambda+\rho,m))-\ell(w)}\,
P_{w,w(\lambda+\rho,m)}(1)\,\, \ch W_\zeta(w(\nu)-\rho) ,
\end{equation}
where the sum runs over the $w\in\EASG_r$ such that $w<w(\lambda+\rho,m)$ and
$w(\nu)-\rho\in P^+$.
\end{theorem}
Note that if $\lambda + \rho$ is a singular weight 
(\ie its stabilizer is non trivial) the coefficient
of a given Weyl module $W_\zeta(\mu)$ in the right-hand side of (\ref{LC})
is an alternating sum of $P_{w,w(\lambda+\rho,m)}(1)$ over the stabilizer 
$\SG_{\nu,m}$ of $\nu$. 
In fact, using the notation of Section~\ref{SECT24} one can rewrite
Eq. (\ref{LC}) as
\begin{equation}\label{EQ28}
\ch L(\lambda)  = \sum_\mu
P^-_{\mu+\rho,\lambda+\rho}(-1)\,\, \ch W_\zeta(\mu) ,
\end{equation}
where the sum is over the $\mu\in P^+$ such that 
$\mu + \rho \in \EASG_r\nu$.

\begin{example}{\rm
Take $r=3$, $n=2$ and $\lambda = (4,0,0)$. Then $\lambda+\rho = (6,1,0)$
and for $m=-n=-2$, one has 
\begin{eqnarray*}
&&C^-_{(6,1,0)} = (6,1,0) - q^{-1} (6,0,1) -q^{-1} (1,6,0)
+q^{-2} (0,6,1) +q^{-2} (1,0,6) -q^{-3} (0,1,6)\\
&&\quad-q^{-1}(5,2,0) + q^{-2} (5,0,2) +q^{-2} (2,5,0)
-q^{-3} (0,5,2) -q^{-3} (2,0,5) +q^{-4} (0,2,5)\\
&&\quad+q^{-2}(4,3,0) - q^{-3} (4,0,3) -q^{-3} (3,4,0)
+q^{-4} (0,4,3) +q^{-4} (3,0,4) -q^{-5} (0,3,4).
\end{eqnarray*}
It follows that the character of $L(4,0,0)$ for $\zeta^2 = -1$
is given by
\[
\ch L(4,0,0) =\ch W_\zeta(4,0,0) -\ch W_\zeta(3,1,0) +\ch W_\zeta(2,2,0) .
\]\finex
}
\end{example}

\subsection{The tensor product theorem}

Let $\Fr$ denote the Frobenius map from $U_\zeta(\gl_r)$
to the (classical) enveloping algebra $U(\gl_r)$ \cite{Lu4,CP}. 
This is the algebra homomorphism defined by  
$\Fr(K_j)=1$ and 
\[
\Fr(E_i^{(k)})=\left\{
\matrix{
E_i^{(k/n)}
\mbox{\ \ if $n$ divides $k$,} \cr
0\hfill
\mbox{\ \ otherwise},
}\right. 
\qquad
\Fr(F_i^{(k)})=\left\{
\matrix{
F_i^{(k/n)}
\mbox{\ \ if $n$ divides $k$,} \cr
0\hfill
\mbox{\ \ otherwise}.
}\right. 
\]
(Here we slightly abuse notation and denote by the same symbols the
Chevalley generators of $U_\zeta(\gl_r)$ and those of $U(\gl_r)$.)
Given a $U(\gl_r)$-module $M$, one can thus define a $U_\zeta(\gl_r)$-module
$M^\Fr$ by composing the action of $U(\gl_r)$ with $\Fr$.
If $M$ is a finite-dimensional module with character the symmetric Laurent
polynomial $\ch M = \varphi(e^{\epsilon_1},\ldots ,e^{\epsilon_r})$,
then 
\[
\ch M^\Fr = p_n(\varphi)(e^{\epsilon_1},\ldots ,e^{\epsilon_r})
:=  \varphi(e^{n\epsilon_1},\ldots ,e^{n\epsilon_r}),
\]
the so-called plethysm of $\varphi$ with the power sum $p_n$ (see \cite{Mcd}).
In particular, the character of the pullback $W(\lambda)^\Fr$
of the classical Weyl module $W(\lambda)$ is the plethysm 
$p_n(s_\lambda)$.
\begin{theorem}[Lusztig \cite{Lu4}]\label{TH2}
Let $\lambda \in P^+$. Write $\lambda = \lambda^{(0)} + n\lambda^{(1)}$, where
$\lambda^{(0)}$ is $n$-restricted, that is,
\[
0\le \lambda_i^{(0)}-\lambda_{i+1}^{(0)} <n \qquad (1\le i \le r-1).
\]
The simple $U_\zeta(\gl_r)$-module $L(\lambda)$ is isomorphic to
the tensor product
\[
L(\lambda) \simeq L(\lambda^{(0)})\otimes W(\lambda^{(1)})^\Fr.
\]
\end{theorem}
Consider now the particular case when $\lambda$ is a partition
whose parts are all divisible by $n$. 
Then, writing $n\lambda$ in place of $\lambda$, we deduce from
Theorem~\ref{TH2} and Eq.~(\ref{EQ28}) that
\begin{equation}\label{EQ29}
p_n(s_\lambda) = \ch L(n\lambda) = 
\sum_\mu P^-_{\mu+\rho,n\lambda+\rho}(-1)\,\, \ch W_\zeta(\mu) 
=\sum_\mu P^-_{\mu+\rho,n\lambda+\rho}(-1)\,\, s_\mu ,
\end{equation}
where the sum is over the $\mu\in P^+$ such that
$\mu + \rho \in \EASG_r(n\lambda + \rho) =  \EASG_r\rho$.

\subsection{Expression of the Littlewood-Richardson coefficients}
\label{SECT33}

Let $\lambda \in \PP^+_r=\{\lambda\in P\ |
\lambda_1\ge \lambda_2\ge \ldots \ge \lambda_r\ge 0\}$,
the set of partitions of length $l(\lambda)\le r$.
It is a well-known result of Littlewood \cite{Lit}
that the coefficients in the expansion of $p_n(s_\lambda)$
on the basis of Schur functions are Littlewood-Richardson
multiplicities.
More precisely, if $\mu\in \PP^+_r$ is such that $\mu + \rho \in \EASG_r\rho$
then there is a unique expression 
\[
\mu+\rho = \gamma + n\alpha, \qquad (\gamma = s\rho,\ s\in \SG_r,\ \alpha \in \N^r)
\]
such that $i<j$ and $\gamma_i \equiv \gamma_j \ (n)$ implies
$\gamma_i > \gamma_j$.
Then for $k \in \{0,1,\ldots ,n-1\}$ the subsequence of $\alpha$
consisting of the $\alpha_i$ such that $\gamma_i \equiv k-r$
is a partition $\mu^{(k)}$ 
(possibly empty), and one has \cite{Lit} 
\begin{equation}\label{EQ30}
\<p_n(s_\lambda)\ ;\,s_\mu\> = (-1)^{\ell(s)} 
\<s_\lambda\ ; \  s_{\mu^{(0)}}\cdots s_{\mu^{(n-1)}}\> 
\end{equation}
where $\< \cdot \,;\, \cdot \>$ denotes the standard scalar product
of the algebra of symmetric functions for which the $s_\lambda$
form an orthonormal basis.
The $n$-tuple of partitions $(\mu^{(0)},\ldots ,\mu^{(n-1)})$
is called the $n$-quotient of $\mu$ and $(-1)^{\ell(s)}$ the
$n$-sign of $\mu$, denoted $\varepsilon_n(\mu)$. 
Conversely, provided that $r$ is large enough,
given an arbitrary $n$-tuple of partitions
$(\mu^{(0)},\ldots ,\mu^{(n-1)})$ there
exists a unique $\mu\in \PP_r^+$ such that $\mu+\rho\in \EASG_r \rho$
and $\mu$ has $(\mu^{(0)},\ldots ,\mu^{(n-1)})$ as $n$-quotient
(see \cite{Mcd,JK}).
\begin{example}{\rm
Let $r=8$, $n=3$, and $\mu =(6,6,4,4,4,3,2,1)$.
Then
\[
\mu + \rho = (13,12,9,8,7,5,3,1)
= (7,6,3,5,4,2,0,1) + 3\,(2,2,2,1,1,1,1,0).
\]
Thus the $3$-quotient of $\mu$ is
\[
(\mu^{(0)},\mu^{(1)},\mu^{(2)})=((1,1),(2,2,1),(2,1)).
\]\finex
}
\end{example}
Let us define the Littlewood-Richardson coefficient
\[
c_{\mu^{(0)},\ldots ,\mu^{(n-1)}}^\lambda := 
\<s_{\mu^{(0)}}\cdots s_{\mu^{(n-1)}}\ ; \  s_\lambda\>
=[W(\mu^{(0)})\otimes \cdots \otimes W(\mu^{(n-1)}) \ : \ 
W(\lambda) ]\,.
\]
Combining (\ref{EQ29}) and (\ref{EQ30}), we have obtained
\begin{theorem}\label{TH3}
Let $\lambda,\mu^{(0)},\ldots ,\mu^{(n-1)}$ be partitions and denote by
$\mu$ the partition with $n$-quotient $(\mu^{(0)},\ldots ,\mu^{(n-1)})$. 
Take $r\ge l(\mu)$, the number of parts of $\mu$.
Then, 
\[
c_{\mu^{(0)},\ldots ,\mu^{(n-1)}}^\lambda = P^-_{\mu+\rho,n\lambda+\rho}(1)
\]
where the right-hand side is a Kazhdan-Lusztig polynomial of parabolic type
for $\EASG_r$ with $m=-n$. 
In other words, setting $\nu = w(n\lambda + \rho,-n)^{-1}(n\lambda+\rho)$,
one has in terms of the (ordinary) Kazhdan-Lusztig polynomials
for $\EASG_r$
\[
c_{\mu^{(0)},\ldots ,\mu^{(n-1)}}^\lambda
=
\sum_{z\in\SG_{\nu,-n}}
(-1)^{\ell(z)}  P_{w(\mu+\rho,-n)z\,,\,w(n\lambda+\rho,-n)}(1).
\]
\end{theorem}
If $l(\lambda)>r$ the polynomial $P^-_{\mu+\rho,n\lambda+\rho}$
is not defined, but in this case $l(\lambda)>l(\mu)$ and it is easy
to see that $c_{\mu^{(0)},\ldots ,\mu^{(n-1)}}^\lambda = 0$.

Note that if $w = \tau^k \sigma,\ w' = \tau^{m} \sigma'$ with
$k, m \in \Z,\ \sigma, \sigma' \in \ASG_r$, then $P_{w,w'}$
is nonzero only if $k=m$ and then
$P_{w,w'}=P_{\sigma,\sigma'}$.
Thus the Kazhdan-Lusztig polynomials above are in fact polynomials
for $\ASG_r$.
\begin{example}{\rm
Take $r=3$ and $n=-m=2$. The dominant weights occuring in
the expansion of $C^-_{(6,3,0)}$ are 
\[
(6,3,0),\ (6,2,1),\ (5,4,0),\ (4,3,2),
\]
with respective coefficients 
\[
1,\ -q^{-1},\ -q^{-1},\ q^{-2}.
\]
This gives the following expressions for some Littlewood-Richardson
coefficients (which are all equal to 1):
\[
c_{(1),(2)}^{(2,1)}=P^-_{(6,3,0),(6,3,0)}(1),\qquad
c_{\emptyset,(2,1)}^{(2,1)}=P^-_{(6,2,1),(6,3,0)}(1),\
\]
\[
c_{(2),(1)}^{(2,1)}=P^-_{(5,4,0),(6,3,0)}(1),\qquad
c_{(1),(1,1)}^{(2,1)}=P^-_{(4,3,2),(6,3,0)}(1).\ 
\]
In terms of ordinary Kazhdan-Lusztig polynomials for $\ASG_3$
we can write for example 
\[
c_{(1),(1,1)}^{(2,1)} = P_{s_2s_0s_2\,,\,s_2s_0s_2s_1s_2s_0s_2}(1)
-P_{s_2s_0s_2s_1\,,\,s_2s_0s_2s_1s_2s_0s_2}(1) = 2 - 1.
\]\finex
}
\end{example}
\begin{example}{\rm
Let us express the coefficient $c_{(2,1),(2,1)}^{(3,2,1)} = 2$
in terms of Kazhdan-Lusztig polynomials.
We take $r=4$, $\lambda = (3,2,1)$ and $\mu = (4,4,2,2)$ so that 
$\mu$ has $2$-quotient $((2,1);(2,1))$.
It follows that 
\[
c_{(2,1),(2,1)}^{(3,2,1)} = P^-_{(7,6,3,2),(9,6,3,0)}(1).
\]
This Kazhdan-Lusztig polynomial corresponds to the following 
elements of $\EASG_4$:
\begin{eqnarray*}
&&w((9,6,3,0),-2) = s_1s_2s_1s_3s_2s_1s_0s_1s_3s_2s_1s_3s_0s_1s_3s_2s_0\tau^{-10},\\
&&w((7,6,3,2),-2) = s_1s_2s_1s_3s_2s_1s_0s_1s_3s_2s_0\tau^{-10}.\hfill
\end{eqnarray*}\finex
}
\end{example}
Observe that if $n\ge r$,
then the $n$-quotient of the partition $n\mu=(n\mu_1,\ldots ,n\mu_r)$
is just $((\mu_1),\ldots ,(\mu_r),\emptyset,\ldots ,\emptyset)$
up to reordering, and therefore Theorem~\ref{TH3} gives 
\[
P^-_{n\mu+\rho,n\lambda+\rho}(1) = 
c_{(\mu_1),\ldots ,(\mu_r)}^\lambda = K_{\lambda,\mu},
\]
the Kostka number. On the other hand, taking into account Lemma~\ref{LEM1} and
the fact that the weight $n\lambda+\rho$ is regular, one also has 
\[
P^-_{n\mu+\rho,n\lambda+\rho}(1) 
= P_{n_{\mu^*}\,\tau^{-r+1},n_{\lambda^*}\,\tau^{-r+1}}(1)
= P_{n_{\mu^*},n_{\lambda^*}}(1).
\]
Hence
\[
P_{n_{\mu},n_{\lambda}}(1)=K_{\lambda^*,\mu^*}=K_{\lambda,\mu}
\]
since the weight multiplicities of the contragredient representation
$W(\lambda^*)$ are equal to those of $W(\lambda)$,
and we recover the expression of \cite{Lu2} for the weight
multiplicities.

Thus we see that the modular Lusztig conjecture with its restriction
$n\ge r$ is enough to express the weight multiplicities in terms
of Kazhdan-Lusztig polynomials, but for what concerns the 
general tensor product multiplicities we need the case $n<r$ and
the quantum Lusztig conjecture.


\section{Littlewood-Richardson coefficients and ribbon tableaux} 
\label{SECT4}

\subsection{Ribbon tableaux}

Let us start from the well-known formula
\begin{equation}\label{EQ31}
h_\mu = \sum_\lambda |\tab(\lambda,\mu)| \, s_\lambda,
\end{equation}
where $h_\mu:=h_{\mu_1}\ldots h_{\mu_r}$ is a product of complete
homogeneous symmetric functions and $\tab(\lambda,\mu)$
denotes the set of semi-standard Young tableaux 
of shape $\lambda$ and weight $\mu$ \cite{Mcd}.
Let $n\in\N^*$.
Semi-standard $n$-ribbon tableaux are combinatorial objects
which replace ordinary Young tableaux when one substitutes
the plethysm $p_n(h_\mu)$ in place of $h_\mu$ in (\ref{EQ31}).
More precisely, denoting by $\tab_n(\lambda,\mu)$ the set of
$n$-ribbon tableaux of shape $\lambda$ and weight $\mu$
(to be defined below), one has
\begin{equation}\label{EQ32}
p_n(h_\mu) = \sum_\lambda \varepsilon_n(\lambda)\, |\tab_n(\lambda,\mu)| 
\, s_\lambda,
\end{equation}
where $\varepsilon_n(\lambda)$ is the $n$-sign of $\lambda$.

A ribbon tableau of weight $\mu=(1,1,\ldots ,1)$ is called
standard. Standard ribbon tableaux were introduced by 
Stanton and White \cite{SW} in relation with generalizations
of the Robinson-Schensted correspondence for the complex
reflection groups $G(n,1,r)=(\Z/n\Z)\wr \SG_r$.
In particular, the case $n=2$ (domino tableaux) is  related 
to Weyl groups of type B, C, D, and therefore to the geometry
of flag manifolds for classical groups \cite{vL} and to the
classification of the primitive ideals of classical enveloping
algebras \cite{BV,Ga}.
Semi-standard domino tableaux were introduced in \cite{CL}
for calculating the multiplicities of the symmetric and 
alternating square of an irreducible representation of $\gl_r$
(see also \cite{KLLT,vL2}).
In an attempt to extend the results of \cite{CL} to higher degree
plethysms, semi-standard $n$-ribbon tableaux were defined in
\cite{LLT2} and several conjectures were formulated.
We shall give a brief review of \cite{LLT2} refering to the paper
for more detail.

\begin{figure}[t]
\begin{center}
\leavevmode
\epsfxsize =4cm
\epsffile{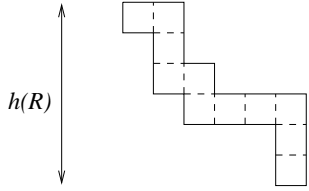}
\end{center}
\caption{\label{FIG1} An 11-ribbon of height $h(R)=6$}
\end{figure}
\begin{figure}[t]
\begin{center}
\leavevmode
\epsfxsize =4cm
\epsffile{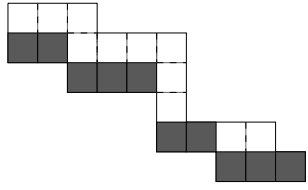}
\end{center}
\caption{\label{FIG1BIS} A skew diagram $\theta$ with its subdiagram 
$\theta\hspace{-3pt}\downarrow$ shaded}
\end{figure}

\begin{figure}[t]
\begin{center}
\leavevmode
\epsfxsize =4cm
\epsffile{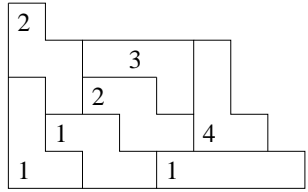}
\end{center}
\caption{\label{FIG2} A 4-ribbon tableau of shape $(8,7,6,6,1)$, weight $(3,2,1,1)$
and spin $9$}
\end{figure}

A ribbon is a connected skew Young diagram of width 1,
\ie which does not contain any $2\times 2$ square
(see Figure~\ref{FIG1}).
The rightmost and bottommost cell is called the origin of the ribbon.
An $n$-ribbon is a ribbon made of $n$ square cells.
Let $\theta$ be a skew Young diagram, and let
$\theta\hspace{-3pt}\downarrow$ be the 
horizontal strip made of the bottom cells of the columns of $\theta$
(see Figure~\ref{FIG1BIS}). 
We say that $\theta$
is a horizontal $n$-ribbon strip of weight $m$ if it can be tiled by
$m$ $n$-ribbons the origins of which lie in 
$\theta\hspace{-3pt}\downarrow$. 
One can check that if such a tiling exists, it is unique (see below Lemma~\ref{LEMCOMB} and
Figure~\ref{FIG7}).
Now, an $n$-ribbon tableau $T$ of shape $\lambda/\nu$ and weight
$\mu=(\mu_1,\,\ldots ,\,\mu_r)$ is defined as a chain of partitions
$$
\nu = \alpha^0 \subset \alpha^1 \subset \cdots \subset \alpha^r=\lambda
$$
such that $\alpha^i/\alpha^{i-1}$ is a horizontal $n$-ribbon strip of weight
$\mu_i$. 
Graphically, $T$ may be described by numbering each $n$-ribbon of
$\alpha^i/\alpha^{i-1}$ with the number $i$ (see Figure~\ref{FIG2}). 
We denote by $\tab_n(\lambda/\nu,\,\mu)$ the
set of $n$-ribbon tableaux of shape $\lambda/\nu$ and weight $\mu$.
Define the spin of a ribbon $R$ as $\sp(R):=h(R)-1$ 
where $h(R)$ is the height of $R$, and the spin of a ribbon tableau $T$
as the sum of the spins of its ribbons. Then
the sign $(-1)^{\sp(T)}$ depends only on the shape $\lambda/\nu$ of $T$
and is equal to the $n$-sign $\varepsilon_n(\lambda)$ when $\nu$
is empty. 
We denote it in general by $\varepsilon_n(\lambda/\nu)$.

\subsection{A $q$-analogue of the Littlewood-Richardson coefficients}

Using a classical formula for multiplying a monomial symmetric
function by a Schur function one can easily derive Eq.~(\ref{EQ32}).
Note that since $h_{s\mu}=h_\mu \ (s\in\SG_r)$, (\ref{EQ32})
implies that
\begin{equation}\label{EQSYM}
|\tab_n(\lambda,s\mu)| = |\tab_n(\lambda,\mu)|, \qquad (s\in \SG_r).
\end{equation} 
Let $\varphi_n$ denote the adjoint of the endomorphism $f \mapsto p_n(f)$
of the space of symmetric functions with respect to  $\< \cdot \,;\, \cdot \>$.
Recall from Section~\ref{SECT33} the definition of the $n$-quotient 
$(\lambda^{(0)},\cdots ,\lambda^{(n-1)})$
of a partition $\lambda$ of length $r$ such that
$\lambda + \rho \in \EASG_r\rho$ (for the action of  $\EASG_r$ 
on weights via $\pi_n$).
Then $(\ref{EQ30})$ is equivalent to 
\begin{equation}\label{EQ34}
\varphi_n(s_\lambda) = 
\varepsilon_n(\lambda) \, s_{\lambda^{(0)}}\cdots s_{\lambda^{(n-1)}},
\end{equation}
where we put $\varepsilon_n(\lambda) = 0$ if 
$\lambda + \rho \not \in \EASG_r\rho$.
By (\ref{EQ32}) we have
\[
|\tab_n(\lambda,\mu)|=
\varepsilon_n(\lambda)\,\<p_n(h_\mu)\,;\,s_\lambda\>
=\varepsilon_n(\lambda)\,\<h_\mu\,;\,\varphi_n(s_\lambda)\> .
\]
Recalling that the basis dual to $\{h_\mu\}$ is the basis $\{m_\mu\}$
of monomial symmetric functions, we thus have 
\begin{equation}
s_{\lambda^{(0)}}\cdots s_{\lambda^{(n-1)}}
=
\sum_{\mu\in P^+} |\tab_n(\lambda,\mu)|\,m_\mu.
\end{equation}
Hence, putting $x^T:=x_1^{\alpha_1}\cdots x_r^{\alpha_r}$
for a ribbon tableau $T$ of weight $\alpha=(\alpha_1,\ldots ,\alpha_r)$,
we get using (\ref{EQSYM})
\begin{equation}\label{EQ36}
s_{\lambda^{(0)}}\cdots s_{\lambda^{(n-1)}}(x_1,\ldots ,x_r)
=
\sum_{\mu \in P^+}\left(\sum_{\beta \in \SG_r\mu} 
\left(\sum_{T\in \tab_n(\lambda,\beta)}
x^T \right)\right)=
\sum_{T\in\tab_n(\lambda,\cdot)} x^T
\end{equation}
where we denote by $\tab_n(\lambda,\cdot)$ the set of $n$-ribbon 
tableaux of shape $\lambda$ (and arbitrary weight).

Now we can introduce a $q$-analogue of (\ref{EQ36}) via
the spin of ribbon tableaux and set
\begin{equation}\label{DEFG}
G(\lambda^{(0)},\ldots,\lambda^{(n-1)};q,x) :=
\sum_{T\in\tab_n(\lambda,\cdot)} q^{\sp(T)} x^T.
\end{equation}
It was proved in \cite{LLT2} that this function is symmetric with
respect to the variables $x_i$. (This is not clear a priori,
and the proof will be recalled below (see Remark~\ref{REMSYM}).)
Thus, expanding  on the basis of Schur functions we get
\begin{equation}\label{EQ38}
G(\lambda^{(0)},\ldots,\lambda^{(n-1)};q,x) = 
\sum_\nu c_{\lambda^{(0)},\ldots ,\lambda^{(n-1)}}^\nu (q) \, s_\nu(x),
\end{equation}
where the $c_{\lambda^{(0)},\ldots ,\lambda^{(n-1)}}^\nu (q)\in\Z[q]$ are 
some $q$-analogues of the Littlewood-Richardson coefficients.
The symmetric function (\ref{EQ38}) is the function 
$\widetilde G_\lambda^{(n)}(x;q)$ of \cite{LLT2} up to the change of
variable $q\mapsto q^{-2}$ and rescaling by an appropriate power of $q$.
\begin{figure}[t]
\begin{center}
\leavevmode
\epsfxsize =14.5cm
\epsffile{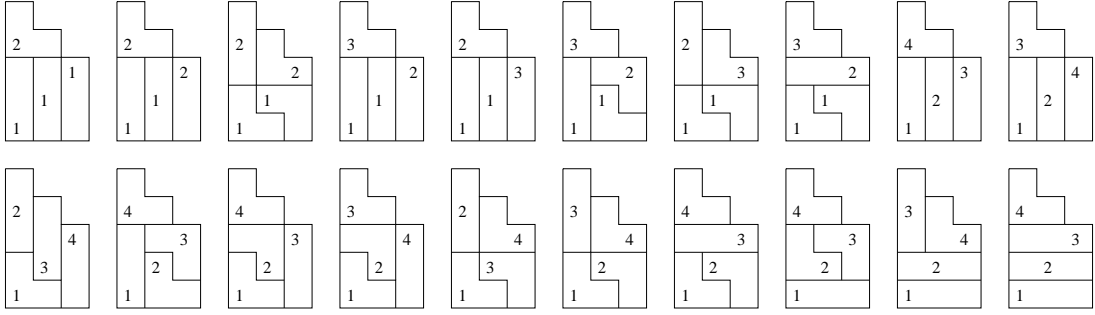}
\end{center}
\caption{\label{FIG3} The 3-ribbon tableaux of shape $(3,3,3,2,1)$ and dominant 
weight}
\end{figure}
\begin{example}{\rm The partition having as $3$-quotient 
$=((1),(1,1),(1))$
is $\mu = (3,3,3,2,1)$. Thus the symmetric function 
$G((1),(1,1),(1);q)$ is calculated by enumerating 
the $3$-ribbon tableaux of shape $\mu$ and dominant weight, 
and counting their spin (see Figure~\ref{FIG3}). One obtains
\begin{eqnarray*}
&&G((1),(1,1),(1);q) = q^7m_{(3,1)}+(q^7+q^5)m_{(2,2)} +
(2q^7+2q^5+q^3)m_{(2,1,1)}\\
&&\hskip 3.7cm +\ (3q^7+5q^5+3q^3+q)m_{(1,1,1,1)}\\
&&\hskip 3.3cm =
q^7s_{(3,1)}+q^5s_{(2,2)}+(q^5+q^3)s_{(2,1,1)}+qs_{(1,1,1,1)}.
\end{eqnarray*}\finex
}
\end{example}
We can now state our main result, which is the $q$-analog of 
Theorem~\ref{TH3}.
\begin{theorem}\label{TH4}
With the notation of Theorem~{\rm\ref{TH3}}
\[
c_{\mu^{(0)},\ldots ,\mu^{(n-1)}}^\lambda(q) = P^-_{\mu+\rho,n\lambda+\rho}(q)
=
\sum_{z\in\SG_{\nu,-n}}
(-1)^{\ell(z)}  P_{w(\mu+\rho,-n)z\,,\,w(n\lambda+\rho,-n)}(q).
\]
\end{theorem}
The next two sections will be devoted to the proof of Theorem~\ref{TH4}.
This proof does not rely on the Lusztig conjecture and thus will
give an independent proof of Theorem~\ref{TH3}.


\section{Canonical bases and Kazhdan-Lusztig polynomials} \label{SECT5}

\subsection{Another basis of $\P$}\label{SECT51}

The basis of $\P$ consisting of the weights $\lambda$ is
adapted to the Coxeter-type presentation of $\EAH_r$ in terms
of the generators $T_0,\ldots ,T_{r-1},\tau$.
There is another natural basis adapted to the Bernstein
presentation in terms of $T_1,\ldots ,T_{r-1},Y_1,\ldots,Y_r$,
which is defined as follows.
Fix $m\in \Z^*$ and consider the action of $\EAH_r$ via $\Pi_m$.
Every $\lambda \in P$ has a unique expression as
$
\lambda = m\beta + \gamma \ (\beta,\gamma \in P,\ 
\gamma \in \SG_r\A_{r,m}).
$
We define $V_\lambda := Y^\beta \gamma$. 
In other words, the basis $\{V_\lambda\}$ is characterized by
\begin{eqnarray}
&& V_\gamma = \gamma \hskip 2.7cm (\gamma \in  \SG_r\A_{r,m}), \\
&& Y^\beta V_\lambda = V_{\lambda+m\beta} \hskip 1.5cm (\lambda,\beta \in P).
\label{EQACTY}
\end{eqnarray}
\begin{example}{\rm \label{EX51}
Take $r=2$ and $m=-2$. Then
\begin{eqnarray*}
&&V_{(-1,-2)}=Y_2(-1,0)=T_1\tau (-1,0)=(-1,-2),\\
&&V_{(-2,-1)}=Y_1(0,-1) = T_0^{-1}\tau (0,-1)=(-2,-1),\\
&&V_{(2,-1)}=Y_1^{-1}(0,-1)=\tau^{-1}T_0(0,-1)=(2,-1),\\
&&V_{(-1,2)}=Y_2^{-1}(-1,0)=\tau^{-1} T_1^{-1} (-1,0)=(-1,2)+(q-q^{-1})(0,1).
\end{eqnarray*}
Take $r=3$ and $m=-3$. Then
\[
V_{(-2,-1,3)}=Y_3^{-1}(-2,-1,0)=\tau^{-1}T_1^{-1}T_2^{-1}(-2,-1,0)
\qquad\qquad
\]\[
=(-2,-1,3)+(q-q^{-1})(0,-1,1)+(q-q^{-1})(-2,0,2)
+(q-q^{-1})^2(-1,0,1).
\]\finex
}
\end{example}
\begin{remark}{\rm Let $n=|m|$.
The basis $\{V_\lambda\}$ can be naturally identified with
the basis of monomial tensors of a certain $U_q(\slchap_n)$-module
(see Section~\ref{SECT71}). \finex 
}
\end{remark}
As illustrated by Example~\ref{EX51}, in some cases the vectors 
$V_\lambda$ and $\lambda$ coincide. This is made more precise in
the following
\begin{proposition}\label{PROP5.1}
If $\lambda =  m\beta + \gamma$ as above with $\beta\in P^-$ then	
$V_\lambda = \lambda$.
In particular, if $m<0$ and $\lambda\in P^+$, or
$m>0$ and $\lambda\in P^-$, then $V_\lambda = \lambda$.
\end{proposition}
\proof
Put $s=w(\gamma,m)$ and $\nu = s^{-1}\gamma$. 
Then by (\ref{EQ24}) and Lemma~\ref{LEM2.2}
\[
V_\lambda = Y^\beta \gamma = T_\beta T_s \nu = T_{y^\beta s}\, \nu .
\]
For $\sigma \not = 1$ in  $\SG_{\nu,m}\subset \SG_r$ one has 
$\ell(s\sigma)>\ell(s)$ (because $s$ is minimal in its coset
$s\SG_{\nu,m}$) and $s\sigma \in \SG_r$. 
Hence by Lemma~\ref{LEM2.2} 
\[
\ell(y^\beta s \sigma) = \ell(y^\beta) + \ell(s\sigma)
>\ell(y^\beta) + \ell(s) = \ell(y^\beta s).
\]
Therefore $y^\beta s$ is also minimal in its coset,
that is
$w(\lambda,m) = y^\beta s$, and 
\[
V_\lambda = T_{y^\beta s}\,\nu = T_{w(\lambda,m)}\, \nu = \lambda.
\]
\cqfd

\noindent
The next proposition gives a key relation between the bar involution
and the basis $V_\lambda$. It will result from the following 
\begin{lemma}\label{LEM2.4}
Let $\beta\in P$ and $s\in\SG_r$. Then
\[
\bar{(Y^\beta T_s)} = T_{w_0}^{-1} Y^{w_0\beta} T_{w_0s}.
\]
\end{lemma}
\proof
Recall that $\ell(w_0s)+\ell(s)=\ell(w_0)$, hence
$T_{w_0s}T_{s^{-1}}=T_{w_0}$ and $\bar{T_s}=T_{w_0}^{-1}T_{w_0s}$.
Write $\beta = \beta'-\beta''$ with $\beta',\ \beta''\in P^+$.
By (\ref{EQBERN}) we have 
$\bar{Y^\beta}=T_{\beta'\rule{0mm}{2.7mm}}T^{-1}_{\beta''}$.
Hence,
$\bar{(Y^\beta T_s)}=T_{\beta'\rule{0mm}{2.7mm}}T^{-1}_{\beta''}
T_{w_0}^{-1}T_{w_0s}$.
Now, using Lemma~\ref{LEM2.2} we see that
$$T_{\beta'\rule{0mm}{2.7mm}}T^{-1}_{\beta''}T_{w_0}^{-1}
= T_{\beta'\rule{0mm}{2.7mm}}T_{w_0}^{-1}T^{-1}_{w_0\beta''}
=T_{w_0}^{-1} T_{w_0\beta'\rule{0mm}{2.7mm}}T^{-1}_{w_0\beta''}$$
because $\beta',\ \beta''\in P^+$.
Now $w_0\beta = (-w_0\beta'')-(-w_0\beta')$, with 
$-w_0\beta'',\ -w_0\beta'\in P^+$. 
Hence, using again (\ref{EQBERN}), 
$\bar{(Y^\beta T_s)}= T_{w_0}^{-1} Y^{w_0\beta} T_{w_0s}$.
\cqfd
\begin{proposition}\label{PROP5.2}
Let $\lambda \in P$ and let $\nu \in \A_{r,m}$ be the point congruent
to $\lambda$. Then
\[
\bar{V_\lambda} = q^{-\ell(w_{0,\nu})} T_{w_0}^{-1} V_{w_0\lambda},
\]
where $w_{0,\nu}$ is the longest element in the stabilizer $\SG_{\nu,m}$.
\end{proposition}

\proof
By Lemma~\ref{LEM2.4}, $\bar{V_\lambda} = \bar{(Y^\beta T_s)}\,\nu
= T_{w_0}^{-1} Y^{w_0\beta} T_{w_0s}\, \nu$.
The minimal length of an element $\sigma \in \SG_r$ such that
$\sigma\,\nu = (w_0s)\nu$ is $\ell(w_0s)-\ell(w_{0,\nu})$.
Hence $T_{w_0s}\, \nu = q^{-\ell(w_{0,\nu})}(w_0s) \nu$, and this proves
the proposition.
\cqfd

\begin{example}{\rm Take $m=-2$ and $\lambda = (2,0)$.
Then,
\begin{eqnarray*}
&&\bar{V_{(2,0)}}= \bar{Y_1^{-1}}(0,0) = 
\tau^{-1}T_0^{-1}(0,0) = (2,0)+(q-q^{-1})(0,2), \\
&&T_1^{-1}V_{(0,2)}=T_1^{-1}Y_2^{-1}(0,0)=T_1^{-1}\tau^{-1}T_1^{-1}(0,0) 
= q(2,0)+(q^2-1)(0,2).
\end{eqnarray*}
\finex}
\end{example} 
%

\subsection{Action of $\EAH_r$ on the basis $V_\lambda$} \label{SECT5BIS}

The next lemma allows to compute
the action of $\EAH_r$ on the basis $\{V_\lambda\}$.
\begin{lemma}\label{LEM5.2}
Let $i\in \{1,\dots ,r-1\}$ and $k\in\Z$. There holds
\[
T_iY_i^k = Y_{i+1}^kT_i + (q-q^{-1})Y_{i+1}{Y_i^k-Y_{i+1}^k\over Y_i-Y_{i+1}}.
\]
In other words,
\[
T_iY_i^k =
\left\{
\matrix{\displaystyle
 Y_{i+1}^kT_i + (q-q^{-1})\sum_{j=1}^kY_i^{k-j}Y_{i+1}^j, \quad (k\ge0), \cr
\displaystyle
 Y_{i+1}^kT_i + (q^{-1}-q)\sum_{j=1}^{-k}Y_i^{-j}Y_{i+1}^{j+k}, \quad (k<0).
}\right.
\]
\end{lemma}
\proof It follows from (\ref{EQB3}) (\ref{EQT3}) by a straightforward
computation. \cqfd

\noindent
Let $\lambda \in P$ and $\ 1\le i \le r-1$.
Write $\lambda = m\beta + \gamma$ with $\beta,\gamma \in P$ and
$\gamma \in  \SG_r \A_{r,m}$.
Then
$
V_\lambda = (\prod_{j\not = i,i+1}Y_j^{\beta_j})
(Y_iY_{i+1})^{\beta_{i+1}} Y_i^{\beta_i-\beta_{i+1}} V_\gamma.
$
Since $T_i$ commutes with $Y_j\ (j\not = i,i+1)$ and $Y_iY_{i+1}$,
we have
\[
T_i V_\lambda =  \left(\prod_{j\not = i,i+1}Y_j^{\beta_j}\right)
(Y_iY_{i+1})^{\beta_{i+1}} T_i Y_i^{\beta_i-\beta_{i+1}} V_\gamma .
\]
Thus to compute $T_i V_\lambda$ we can use the commutation relation of
Lemma~\ref{LEM5.2} with $k=\beta_i-\beta_{i+1}$ together with
the fact that since
$V_\gamma = \gamma$,
we have
\[
T_i V_\gamma =
\left\{
\matrix{ V_{s_i\gamma} \hfill&\quad \mbox{ \rm if }
\desc(\gamma,i,m)=\sgn(m),\hfill \cr
q^{-1}V_{s_i\gamma}\hfill& \quad \mbox{ \rm if }\desc(\gamma,i,m)=0,\hfill \cr
V_{s_i\gamma} + (q^{-1}-q)V_{\gamma}& \quad \mbox{ \rm if }
\desc(\gamma,i,m)=-\sgn(m).\hfill \
}\right.
\]

\subsection{Projection on the positive Weyl chamber}\label{SECT53}

From now on we fix $n \ge 2$ and we assume that $\EAH_r$ acts 
on $\P_r$ via $\Pi_{-n}$.
Introduce the $\Z[q,q^{-1}]$-submodule 
\[
\I_r:=\sum_{i=1}^{r-1} \Im C'_i \subset \P_r,
\]
and define $\F_r:=\P_r / \I_r$.
The image of $\lambda \in P$ in $\F_r$ under the natural projection 
\[
\pr : \P_r \longrightarrow \F_r
\]
will be denoted by $[\lambda]=[\lambda_1,\ldots ,\lambda_r]$.
For $v\in \P_r$ we have by definition 
\[
\pr(C'_iv)=0=\pr(T_iv)+q\,\pr(v).
\]
Hence taking $v=\lambda\in P$, we obtain that
if $\lambda_i<\lambda_{i+1}$ then	 
$[\lambda] = -q^{-1}[s_i\lambda]$, and
if $\lambda_i=\lambda_{i+1}$ then $[\lambda]=0$.
This implies that a spanning set of $\F_r$ is given by the $[\lambda]$
such that $\lambda_1>\lambda_2> \ldots >\lambda_r$. 
We put $P^{++}:=\{\lambda\in P \ | \ \lambda_1>\lambda_2> \ldots >\lambda_r\}.$
\begin{lemma} \label{LEMBASE}
$\{[\lambda]\ |\ \lambda \in P^{++}\}$ is a basis of $\F_r$.
\end{lemma}
\proof
Suppose that $\sum_{\lambda \in P^{++}} a_\lambda [\lambda] = 0$.
Then $\sum_{\lambda \in P^{++}} a_\lambda \lambda \in \I_r$.
Recall that 
\[
C_{w_0}=\sum_{s\in\SG_r}(-q)^{\ell(s)-\ell(w_0)}T_s
=\bar{C_{w_0}}=\sum_{s\in\SG_r}(-q)^{-\ell(s)+\ell(w_0)}T_s^{-1}
\]
satisfies $C_{w_0}C'_i=0 \ (1\le i \le r-1)$.
Hence $\I_r \subset \ker C_{w_0}$. Thus
\[
C_{w_0}(\sum_{\lambda\in P^{++}} a_\lambda \lambda)
= \sum_{\lambda\in P^{++},\ s\in\SG_r} 
a_\lambda (-q)^{-\ell(s)+\ell(w_0)}s\lambda = 0,
\]
which implies that $a_\lambda = 0$ for all $\lambda \in P^{++}$.
\cqfd

\noindent
Note that for $v\in\P_r$, $\bar{C'_i v} = C'_i\bar{v}$.
Hence $\bar{\I_r}\subset \I_r$ and one can define
a semi-linear involution on $\F_r$ by
\begin{equation}\label{EQINV}
\bar{\pr(v)}:=\pr(\bar{v})\qquad  (v\in \P_r).
\end{equation}
Let us define
\begin{equation}\label{DEFWEDGE}
|\lambda\>:=q^{-\ell(w_0)}\,\pr(V_\lambda).
\end{equation}
Then, by Proposition~\ref{PROP5.1}, for $\lambda\in P^{++}$ we
have $|\lambda\> = q^{-\ell(w_0)}\,[\lambda]$, so that
$\{|\lambda\>\ |\ \lambda \in P^{++}\}$ is also a basis of $\F_r$.
The next proposition shows that it is also useful to work with
the vectors $|\lambda\>$ associated with arbitrary weights $\lambda\in P$,
which can be thought of as some $q$-wedge products (see below
Section~\ref{SECT75}).

\begin{proposition}\label{PROP5.4}
For $\lambda\in P$, we have 
$\bar{|\lambda\>}=(-1)^{\ell(w_0)} q^{\ell(w_0)-\ell(w_{0,\nu})} 
\, |w_0\lambda\>.$
\end{proposition}
\proof
By Proposition~\ref{PROP5.2} we have
$
\bar{V_\lambda} = q^{-\ell(w_{0,\nu})} T_{w_0}^{-1} V_{w_0\lambda}.
$
But for all $v\in \P_r$,
\[\pr(T_{w_0}^{-1}v) = (-q)^{-\ell(w_0)}\,\pr(v).\] 
Thus,
\[
\bar{|\lambda\>}=
q^{\ell(w_0)}\,\pr(\bar{V_\lambda})=
(-1)^{\ell(w_0)} q^{-\ell(w_{0,\nu})} \pr(V_{w_0\lambda})
=(-1)^{\ell(w_0)} q^{\ell(w_0)-\ell(w_{0,\nu})} \, |w_0\lambda\>.
\]
\cqfd
\begin{remark}{\rm \label{REM58}
It is easy to check that
the exponent $\ell(w_0)-\ell(w_{0,\nu})$ of $q$ 
is equal to the number of pairs $(i,j)$ with 
$1\le i<j\le r$ such that $\lambda_i-\lambda_j$ is not 
divisible by $n$.
\finex}
\end{remark}

\noindent
The next proposition gives a set of straightening rules
to express an element $|\mu\>$ with $\mu \not \in P^{++}$
on the basis $\{|\lambda\>\ |\ \lambda \in P^{++}\}$.
\begin{proposition}\label{PROP5.5}
Let $\mu \in P$ be such that $\mu_i<\mu_{i+1}$.
Write $\mu_{i+1}=\mu_i+kn+j$ with $k\ge 0$ and $0\le j<n$. Then 
\begin{eqnarray}
&&|\mu\> = -|s_i\mu\> \hskip 7.1cm \mbox{ if } j=0, \label{STR1}\\
&&|\mu\> = -q^{-1}|s_i\mu\> \hskip 6.5cm \mbox{ if } k=0, \label{STR2}\\
&&|\mu\> = -q^{-1}|s_i\mu\> - |y_i^{-k}y_{i+1}^k\mu\>
-q^{-1}|y_i^ky_{i+1}^{-k}s_i\mu\> \hskip 1.5cm \mbox{otherwise.} \label{STR3}
\end{eqnarray}
\end{proposition}
\proof
To simplify the notation, let us write $l=\mu_i$ and $m=\mu_{i+1}$.
Since the relations only involve components $i$ and $i+1$ we 
shall also use the shorthand notations 
$(k,l)$ and $|k,l\>$ in place of 
$V_{(\mu_1,\ldots,\mu_{i-1},k,l,\mu_{i+2},\ldots,\mu_r)}\in \P_r$ and
$|(\mu_1,\ldots,\mu_{i-1},k,l,\mu_{i+2},\ldots,\mu_r)\>\in \F_r$.

Suppose $j=0$. 
It follows from Section~\ref{SECT5BIS} that $T_i (l,l) = q^{-1} (l,l)$.
Hence $(l,l)\in \Im C'_i$. 
Since $(Y_i^{-k}+Y_{i+1}^{-k})C'_i = C'_i(Y_i^{-k}+Y_{i+1}^{-k})$ we
also have $(Y_i^{-k}+Y_{i+1}^{-k})(l,l) = (m,l)+(l,m) \in \Im C'_i$,
and thus $|l,m\> + |m,l\> = 0$.

Suppose $k=0$. Then 
$T_i (l,m) = (m,l)$ by Section~\ref{SECT5BIS}, and 
$C'_i(l,m) = (m,l)+q(l,m) \in \Im C'_i$, which gives
$|l,m\> = -q^{-1}|m,l\>$.

Finally suppose that $j,k>0$. By the previous case
$(m,l+kn)+q(l+kn,m)\in \Im C'_i$. Applying $Y_i^k + Y_{i+1}^k$
we get that 
$(m,l) + (m-kn,l+kn) + q (l,m) + q(l+kn,m-kn) \in \Im C'_i$,
which gives the third claim. \cqfd
\begin{example}{\rm Take $r=2$ and $n=2$. Then
\[
|1,4\> = -q^{-1}\,|4,1\> - |3,2\> - q^{-1}\, |2,3\>,
\]
by Eq.~(\ref{STR3}), and 
$ |2,3\> = -q^{-1}\,|3,2\>$ by Eq.~(\ref{STR2}). Thus
\[
|1,4\> = -q^{-1}\,|4,1\> +(q^{-2}-1)\, |3,2\>.
\]
Hence, by Proposition~\ref{PROP5.4}, 
$\ \bar{|4,1\>} = |4,1\> +(q-q^{-1})\, |3,2\>.$
\finex}
\end{example}
For $\mu \in P^{++}$ write 
$\bar {|\mu\>} = \sum_{\lambda\in P^{++}} a_{\lambda\mu}(q)\,|\lambda\>.$
Using Proposition~\ref{PROP5.4} and Proposition~\ref{PROP5.5},
we easily see that the coefficients $a_{\lambda\mu}(q)$ satisfy the 
following properties 

\begin{corollary}\label{COR}
{\rm (i)} \ The coefficients $a_{\lambda\mu}(q)$ are invariant under
translation of $\lambda$ and $\mu$ by $\epsilon_1 + \cdots + \epsilon_r$.
Hence it is enough to describe the $a_{\lambda\mu}(q)$ for which
$\lambda-\rho$ and $\mu-\rho$ have non-negative components,
\ie $\lambda-\rho$ and $\mu-\rho$ are partitions. 

\smallskip
\noindent
{\rm (ii)} \ If $a_{\lambda\mu}(q) \not = 0$ then $\lambda \in \ASG_r \mu$.
In particular, if $\lambda-\rho$ and $\mu-\rho$ are partitions, 
they are partitions of the same integer $k$.

\smallskip
\noindent
{\rm (iii)} \ The matrix ${\bf A}_k$ with entries the $a_{\lambda\mu}(q)$
for which $\lambda-\rho$ and $\mu-\rho$ are partitions of $k$
is lower 
unitriangular if the columns and rows are indexed in decreasing 
lexicographic order.
\end{corollary}

\begin{example}{\rm
For $n=2$ and $r=3$, the matrices ${\bf A}_k$ for $k=2,3,4$ are
\[
\begin {array}{cc} \mbox{\small (4,1,0)}&\mbox{\small (3,2,0)}\\[3mm]
1&0\\q-q^{-1}&1\end {array}
\qquad
\begin {array}{ccc} \mbox{\small (5,1,0)}&\mbox{\small (4,2,0)}&
\mbox{\small (3,2,1)}\\[3mm]
1&0&0\\0&1&0\\q-q^{-1}&0&1
\end {array}
\qquad
\begin {array}{cccc} \mbox{\small (6,1,0)}&\mbox{\small (5,2,0)}&
\mbox{\small (4,3,0)}&\mbox{\small (4,2,1)}\\[3mm]
1&0&0&0\\
q-q^{-1}&1&0&0\\
q^{-2}-1&q-q^{-1}&1&0\\
0&q^{2}-1& q-q^{-1}&1\\
\end{array}
\]
\finex}
\end{example}

\subsection{Canonical bases of $\F_r$}\label{SECT54}

Let $\L^+$ (resp. $\L^-$) be the $\Z[q]$ (\resp $\Z[q^{-1}]$)-lattice 
in $\F_r$ with basis $\{|\lambda\>\ | \ \lambda \in P^{++} \}$. 
The fact that the matrix of the bar involution is unitriangular
on the basis $\{|\lambda\>\ | \ \lambda \in P^{++} \}$ implies
by a classical argument (see \cite{Lu6}, 7.10 and \cite{Du})
that 
\begin{theorem}
There exist bases $\{G^+_\lambda \ | \ \lambda \in P^{++}\}$, 
$\{G^-_\lambda\ | \ \lambda \in P^{++}\}$ of $\F_r$
characterized by:
\begin{quote}
{\rm (i)} \quad $\overline{G^+_\lambda}=G^+_\lambda$, \quad
$\overline{G^-_\lambda}=G^-_\lambda$,

{\rm (ii)} \quad $G^+_\lambda \equiv |\lambda\> \mod q\L^+$,
\quad $G^-_\lambda \equiv |\lambda\> \mod q^{-1}\L^-$.
\end{quote}
\end{theorem}
These bases were introduced in \cite{LT} (in the limit 
$r\rightarrow \infty$, \cf Section~\ref{SECT7}), using Proposition~\ref{PROP5.4} as the definition
of the bar involution on $\F$.
Set
\[
G^+_\mu =\sum_\lambda c_{\lambda,\mu}(q)\,|\lambda\>\,,\quad\quad
G^-_\lambda =\sum_\mu l_{\lambda,\mu}(-q^{-1})\, |\mu\>\,.
\]
Let ${\bf C}_k$ and ${\bf L}_k$ denote respectively the matrices 
with entries the coefficients $c_{\lambda\mu}(q)$ and
$l_{\lambda\mu}(q)$
for which $\lambda-\rho$ and $\mu-\rho$ are partitions of $k$.
\begin{example}{\rm For $r=3$ and $n=2$ we have
\[
{\bf C}_4=
\begin {array}{cccc} \mbox{\small (6,1,0)}&\mbox{\small (5,2,0)}&
\mbox{\small (4,3,0)}&\mbox{\small (4,2,1)}\\[3mm]
1&0&0&0\\
q&1&0&0\\
0&q&1&0\\
q&q^{2}& q&1\\
\end{array}
\qquad
{\bf L}_4=
\begin {array}{cccc} \mbox{\small (6,1,0)}&\mbox{\small (5,2,0)}&
\mbox{\small (4,3,0)}&\mbox{\small (4,2,1)}\\[3mm]
1&q&q^2&0\\
0&1&q&0\\
0&0&1&q\\
0&0& 0&1\\
\end{array}
\]
\finex}
\end{example}
Clearly, if $c_{\lambda,\mu}$ or $l_{\lambda,\mu}$ $\not = 0$, then
$\lambda$ and $\mu$ lie on the same orbit under $\EASG_r$.
Let $\nu$ be the point of $\A_{r,-n}$ on this orbit.
Write $\widehat{w}_\lambda := w(w_0\lambda,-n)w_{0,\nu}$
and similarly $\widehat{w}_\mu := w(w_0\mu,-n)w_{0,\nu}$.
The main result of this section is
\begin{theorem}[Varagnolo, Vasserot \cite{VV}]\label{THVV}
With the above notation, we have
\begin{equation}\label{EQE}
l_{\lambda,\mu} = P^-_{\mu,\lambda},
\end{equation}
a parabolic Kazhdan-Lusztig polynomial for the action of $\EASG_r$
on $P_r$ via $\pi_{-n}$, and
\begin{equation}\label{EQD}
c_{\lambda,\mu} = \sum_{s\in\SG_r}(-q)^{\ell(s)}
P_{s \widehat{w}_\lambda\, , \, \widehat{w}_\mu}.
\end{equation}
\end{theorem}

\begin{remark}{\rm (i) In view of Theorem~\ref{THDEO}, it follows from
Eq.~(\ref{EQD}) that $c_{\lambda,\mu}$ is also a parabolic 
Kazhdan-Lusztig polynomial of negative type with respect to 
the parabolic subgroup
$\SG_r$ of $\EASG_r$, (but for the right $\EAH_r$-module
${\bf 1}_{q^{-1}}\otimes_{H_r} \EAH_r$).
This agrees with the expression obtained by 
Goodman and Wenzl when $\mu - \rho$ is a $n$-regular 
partition \cite{GW} (see below Remark~\ref{REMTILT}).

\smallskip\noindent
(ii) Let $\bar{\F}_r$ denote the specialization of $\F_r$ at $q=1$.
Define a $\Z$-linear map $\iota$ from the Grothendieck group
of finite-dimensional representations of $U_\zeta(\gl_r)$ to
$\bar{\F}_r$ by $$\iota[W(\lambda)] = |\lambda+\rho\> \qquad
 (\lambda \in P^+_r).$$
Then comparing Theorem~\ref{THVV} and the Lusztig conjecture 
(\ref{EQ28}) we see that $\iota[L(\lambda)] = G^-_{\lambda+\rho}$.
}
\end{remark}
\proof
Consider the element 
$D_\lambda := \pr(C^-_\lambda)\in \F_r$.
Then $\bar{D_\lambda} = D_\lambda$ by (\ref{EQINV}).
Since $\lambda \in P^{++}$, $\desc(\lambda,i,-n)=1$ for 
all $i=1,\ldots,r-1$. Therefore using (\ref{OBSERV})
we see that 
\[
D_\lambda =  [r]! \sum_{\mu\in P^{++}} P^-_{\mu,\lambda}(-q^{-1})\,|\mu\>.
\]
Hence $(1/[r]!)D_\lambda$ is bar invariant and congruent to
$|\lambda\>$ modulo $q^{-1}\L^-$. Thus $D_\lambda = [r]!\,G^-_\lambda$
and (\ref{EQE}) is proved.

Next put $E_\mu := \pr(C^+_{w_0\mu})\in \F_r$.
Then $\bar{E_\mu} = E_\mu$.
We have
\[
E_\mu = 
\pr\left(\sum_{\alpha \in \EASG_r \nu} P^+_{\alpha,w_0\mu}\,\alpha\right)
=\sum_{\lambda\in P^{++}}\left(
\sum_{s\in\SG_r} (-q)^{-\ell(s)}P^+_{s\lambda,w_0\mu}\right)
q^{\ell(w_0)} \,|\lambda\>
.
\]
This shows that $E_\mu \equiv (-1)^{\ell(w_0)}\,|\mu\> \mod q\L^+$.
Hence, $E_\mu = (-1)^{\ell(w_0)}\,G^+_\mu$.
It follows that 
\[
c_{\lambda,\mu} = \sum_{s\in\SG_r} (-q)^{\ell(w_0)-\ell(s)}
P^+_{s\lambda,w_0\mu}
=
 \sum_{\sigma\in\SG_r} (-q)^{\ell(\sigma)}
P_{w(\sigma w_0 \lambda,-n)w_{0,\nu}\, , \, w(w_0\mu,-n)w_{0,\nu}}
\]
by Theorem~\ref{THDEO}. Finally, since $w_0\lambda \in P^{-\,-}$ we have
$w(\sigma w_0 \lambda,-n) = \sigma w(w_0 \lambda,-n)$
for all $\sigma\in\SG_r$, and we get (\ref{EQD}).
\cqfd


\section{A $q$-analogue of the tensor product theorem} \label{SECT6}

\subsection{Action of $Z(\EAH_r)$ on $\F_r$}
\label{SECT55}

By a result of Bernstein  (see \cite{Lu2}, Th.~8.1),
the center $Z(\EAH_r)$ of $\EAH_r$ 
is the algebra of symmetric Laurent polynomials in the elements $Y_i$.
Clearly, $Z(\EAH_r)$ leaves invariant the submodule $\I_r$.
It follows that $Z(\EAH_r)$ acts on $\F_r = \P_r/\I_r$.
This action can be computed via (\ref{EQACTY}) and (\ref{DEFWEDGE}).
In particular $B_k = \sum_{i=1}^r Y_i^k$ acts by
\begin{equation}
B_k \,|\lambda\>
= \sum_{j=1}^r|\lambda-nk\epsilon_j\>, \qquad\qquad (k\in\Z^*).\label{EQ50}
\end{equation}
Note that the right-hand side of (\ref{EQ50}) may involve terms
$|\mu\>$ with $\mu\not \in P^+$ which have to be expressed on the
basis $\{|\lambda\>\ |\ \lambda \in P^{++}\}$ by repeated applications
of Proposition~\ref{PROP5.5}.

\begin{example}\label{EX63}{\rm Take $r=4$ and $n=2$. We have 
\[
B_{-2}\,|3,2,1,0\>
= |7,2,1,0\> + |3,6,1,0\> + |3,2,5,0\> + |3,2,1,4\>.
\]
By Proposition~\ref{PROP5.5}, 
\begin{eqnarray*}
&&|3,6,1,0\> = -q^{-1}\,|6,3,1,0\> + (q^{-2}-1)\,|5,4,1,0\>,\\
&&|3,2,5,0\> = -q^{-1}\,|3,5,2,0\> + (q^{-2}-1)\,|3,4,3,0\>
=q^{-1}\,|5,3,2,0\>,\\
&&|3,2,1,4\> = -q^{-1}\,|3,2,4,1\> + (q^{-2}-1)\,|3,2,3,2\>
=-q^{-2}\,|4,3,2,1\>,
\end{eqnarray*}
which yields
\begin{eqnarray*}
&&
B_{-2}\,|3,2,1,0\>
=
|7,2,1,0\>  -q^{-1}\,|6,3,1,0\>\hskip 2.5cm\\ 
&&\hskip 4cm +\, (q^{-2}-1)\,|5,4,1,0\>
+ q^{-1}\,|5,3,2,0\> -q^{-2}\,|4,3,2,1\>.
\end{eqnarray*}
\finex
}
\end{example}
The compatibility of the bar involution with this action is
given by the next
\begin{proposition}\label{PROPBARFOCKB}
For $u\in \F_r$ and $z\in Z(\EAH_r)$ one has
\[
\bar{z u} = z \bar{u}.
\]
\end{proposition}
\proof
Since $z$ is a symmetric Laurent polynomial in the $Y_i$, we see
using Lemma~\ref{LEM2.4} that $\bar{z} = T_{w_0}^{-1} z T_{w_0} = z$.
\cqfd

\subsection{Action of $Z(\EAH_r)$ and ribbon tableaux}

We shall now show that the straightening relations can be avoided
provided that one uses appropriate linear bases of $Z(\EAH_r)$.
For $d \in [1,r] := \{1,2,\ldots ,r\}$ and $m\in \N^*$ define 
\begin{eqnarray}
&&\U_d := \sum_{1\le i_1< i_2< \ldots < i_d \le r}
Y_{i_1}Y_{i_2}\cdots Y_{i_d},\\[2mm]
&&\V_d := \sum_{1\le i_1< i_2< \ldots < i_d \le r}
Y_{i_1}^{-1}Y_{i_2}^{-1}\cdots Y_{i_d}^{-1},\label{DEFV}\\[2mm]
&&\UU_m := \sum_{1\le i_1\le i_2\le \ldots \le i_m \le r}
Y_{i_1}Y_{i_2}\cdots Y_{i_m},\\[2mm]
&&\VV_m := \sum_{1\le i_1\le i_2\le \ldots \le i_m \le r}
Y_{i_1}^{-1}Y_{i_2}^{-1}\cdots Y_{i_m}^{-1}.
\end{eqnarray}
For $\alpha \in [1,r]^s$ set 
$\U_\alpha := \U_{\alpha_1}\cdots \U_{\alpha_s},\ 
\V_\alpha := \V_{\alpha_1}\cdots \V_{\alpha_s}$, 
and for $\beta \in {\N^*}^s$ set
$\UU_\beta := \UU_{\beta_1}\cdots \UU_{\beta_s},\
\VV_\beta := \VV_{\beta_1}\cdots \VV_{\beta_s}$.
In other words, using the notation of \cite{Mcd}
for symmetric functions,
\begin{eqnarray*}
&&\U_\alpha =  e_\alpha(Y_1,\ldots,Y_r), \qquad 
\V_\alpha = e_\alpha(Y_1^{-1},\ldots,Y_r^{-1}), \\
&&\UU_\beta = h_\beta(Y_1,\ldots,Y_r), \qquad
\VV_\beta = h_\beta(Y_1^{-1},\ldots,Y_r^{-1}).
\end{eqnarray*}
The following formulas were obtained in \cite{LLT2}.
They will allow us to relate ribbon tableaux to Kazhdan-Lusztig 
polynomials.
Put
\begin{equation}
L_{\lambda/\nu\,,\,\mu}^{(n)}(q) :=  \sum_{T\in\tab_n(\lambda/\nu\,,\,\mu)}
q^{\sp(T)}.
\end{equation}
\begin{theorem}\label{THRIB}
Let $\nu \in \PP^+_r$ and $\alpha \in [1,r]^s$. Set
$k=|\alpha|:=\alpha_1 + \cdots + \alpha_s$.
We have
\begin{eqnarray}
&&\U_\alpha \, |\nu+\rho\> = (-q)^{-(n-1)k}\, \label{EQU}
\sum_{\mu\in\PP^+_r}
L_{\nu'/\mu'\,,\,\alpha}^{(n)}(-q) \,\, |\mu+\rho\>, \\
&&\V_\alpha \, |\nu+\rho\> = (-q)^{-(n-1)k}\,\sum_{\lambda\in\PP^+_r}
L_{\lambda'/\nu'\,,\,\alpha}^{(n)}(-q) \,\, |\lambda+\rho\>,\label{EQV}
\end{eqnarray}
where for $\lambda\in \PP^+_r$, $\lambda'$ denotes the conjugate
partition.
\end{theorem}
Note that in (\ref{EQU}) (\ref{EQV}) $\lambda'$, $\mu'$, $\nu'$
may be partitions of length $s>r$.
\begin{figure}[t]
\begin{center}
\leavevmode
\epsfxsize =8cm
\epsffile{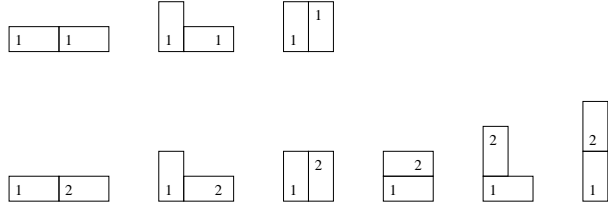}
\end{center}
\caption{\label{FIG6} The domino tableaux of weight $(2)$ and $(1,1)$}
\end{figure}

\begin{example}{\rm Let us redo the calculation of Example~\ref{EX63}
using domino tableaux.
Clearly, $B_{-2} = \V_{(1,1)}-2\V_{(2)}$.
Now, applying Theorem~\ref{THRIB} we have  
\begin{eqnarray*}
&&\V_{(2)}\,|\rho\> = q^{-2}\,|(1,1,1,1)+\rho\>
-q^{-1}\,|(2,1,1)+\rho\>
+|(2,2)+\rho\>, \hskip2.5cm \\[2mm]
&&\V_{(1,1)}\,|\rho\> = 
q^{-2}\,|(1,1,1,1)+\rho\>
-q^{-1}\,|(2,1,1)+\rho\>
+(1+q^{-2})\,|(2,2)+\rho\>\hskip1cm \\
&&\hskip 8cm
-q^{-1}\,|(3,1)+\rho\>
+|(4)+\rho\>
\end{eqnarray*}
as illustrated by Figure~\ref{FIG6},
and we recover the result of Example~\ref{EX63}.
}\finex
\end{example}
The proof of Theorem~\ref{THRIB} is based on the following
simple combinatorial lemma.
\begin{lemma}\label{LEMCOMB}
Let $\lambda,\ \nu\in\PP^+_r$ and $k\in[1,r]$.
Put $\beta = \epsilon_1+\cdots +\epsilon_k$. 
The skew Young diagram $\lambda'/\nu'$ is a horizontal $n$-ribbon 
strip of weight $k$ if and only if there exist $s,\sigma \in\SG_r$ such
that $\nu+\rho+s(n\beta)=\sigma (\lambda+\rho)$.
If this is the case, 
\[
\ell(\sigma) = (n-1)k - \sp(\lambda'/\nu').
\]
\end{lemma}
\proof The proof is elementary and is left to the reader. \cqfd
\begin{figure}[t]
\begin{center}
\leavevmode
\epsfxsize =4cm
\epsffile{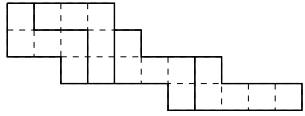}
\end{center}
\caption{\label{FIG7} A horizontal 5-ribbon strip of weight 4 and spin 7}
\end{figure}

\begin{example}{\rm Take $r=11$, $\lambda = (4,4,4,4,3,2,2,2,1,1,1)$
and $\nu = (2,2,1,1,1,1)$.
Then $\lambda'/\nu' = (11,8,5,4)/(6,2)$ is a horizontal 5-ribbon strip
of weight 4.
Indeed 
\[
(12,11,14,13,7,6,9,3,2,1,5) = \nu + \rho + (0,0,5,5,0,0,5,0,0,0,5)
\]
is a permutation of $\lambda + \rho$. This permutation has length
9, thus $\sp(\lambda'/\nu') = 4.4-9 = 7$, as can be checked on 
Figure~\ref{FIG7}.
\finex}
\end{example}

\medskip\noindent {\it Proof of Theorem~{\rm\ref{THRIB}}--- \ }
Since $\V_\alpha := \V_{\alpha_1}\cdots \V_{\alpha_s}$, it is enough to
prove the theorem in the case $\alpha = (k)$.
Let $\beta = \epsilon_1+\cdots +\epsilon_k$.
Observe that we can reformulate (\ref{DEFV}) as
$\V_k= \sum_{\zeta \in \SG_r\beta} Y^{-\zeta}$.
Hence we have
\[
\V_k\,|\nu + \rho\>
= 
\sum_{\gamma \in \SG_r n\beta} |\nu + \rho + \gamma\>.
\]
If $\xi := \nu + \rho + \gamma \not \in P^{++}$ we have to use 
the straightening relations of Proposition~\ref{PROP5.5} to express
$|\xi\>$ on the basis $\{|\lambda\>\ | \ \lambda\in P^{++}\}$.
But if $\xi_i<\xi_{i+1}$ then clearly we must have
$\xi_i<\xi_{i+1}<\xi_i + n$, and we need only the simple relation (\ref{STR2}).
It follows that $|\xi\> = (-q)^{-\ell(\sigma)}|\lambda+\rho\>$,
where $\lambda+\rho$ is the decreasing reordering of $\xi$ and $\sigma$
is the permutation mapping $\xi$ into $\lambda+\rho$.
By Lemma~\ref{LEMCOMB}, $\ell(\sigma) = (n-1)k - \sp(\lambda'/\nu')$
and we are done.
The proof for $\U_k$ is similar. \cqfd

\noindent
We now deduce from Theorem~\ref{THRIB} similar formulas for the operators 
$\UU_\beta$ and $\VV_\beta$.
\begin{theorem}\label{THRIB2}
Let $\nu \in \PP^+_r$ and $\beta \in {\N^*}^s$. 
We have
\begin{eqnarray}
&&\UU_\beta \, |\nu+\rho\> = \sum_{\mu\in\PP^+_r}
L_{\nu/\mu\,,\,\beta}^{(n)}(-q^{-1})
\, \, |\mu+\rho\>, \\
&&\VV_\beta \, |\nu+\rho\> = \sum_{\lambda\in\PP^+_r}
L_{\lambda/\nu\,,\,\beta}^{(n)}(-q^{-1})
\,\, |\lambda+\rho\>.
\end{eqnarray}
\end{theorem}
\proof
Again, it is enough to prove this for $\beta = (k)$.
Recall that a composition of $k\in\N$ is an ordered partition of $k$,
that is, a sequence $\alpha = (\alpha_1,\ldots ,\alpha_s)$ of positive integers
such that $\sum_i \alpha_i = k$. We denote this by $\alpha\models k$
and we call $s$ the length $l(\alpha)$ of $\alpha$.
There is a classical formula for expressing the symmetric function $h_k$
in terms of the $e_\alpha$, namely 
\[
h_k = \sum_{\alpha\models k} (-1)^{k-l(\alpha)}\, e_\alpha.
\]
Thus by Theorem~\ref{THRIB}, we have
\[
\VV_k\,|\nu + \rho \> = (-q)^{-(n-1)k}
\sum_\lambda \left(
\sum_{\alpha\models k} (-1)^{k-l(\alpha)} 
L^{(n)}_{\lambda'/\nu'\,,\,\alpha}(-q) \right) |\lambda+\rho\> .
\]
Recall that for a ribbon tableau $T$, 
$(-1)^{\,\sp(T)}=\varepsilon_n(\lambda/\nu)$ depends
only on the shape $\lambda/\nu$ of $T$. 
It is clear that 
$\varepsilon_n(\lambda'/\nu') = (-1)^{(n-1)k}
\varepsilon_n(\lambda/\nu)$, hence we are reduced to prove that
\[
q^{-(n-1)k}\sum_{\alpha\models k} (-1)^{k-l(\alpha)} 
L^{(n)}_{\lambda'/\nu'\,,\,\alpha}(q)
=
L^{(n)}_{\lambda/\nu\,,\,k}(q^{-1}).
\]
\begin{figure}[t]
\begin{center}
\leavevmode
\epsfxsize =9cm
\epsffile{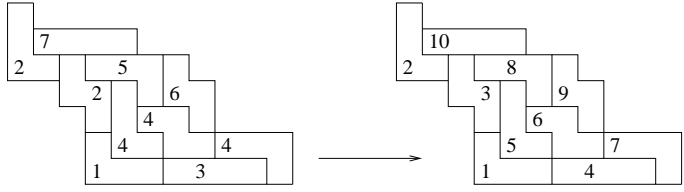}
\end{center}
\caption{\label{FIG8} Standardization $T \rightarrow \T$ of a ribbon tableau }
\end{figure}
To do this, we associate with each ribbon tableau T of weight $\alpha$
a standard ribbon tableau $\T$ of weight $(1,\ldots ,1)$ as follows.
Consider two ribbons $R$ and $R'$ of $T$, numbered $i$ and $i'$
respectively. We say that $R<R'$ if $i<i'$, or $i=i'$ and $R$ is to the left
of $R'$. 
Clearly this is a total order.
Now $\T$ is the tableau with the same shape and inner ribbon structure 
as $T$, whose ribbons are numbered $1,2,3,\ldots$ in the
previous total order (see Figure~\ref{FIG8}).

Let us fix a skew shape $\lambda'/\nu'$ and consider the set $\E$ of all 
ribbon tableaux of this shape and of arbitrary weight $\alpha\models k$.
For $T\in\E$ of weight $\alpha$, write $\epsilon(T):=(-1)^{k-l(\alpha)}$.
We want to prove that
\begin{equation}\label{TOPROVE}
\sum_{T\in\E} \epsilon(T)\,q^{\sp(T)} =
\left\{
\begin{array}{ll}
q^{(n-1)k}\,L^{(n)}_{\lambda/\nu\,,\,k}(q^{-1}) & \mbox{\rm if
$\lambda/\nu$ is a horizontal $n$-ribbon strip,}\\ [2mm]
0  & \mbox{\rm otherwise.}
\end{array}
\right.
\end{equation}
Let $\T\in \E$ be a standard tableau, and let $\E_\T\subset \E$ denote
the subset consisting of those tableaux $T$ whose standardization gives $\T$.
We say that $\T$ is a column if for all $i=1,\ldots , k-1$ the ribbon $R_{i+1}$
numbered $i+1$ lies above the ribbon $R_i$ numbered $i$, that is,
if the origin of $R_{i+1}$ lies in a row strictly above
the origin of $R_i$.
Eq.~(\ref{TOPROVE}) will follow from the more precise statement
\begin{equation}\label{TOPROVE2}
\sum_{T\in\E_\T} \epsilon(T)\,q^{\sp(T)} =
\left\{
\begin{array}{ll}
q^{\sp(\T)} & \mbox{\rm if
$\T$ is a column,}\\ [2mm]
0  & \mbox{\rm otherwise.}
\end{array}
\right.
\end{equation}
Now this is very easy. First, by definition all $T\in \E_\T$ have the
same inner ribbon structure, hence the same spin, and we can simplify
the powers of $q$ of both sides of (\ref{TOPROVE2}).
Then we only have to observe that a tableau $T\in\E_\T$ is specified
by the numbering of its ribbons, \ie by a map 
$f_T:\ [1,k] \longrightarrow [1,k]$ satisfying 
\begin{quote}
(i) $f_T(i+1) = f_T(i)$ or $f_T(i+1) = f_T(i)+1$,

(ii) if $R_{i+1}$ lies above $R_i$ in $\T$ then $f_T(i+1) = f_T(i)+1$.
\end{quote}
Let $a(\T)$ be the number of $i$'s such that $R_{i+1}$ is not above
$R_i$. 
Then clearly $|\E_\T|=2^{a(\T)}$ and more precisely the number
of $T\in\E_T$ such that $f_T(k)=j$ (\ie $\epsilon(T) = (-1)^{k-j}$)
is equal to ${a(\T)\choose j}$.
Hence by the binomial theorem
\[
\sum_{T\in\E_\T} \epsilon(T) =
\left\{
\begin{array}{ll}
1 & \mbox{\rm if
$a(\T)=0$, \ie $\T$ is a column,}\\ [2mm]
0  & \mbox{\rm otherwise.}
\end{array}
\right.
\]
To finish the proof we need only note that $\lambda/\nu$ is a
horizontal $n$-ribbon strip if and only if there exists a 
(necessarily unique) column tableau $\T$ of shape $\lambda'/\nu'$,
and in this case $\sp(\T) = (n-1)k - \sp(\lambda/\nu)$.
\cqfd

\begin{remark}\label{REMSYM}{\rm Since the $\VV_m$ commute, 
$\VV_\beta$ is invariant under permutation of $\beta$. 
Hence Theorem~\ref{THRIB2} implies that $L_{\lambda/\nu\,,\,\beta}^{(n)}(q)$
is also invariant under permutation of $\beta$.
This proves that the polynomial (\ref{DEFG}) is indeed symmetric.\finex
}
\end{remark}

\subsection{Action of $Z(\EAH_r)$ on the canonical basis $\{ G^-_\lambda\}$}

For $\lambda\in\PP^+_r$, define 
\[
S_\lambda := s_{\lambda}(Y_1^{-1},\ldots,Y_r^{-1})  \in Z(\EAH_r),
\]
where $s_\lambda$ is the Schur function.
The following theorem is a formal analogue of Theorem~\ref{TH2}.
\begin{theorem}\label{THSTPR}
Let $\lambda \in P^+$. Write $\lambda = \lambda^{(0)} + n\lambda^{(1)}$, where
$\lambda^{(0)}$ is $n$-restricted, that is,
\[
0\le \lambda_i^{(0)}-\lambda_{i+1}^{(0)} <n \qquad (1\le i \le r-1),
\]
and $\lambda_r^{(1)}\ge 0$.
Then
$
G^-_{\lambda+\rho} = S_{\lambda^{(1)}}\,G^{-}_{\lambda^{(0)}+\rho}.
$
\end{theorem}
\begin{figure}[t]
\begin{center}
\leavevmode
\epsfxsize =9cm
\epsffile{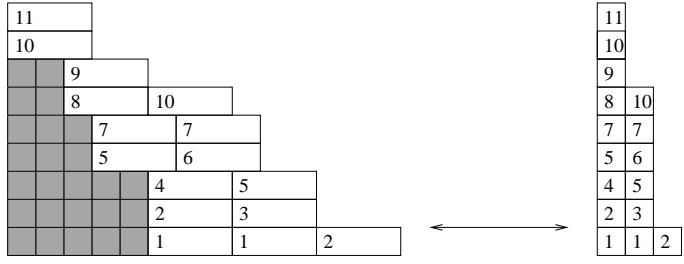}
\end{center}
\caption{\label{FIG9} Correspondence between $n$-ribbon tableaux of 
spin 0 and $n$-restricted inner shape, and ordinary tableaux}
\end{figure}
\proof
By definition of the basis $G^-$, we have to prove that 
$F_\lambda:=S_{\lambda^{(1)}}\,G^{-}_{\lambda^{(0)}+\rho}$
satisfies
\[
\bar{F_\lambda} = F_\lambda \qquad \mbox{\rm and} \qquad
F_\lambda \equiv |\lambda + \rho\> \quad \mod q^{-1}\L^-.
\]
The first property is clear by Proposition~\ref{PROPBARFOCKB}.
Indeed, $S_\lambda$ is a $\Q$-linear combination of products
of elements $B_{-i}$.
To prove the second one we observe that by Theorem~\ref{THRIB2}
for all $\nu \in \PP^+_r$ and $\alpha \in {\N^*}^s$,
$\VV_\alpha \, |\nu+\rho\> \in \L^-$.
Since $G^-_{\lambda^{(0)}+\rho} \equiv |\lambda^{(0)} + \rho\> \mod q^{-1}\L^-$,
and $S_{\lambda^{(1)}}$ is a $\Z$-linear combination of operators
$\VV_\alpha$ we thus have 
\[
F_\lambda \equiv S_{\lambda^{(1)}} |\lambda^{(0)} + \rho\> 
\quad \mod q^{-1} \L^-.
\]
In fact Theorem~\ref{THRIB2} implies \[
\VV_\alpha \, |\nu+\rho\> \equiv \sum_{T} 
|{\rm sh}(T)+\rho\>\quad \mod q^{-1} \L^-,
\]
where the sum is over the $n$-ribbon tableaux of weight $\alpha$,
spin $0$ and inner shape $\nu$, and ${\rm sh}(T)$ stands 
for the outer shape of $T$.
Therefore for all $\alpha$
\[
\VV_\alpha \,  |\lambda^{(0)} + \rho\>
\equiv 
 \sum_{T'} |{\rm sh}(T')+\rho\>\quad \mod q^{-1} \L^-,
\]
where the sum is over the $n$-ribbon tableaux $T'$ of weight $\alpha$
with inner shape ${\lambda^{(0)}}$ whose ribbons are all horizontal.
Now $\lambda^{(0)}$ being $n$-restricted, 
there is an obvious
bijection between the set of these tableaux $T'$ and the
set $\tab(\cdot,\alpha)$ of ordinary Young tableaux of weight $\alpha$ 
(see Figure~\ref{FIG9}).
Hence, for all $\alpha$
\[
\VV_\alpha \,  |\lambda^{(0)} + \rho\>
\equiv 
 \sum_\beta |\tab(\beta,\alpha)| \, |\lambda^{(0)}+n\beta+\rho\>
\quad \mod q^{-1} \L^-.
\]
Comparing with the well-known formula
$
h_\alpha = \sum_\beta |\tab(\beta,\alpha)| \, s_\beta 
$
which yields
\[
\VV_\alpha = \sum_\beta |\tab(\beta,\alpha)| \, S_\beta,
\]
we deduce that for all $\beta$,
\[
S_\beta \,  |\lambda^{(0)} + \rho\>
\equiv 
|\lambda^{(0)}+n\beta+\rho\> \quad \mod q^{-1} \L^-,
\]
and putting $\beta = \lambda^{(1)}$ we are done. \cqfd

\subsection{Proof of Theorem~\ref{TH4}}

Let us write in the ring of
symmetric functions $s_\lambda = \sum_\nu \kappa_{\lambda,\nu}\, h_\nu$.
Then we also have $m_\nu =  \sum_\lambda \kappa_{\lambda,\nu}\, s_\lambda$.
Hence
\[
G(\mu^{(0)},\ldots,\mu^{(n-1)};q)
:=\sum_\nu L^{(n)}_{\mu,\nu}(q)\, m_\nu
=
\sum_\lambda
\left(
\sum_\nu  \kappa_{\lambda,\nu}\, L^{(n)}_{\mu,\nu}(q) \right) s_\lambda,
\]
which gives
\[
c_{\mu^{(0)},\ldots ,\mu^{(n-1)}}^\lambda(q) =
\sum_\nu  \kappa_{\lambda,\nu}\, L^{(n)}_{\mu,\nu}(q).
\]
On the other hand, by Theorem~\ref{THSTPR} and Theorem~\ref{THVV} we have
\[
S_\lambda\,|\rho\> = G^-_{n\lambda + \rho}
= \sum_\mu P^-_{\mu+\rho,n\lambda + \rho}(-q^{-1})\,|\mu+\rho\>.
\]
Finally, using Theorem~\ref{THRIB2} we get
\begin{eqnarray*}
&&S_\lambda\,|\rho\> =  \sum_\nu \kappa_{\lambda,\nu}\, \VV_\nu \, |\rho\>
=\sum_\mu\left(
\sum_\nu  \kappa_{\lambda,\nu}\, L^{(n)}_{\mu,\nu}(-q^{-1})\right)
|\mu+\rho\>\\[1mm]
&&\hskip 1.1cm
=\sum_\mu c_{\mu^{(0)},\ldots ,\mu^{(n-1)}}^\lambda(-q^{-1}) \, |\mu+\rho\>,
\end{eqnarray*}
and by comparing the coefficients of $|\mu+\rho\>$ we have
\[
c_{\mu^{(0)},\ldots ,\mu^{(n-1)}}^\lambda(q) = P^-_{\mu+\rho,n\lambda+\rho}(q).
\]
\cqfd

\section{An inversion formula for Kazhdan-Lusztig polynomials }
\label{SECT7}

In this section we extend the coefficients to $\Q(q)$ and work with
\[
\PB_r:=\Q(q)\otimes_{\Z[q,q^{-1}]}\P_r,\quad
\FB_r:=\Q(q)\otimes_{\Z[q,q^{-1}]}\F_r,\quad
\EAHB_r:=\Q(q)\otimes_{\Z[q,q^{-1}]}\EAH_r.
\]

\subsection{Action of $U_q(\widehat{\hbox{\twelvegoth sl}}_n)$ on
the weight lattice of $\hbox{\twelvegoth gl}_r$}
\label{SECT71}

Let $U_q(\slchap_n)$ be the quantum enveloping algebra of the
affine Lie algebra $\slchap_n$.
This is a $\Q(q)$-algebra with generators $e_i,\ f_i,\ q^{\pm h_i}\
(0\le i \le n-1)$. The standard relations can be found
for example in \cite{LLT} and will be omitted.
There is a canonical involution $x\mapsto \overline x$ of $U_q(\slchap_n)$
defined as the unique ring homomorphism such that $\bar{q} = q^{-1}$,
$\bar{e_i} = e_i$, and $\bar{f_i} = f_i$.

Using the basis $\{V_\lambda\}$ for $m=-n$ one can define an action of
$U_q(\slchap_n)$ on $\PB_r$.
First we start with the trivial case $r=1$, where $\PB_r$ reduces to
$\PB_1 = \bigoplus_{l\in\Z} \Q(q) V_l.$
It is immediate to check that the formulas
\[
f_iV_l := \delta_{l\equiv i}\, V_{l+1}, \qquad
e_iV_l := \delta_{l\equiv i+1}\, V_{l-1},\qquad
q^{\pm h_i}V_l := q^{\pm(\delta_{l\equiv i}-\delta_{l\equiv i+1})}\,V_l,
\]
extend to an action of $U_q(\slchap_n)$ on $\PB_1$.
Here,  $a\equiv b$ means congruence modulo $n$ and $\delta_{a\equiv b}$
is the Kronecker $\delta$ equal to 1 if $a\equiv b$ and to 0 otherwise.
Then using the comultiplication
\begin{equation}\label{COM}
\Delta f_i = f_i\otimes 1 + q^{h_i}\otimes f_i , \quad
\Delta e_i = e_i\otimes q^{-h_i} + 1\otimes e_i , \quad
\Delta q^{\pm h_i} = q^{\pm h_i} \otimes q^{\pm h_i},
\end{equation}
and identifying $\PB_r$ with $\PB_1^{\otimes r}$ by
$V_\lambda \mapsto V_{\lambda_1}\otimes \cdots \otimes V_{\lambda_r}$,
we obtain the following formulas
\begin{eqnarray}
&&f_iV_\lambda := \sum_{\stackrel{\scriptstyle j=1}{\lambda_j\equiv i}}^r
q^{\sum_{k=1}^{j-1}(\delta_{\lambda_k\equiv i}-\delta_{\lambda_k\equiv i+1})}
\,V_{\lambda+\epsilon_j},\label{EQ44}\\
&&e_iV_\lambda := \sum_{\stackrel{\scriptstyle j=1}{\lambda_j\equiv i+1}}^r
q^{-\sum_{k=j+1}^r(\delta_{\lambda_k\equiv i}-\delta_{\lambda_k\equiv i+1})}
\,V_{\lambda-\epsilon_j}.\label{EQ45}
\end{eqnarray}
\begin{proposition} \label{PROPCOMM}
This action of $U_q(\slchap_n)$ on $\PB_r$
commutes with the action of $\EAHB_r$ via $\Pi_{-n}$.
\end{proposition}
\proof
It is clear from (\ref{EQ44}) (\ref{EQ45}) that
\[
f_i Y^\mu V_\lambda = f_i V_{\lambda - n\mu} = Y_\mu f_i V_\lambda, \qquad
e_i Y^\mu V_\lambda = e_i V_{\lambda - n\mu} = Y_\mu e_i V_\lambda,
\]
that is, the action of $U_q(\slchap_n)$ commutes with
the operators $Y^\mu$.
Hence, recalling the discussion of Section~\ref{SECT5BIS},
we see that it is enough to prove that
$f_i T_j V_\gamma = T_j f_i V_\gamma$ for $\gamma\in \SG_r \A_{r,-n}$
and $1\le j \le r-1$.
Moreover, since $T_j$ only acts on components $j$ and $j+1$ of $\gamma$,
we can assume that $r=2$.
Then the claim is verified by a direct computation.
For example on the one hand
\[
f_0 T_1 \, V_{(1-n,0)} = f_0 \, V_{(0,1-n)} = V_{(1,1-n)},
\]
and on the other hand
\[
T_1 f_0 \,  V_{(1-n,0)} = q^{-1} T_1 V_{(1-n,1)}
= q^{-1} T_1 Y_2^{-1} V_{(1-n,1-n)}
\]
\[
\hskip 2cm
= q^{-1}(Y_1^{-1}T_1+(q-q^{-1})Y_1^{-1}) V_{(1-n,1-n)}
=V_{(1,1-n)}.
\]
\cqfd

\begin{remark}{\rm
This action of $U_q(\slchap_n)$ does not commute with the positive level
action $\Pi_n$ of $\EAHB_r$. For example if $r=2$ and $n=3$
\[
f_2\Pi_3(T_1)(V_{(2,0)}) = q^{-1}V_{(0,3)}, \qquad
\Pi_3(T_1)(f_2 V_{(2,0)}) = qV_{(0,3)}.
\]
However, one can easily obtain an action commuting with $\Pi_n$ by
simply replacing the comultiplication $\Delta$ of (\ref{COM})
by its opposite
\[
\Delta^\op f_i = f_i\otimes q^{h_i} + 1 \otimes f_i , \quad
\Delta^\op e_i = e_i\otimes 1 + q^{-h_i} \otimes e_i.
\]
\finex
}
\end{remark}
The action of $U_q(\slchap_n)$ is compatible with the bar involution
of $\PB_r$ in the following sense.
\begin{proposition}\label{BARF}
For $x\in U_q(\slchap_n)$ and $v\in \PB_r$ one has
$\bar{(xv)} = \bar{x}\,\,\bar{v}$.
In other words,
\[
\bar{f_i v} = f_i \bar{v},\qquad
\bar{e_i v} = e_i \bar{v}\qquad\qquad (0\le i \le n-1).
\]
\end{proposition}
\proof
We can assume that $v=V_\lambda$.
Then by (\ref{EQ51}) and Proposition~\ref{PROP5.4} we have
\begin{equation}\label{EQF1}
\bar{f_i V_\lambda} :=  q^{-\ell(w_{0,\xi})}\,
T_{w_0}^{-1}
\left(
\sum_{\stackrel{\scriptstyle j=1}{\lambda_j\equiv i}}^r
q^{\sum_{k=1}^{j-1}(\delta_{\lambda_k\equiv i}-\delta_{\lambda_k\equiv i+1})}
\,V_{w_0(\lambda+\epsilon_j)}
\right) .
\end{equation}
Here, $\xi\in\A_{r,-n}$ is the point congruent to $\lambda+\epsilon_j$,
which does not depend on $j$ because $\lambda_j$ is required to be
$\equiv i$.
On the other hand, since $f_i$ commutes with $T_{w_0}^{-1}$ by
Proposition~\ref{PROPCOMM},
\begin{eqnarray}
&&f_i \bar{V_\lambda} =
q^{-\ell(w_{0,\nu})}\,T_{w_0}^{-1}
\left(\sum_{\stackrel{\scriptstyle j=1}{\lambda_{r+1-j}\equiv i}}^r
q^{\sum_{k=1}^{j-1}(\delta_{\lambda_{r+1-k}\equiv i}-
\delta_{\lambda_{r+1-k}\equiv i+1})}
\,V_{w_0\lambda+\epsilon_j}\right)\qquad  \\
&&\hskip 0.9cm \label{EQF2}
= q^{-\ell(w_{0,\nu})}\,T_{w_0}^{-1}
\left(\sum_{\stackrel{\scriptstyle j=1}{\lambda_j\equiv i}}^r
q^{\sum_{k=1}^{r-j}(\delta_{\lambda_{r+1-k}\equiv i}-
\delta_{\lambda_{r+1-k}\equiv i+1})}
\,V_{w_0(\lambda+\epsilon_j)}\right).
\end{eqnarray}
It remains to check that the coefficients of
$T_{w_0}^{-1} V_{w_0(\lambda+\epsilon_j)}$
in (\ref{EQF1}) and (\ref{EQF2}) are equal, which is equivalent to
\[
\sum_{k=1}^r (\delta_{\lambda_k\equiv i}-\delta_{\lambda_k\equiv i+1})
- 1 = \ell(w_{0,\nu})-\ell(w_{0,\xi}).
\]
This is elementary, using for instance Remark~\ref{REM58}.
The formula for $e_i$ is similar and its proof is omitted.
\cqfd

\subsection{Action of $U_q(\widehat{\hbox{\twelvegoth sl}}_n)$
on $\FB_r$}
\label{SECT75}

Since the action of $U_q(\slchap_n)$ on $\PB_r$ commutes with
the action of $\EAHB_r$, the subspace
$\IB_r:=\Q(q)\otimes_{\Z[q,q^{-1}]}\I_r$ is stable under
$U_q(\slchap_n)$ and we obtain an induced action of
$U_q(\slchap_n)$ on $\FB_r$.
As explained in Section~\ref{SECT71}, the vector 
$V_\lambda$ should be regarded as a monomial tensor
$V_\lambda \equiv V_{\lambda_1}\otimes \cdots \otimes V_{\lambda_r}$.
Hence its projection $|\lambda\>$ on $\FB_r$
should be thought of as some $q$-wedge product 
$|\lambda\> \equiv V_{\lambda_1}\wedge_q \cdots \wedge_q V_{\lambda_r}$
with the anticommutation relations replaced by the straightening
rules of Proposition~\ref{PROP5.5}.
The action on $|\lambda\>$ of the generators of $U_q(\slchap_n)$ is obtained
by projecting (\ref{EQ44}), (\ref{EQ45}): 
\begin{eqnarray}
&&f_i |\lambda\> := \sum_{\stackrel{\scriptstyle j=1}{\lambda_j\equiv i}}^r
q^{\sum_{k=1}^{j-1}(\delta_{\lambda_k\equiv i}-\delta_{\lambda_k\equiv i+1})}
\,|\lambda+\epsilon_j\>, \hskip 1.7cm (0\le i \le n-1), \label{EQ51}\\
&&e_i|\lambda\> := \sum_{\stackrel{\scriptstyle j=1}{\lambda_j\equiv i+1}}^r
q^{-\sum_{k=j+1}^r(\delta_{\lambda_k\equiv i}-\delta_{\lambda_k\equiv i+1})}
\,|\lambda-\epsilon_j\>, \hskip 0.7cm (0\le i \le n-1).\label{EQ52}
\end{eqnarray}
Note that if $\lambda\in P^{++}$ then $\lambda \pm \epsilon_j \in P^+$.
It follows that either $|\lambda \pm \epsilon_j\>$
belongs to the basis $\{|\lambda\>\ |\ \lambda \in P^{++}\}$,
or $|\lambda \pm \epsilon_j\> = 0$.
Hence, Eq.~(\ref{EQ51}) (\ref{EQ52}) require no straightening relation and
are very simple to use in practice.
The compatibility of the bar involution with this action is
given by the next
\begin{proposition}\label{PROPBARFOCK}
For $u\in \FB_r$ and $0\le i \le n-1$ one has
\[
\bar{f_i u} = f_i \bar{u},\qquad
\bar{e_i u} = e_i \bar{u}.
\]
\end{proposition}
\proof
This follows immediately from (\ref{EQINV}) and Proposition~\ref{BARF}.
\cqfd

\subsection{The Fock space $\FI$}

For $s\ge r$ define a linear map $\varphi_{r,s} : \FB_r \longrightarrow \FB_s$
by
\[
\varphi_{r,s}(|\lambda\>) := 
|\lambda_1,\ldots,\lambda_r,-r,-r-1,\ldots,-s+1\> \qquad
(\lambda\in P^+_r).
\]
Then clearly $\varphi_{s,t} \circ \varphi_{r,s} = \varphi_{r,t}$.
Let $\FI := {\displaystyle\lim_\rightarrow}\, \FB_r$ 
be the direct limit of the
vector spaces $\FB_r$ with respect to the maps $\varphi_{r,s}$.
Each $|\lambda\>$ in $\FB_r$ has an image $\varphi_r(|\lambda\>)\in \FI$,
which should be thought of as some infinite $q$-wedge
\[
\varphi_r(|\lambda\>)\equiv V_{\lambda_1}\wedge_q \cdots
\wedge_q V_{\lambda_r}\wedge_qV_{-r}\wedge_q V_{-r-1} \wedge_q \cdots
\]
\begin{lemma}\label{LEMBASE2}
{\rm (i)} If $\lambda_r \le -r$ then $\varphi_r(|\lambda\>) = 0$.

\noindent
{\rm (ii)} If $\lambda\in P^{++}_r$ and $\lambda_r>-r$ then 
$\varphi_r(|\lambda\>) \not = 0$.
\end{lemma}
\proof
(i) Write $\lambda_r = k \le -r$ and consider the element
\[
\varphi_{r,-k+1}(|\lambda\>) = 
|\lambda_1,\ldots,\lambda_r,-r,-r-1,\ldots,k\>. 
\]
By applying Proposition~\ref{PROP5.5} one checks easily that
$
|k,-r,-r-1,\ldots,k\> 
$
straightens to 0 in $F_{-k-r+2}$.  
Therefore $\varphi_{r,-k+1}(|\lambda\>) = 0$, 
hence $\varphi_r(|\lambda\>) = 0$.
(ii) By Lemma~\ref{LEMBASE} if $\lambda \in P^{++}_r$ and
$\lambda_r > -r$ then
$\varphi_{r,s}(|\lambda\>) \not = 0$ for all $s>r$.
Hence $\varphi_r(|\lambda\>) \not = 0$.
\cqfd

\noindent
Let $\PP^+$ denote the set of all partitions, \ie of all finite
non-increasing sequences of positive integers.
Put $\rho_r^*:=(0,-1,\ldots ,-r+1)$, and for 
$\alpha = (\alpha_1,\ldots ,\alpha_s)\in \PP^+$ define
\[
|\alpha) := \varphi_s(|\alpha+\rho_s^*\>).
\]
It readily follows from Lemma~\ref{LEMBASE} and Lemma~\ref{LEMBASE2}
that $\{ |\alpha)\,|\,\alpha\in \PP^+\}$ is a basis of~$\FI$.
We define a grading on $\FI$ by requiring that 
\[
\deg |\alpha) := \sum_{i=1}^s \alpha_i.
\]
Then for all $\lambda\in P_r$, $\varphi_r(|\lambda\>)$ is
homogeneous of degree
\[
\deg \varphi_r(|\lambda\>) = \sum_{i=1}^r (\lambda_i+i-1).
\] 
In particular, if $\sum_{i=1}^r (\lambda_i+i-1) < 0$ then
$\varphi_r(|\lambda\>) = 0$.
\subsection{Action of $U_q(\widehat{\hbox{\twelvegoth sl}}_n)$ on $\FI$}

Let $\lambda \in \P_r$. 
It follows easily from (\ref{EQ51}) that 
\begin{equation}\label{STAB}
f_i \, \varphi_{r,s} (|\lambda\>)
= \varphi_{r+1,s} \, f_i \, \varphi_{r,r+1}\, (|\lambda\>)
\end{equation}
for all $s>r$.
Hence one can define an endomorphism $f_i$ of $\FI$ by
\begin{equation}\label{ACTFDEF}
f_i \, \varphi_r (|\lambda\>) =  \varphi_{r+1} \, f_i \, \varphi_{r,r+1}\, (|\lambda\>)
\end{equation}
and thus get an action of $U_q^-(\slchap_n)$ on $\FI$. 

On the basis $\{|\alpha)\ | \ \alpha\in \PP^+\}$ this action is expressed
as follows.
Let $\alpha$ and $\beta$ be two Young diagrams 
such that $\beta$ is obtained from $\alpha$ by adding a  cell $\gamma$ 
whose content is $\equiv i \mod n$.
Such a cell is called a removable $i$-cell of $\beta$, or
an indent $i$-cell of $\alpha$.
Let $I_i^r(\alpha,\beta)$ (\resp $R_i^r(\alpha,\beta)$)
be the number of indent $i$-cells of $\alpha$ 
(\resp of removable $i$-cells of $\alpha$)
situated to the right of $\gamma$ ($\gamma$ not included).
Set
$N_i^r(\alpha,\beta)=I_i^r(\alpha,\beta)-R_i^r(\alpha,\beta)$. Then
Eq.~(\ref{EQ51}) gives
\begin{equation}\label{ACTF}
f_i |\alpha ) = \sum_\beta q^{N_i^r(\alpha,\beta)} |\beta), 
\end{equation}
where the sum is over all partitions $\beta$ such that
$\beta/\alpha$ is an $i$-cell.

Defining an action of $U^+_q(\slchap_n)$ is not as straightforward since there is no 
formula like (\ref{STAB}) for $e_i$.
For example if $n=2$, 
\begin{eqnarray*}
&& e_1\,|2\> = |1\> , \\
&& e_1\,|2,-1\> = q^{-1}\,|1,-1\> , \\
&& e_1\,|2,-1,-2\> = |1,-1,-2\> + |2,-1,-3\>, \\
&& e_1\,|2,-1,-2,-3\> = q^{-1}\,|1,-1,-2,-3\> , 
\end{eqnarray*}
and in general
\begin{eqnarray*}
&& e_1\varphi_{1,2r}\,|2\> = q^{-1}\varphi_{1,2r} e_1 |2\>, \\
&& e_1\varphi_{1,2r+1}\,|2\> = \varphi_{1,2r+1} e_1 |2\> + |2,-1,\ldots ,-2r+1,-2r-1\>. 
\end{eqnarray*}
However, one can check that putting 
$e_i \varphi_r (|\lambda\>) := 
q^{-\delta_{i\equiv r}}\,\varphi_r (e_i |\lambda\>)$
one gets a well-defined action of $U^+_q(\slchap_n)$ compatible with 
(\ref{ACTFDEF}) (see \cite{KMS}).
Its combinatorial description is given by
\begin{equation}\label{ACTE}
e_i |\beta ) = \sum_\alpha q^{-N_i^l(\alpha,\beta)} |\alpha), 
\end{equation}
where the sum is over all partitions $\alpha$ such that
$\beta/\alpha$ is an $i$-cell, and $N_i^l(\alpha,\beta)$ is defined
as $N_i^r(\alpha,\beta)$ but replacing right by left.

In contrast to $\FB_r$, the representation $\FI$ has primitive 
vectors, \ie vectors annihilated by all $e_i$. In particular 
the vector $|0)$ labelled by the unique partition of 0 is
primitive. 
In fact $\FI$ is a level 1 integrable representation of 
$U_q(\slchap_n)$, while $\FB_r$ is a level 0 non-integrable
representation.
As shown by Kashiwara, Miwa and Stern \cite{KMS}, 
the decomposition of $\FI$ into simple $U_q(\slchap_n)$-modules
is obtained by considering the limit $r\rightarrow \infty$
of the action of $Z(\EAHB_r)$ on $\FB_r$.

\subsection{Action of $\HI$ on $\FI$}

Let $\lambda \in P_r$. It follows from the easily checked relations
\begin{equation}\label{EQREDRESS}
|-s,-r,-r-1,\ldots,-s\> = 0,\quad 
|-r,-r-1,\ldots,-s,-r\> =0, \quad (s\ge r\ge 0)
\end{equation}
that the vector $\varphi_s \, B_k \, \varphi_{r,s}(|\lambda\>)$
is independent of $s$ for $s> r$ large enough.
Hence one can define 
endomorphisms $B_k$ of $\FI$ by 
\begin{equation}\label{ACTBk}
B_k \varphi_r (|\lambda\>) := \varphi_s \, B_k \, \varphi_{r,s}(|\lambda\>) 
\qquad (k\in\Z^*,\ s\gg 1).
\end{equation}
By construction, these endomorphisms commute with the action of 
$U_q(\slchap_n)$ on $\FI$.
However they no longer generate a commutative
algebra but a Heisenberg algebra. 
Indeed, it was proved by Kashiwara, Miwa and Stern \cite{KMS}
that 
\begin{equation}\label{EQCOMMUT}
[B_k,B_l] = 
\left\{
\begin{array}{ll}\displaystyle
k\,{1-q^{-2nk}\over 1-q^{-2k}} &\mbox{if $k=-l$,}\\[1mm]
0 & \mbox{otherwise.}
\end{array}
\right.
\end{equation}
We shall denote this Heisenberg algebra by $\HI$.
The elements $\U_\beta,\ \V_\beta,\ 
\UU_\beta,\ \VV_\beta$ of $Z(\EAH_r)$ also give 
rise to well-defined elements of $\HI$ that we still denote 
by $\U_\beta,\ \V_\beta,\ \UU_\beta,\ \VV_\beta$.
By Theorem~\ref{THRIB} and Theorem~\ref{THRIB2},
their action on the basis $\{|\nu),\ \nu\in \PP^+\}$ of  $\FI$ is given by
\begin{eqnarray}
&&\U_\beta \, |\nu) = q^{-(n-1)k}\sum_{\mu\in\PP^+}
L_{\nu'/\mu'\,,\,\beta}^{(n)}(-q)\label{ACTUT}
\, \, |\mu),\\
&&\V_\beta \, |\nu) = q^{-(n-1)k}\sum_{\lambda\in\PP^+}
L_{\lambda'/\nu'\,,\,\beta}^{(n)}(-q)\label{ACTVT}
\,\, |\lambda),\\
&&\UU_\beta \, |\nu) = \sum_{\mu\in\PP^+}
L_{\nu/\mu\,,\,\beta}^{(n)}(-q^{-1})\label{ACTU}
\, \, |\mu), \\
&&\VV_\beta \, |\nu) = \sum_{\lambda\in\PP^+}
L_{\lambda/\nu\,,\,\beta}^{(n)}(-q^{-1})
\,\, |\lambda),\label{ACTV}
\end{eqnarray}
where $k=|\beta|$.
It was shown in \cite{KMS} that $\FI$ is irreducible under the
commuting actions of $U_q(\slchap_n)$ and $\HI$. 
It follows that $\{\VV_\beta \,|0) , \,  \beta \in \PP^+\}$
is a basis of the space of primitive vectors of $\FI$
for $U_q(\slchap_n)$. 

\subsection{The bar involution of $\FI$}

Before introducing the involution we need the following
simple lemmas.
\begin{lemma}\label{LEMREDR}
Let $\mu\in P_{m+1}$ such that
$\mu_i > -m \ (i=1,\ldots , m+1)$ and 
$\sum_i(\mu_i+i-1) \le m$.
Then $|\mu\> = 0$.
\end{lemma}
\proof
We have $$|\mu\> = \sum_{\lambda\in P_{m+1}^{++}} x_\lambda\, |\lambda\>$$
for some coefficients $x_\lambda$.
Because of the hypothesis $\mu_i > -m $ and of the form of the
straightening relations, the components of the weights 
$\lambda$ occuring in this sum must all be $>-m$. 
On the other hand, setting $\alpha_i=\lambda_i+i-1$ we see that
$\alpha$ is a partition with $|\alpha| \le m$, hence $l(\alpha)\le m$.
Thus the last component of all the $\lambda$ must be $= -m$, 
and the sum is empty. \cqfd

\begin{lemma}\label{LEMINV}
Let $\lambda \in P_r$ and let $m\ge r$.
Assume that $\lambda_i > -m\ (i=1,\ldots ,r)$ and
$\sum_i(\lambda_i+i-1) \le m$.
Then 
\[
|-m,\lambda_1,\ldots,\lambda_r,-r,\ldots ,-m+1\>
=
(-1)^m q^{-a(\lambda)}\,
|\lambda_1,\ldots,\lambda_r,-r,\ldots ,-m+1,-m\>,
\] 
where 
$a(\lambda)=\sharp\{j \le r\ |\ \lambda_j \not \equiv -m\}
+ \sharp\{-r\ge j \ge -m+1\ |\ j \not \equiv -m\}$.
\end{lemma}
\proof
Consider the straightening of 
$\nu = |-m,\lambda_1,\ldots,\lambda_r,-r,\ldots ,-m+1\>$
computed by means of Proposition~\ref{PROP5.5}. 
At each step, if the third rule (\ref{STR3}) has to be used, 
then only the first term of the right-hand side may be non-zero.
Indeed the two other terms involve weights $\mu$ which satisfy 
the hypothesis of Lemma~\ref{LEMREDR}. 
Therefore the straightening of $\nu$ is simply obtained by
reordering its components and multiplying by the appropriate sign
and power of $q$.
\cqfd

\noindent
If $\lambda$ satisfies the hypothesis of Lemma~\ref{LEMINV}, 
then repeated applications of this lemma show that for $p\ge m$,
\begin{eqnarray*}
&&(-1)^{p\choose 2} q^{b(\lambda,p)}\,
|-p,\ldots , -r,\lambda_r,\ldots ,\lambda_1\>
=\\
&& \hskip 3cm (-1)^{m\choose 2} q^{b(\lambda,m)}\,
|-m,\ldots , -r,\lambda_r,\ldots ,\lambda_1,-m-1,\ldots ,-p\>
\end{eqnarray*}
Here $b(\lambda,p)$ is the number of pairs $(i,j)$ of 
components of $(\lambda_1,\ldots ,\lambda_r,-r,\ldots , -p)$
with $i\not \equiv j \mod n$. 
In other words, using Proposition~\ref{PROP5.4} and Remark~\ref{REM58}
\[
\bar{\varphi_{r,p}(|\lambda\>)}
= 
 \varphi_{m,p}(\bar{\varphi_{r,m}(|\lambda\>)}).
\]
Thus we can define a semi-linear involution on $\FI$ by putting
\begin{equation}\label{DEFINVOLU}
\bar{\varphi_r(|\lambda\>)} := \varphi_m(\bar{\varphi_{r,m}(|\lambda\>)})
\quad
(\lambda\in P_r,\ \deg \varphi_r(|\lambda\>) = m,\ \lambda_i > -m).
\end{equation}
In particular, for $\alpha\in \PP^+$ and $s\ge |\alpha|$, we have  
$ \bar{|\alpha)} = \varphi_s(\bar{|\alpha+\rho^*_s\>}) $.
\begin{proposition}\label{INVOLU}
For $\alpha\in \PP^+$, $0\le i \le n-1$ and $k\in\N^*$ we have
\[
\bar{f_i\,|\alpha)} = f_i \,\bar{|\alpha)},\quad
\bar{e_i\,|\alpha)} = e_i \,\bar{|\alpha)},\quad
\bar{B_{-k}\,|\alpha)} = B_{-k}\, \bar{|\alpha)},\quad
\bar{B_k\,|\alpha)} = q^{2(n-1)k} B_k\, \bar{|\alpha)}.
\]
\end{proposition}
\proof For $f_i$ and $B_{-k}$, the proof readily follows from
Proposition~\ref{PROPBARFOCK}, Proposition~\ref{PROPBARFOCKB}  and
(\ref{ACTFDEF}) (\ref{ACTBk}) (\ref{DEFINVOLU}).
(Note that the condition $\lambda_i > -m$ in (\ref{DEFINVOLU})
is preserved by the action of these lowering operators.)
Let us prove the statement for $B_k$.
We argue by induction on $\deg |\alpha)$. In degree 0, the unique
basis vector is $|0)$ and we have $B_k|0) = B_k\bar{|0)} = 0$,
so the claim is trivially verified.
Let us assume that the result is proved  
for all $|\alpha)$ of degree $\le m$.
Since the action of the operators $B_{-l}$ and $f_i$ on $|0)$ generates
the whole Fock space, it is enough to prove that 
\[
\bar{B_k\,f_i v} = q^{2(n-1)k} B_k\, \bar{f_i v}, \qquad
\bar{B_k\,B_{-l} v} = q^{2(n-1)k} B_k\, \bar{B_{-l} v}
\]
for all $v$ of degree $\le m$.
Now $B_k$ and $f_i$ commute, so 
\[
\bar{B_k f_i v} = \bar{f_i B_k v} = f_i \bar{B_k v}
= f_i (q^{2(n-1)k} B_k \bar{v}) = q^{2(n-1)k} B_k f_i \bar{v}
= q^{2(n-1)k} B_k \bar{f_i v}.
\]
If $l \not = k$ we know that $B_k$ and $B_{-l}$ commute and we can
argue similarly.
Finally if $l=k$, by (\ref{EQCOMMUT}),
\begin{eqnarray*}
&&\bar{B_k\,B_{-k} v} = \bar{B_{-k}\,B_k v} + k\,{1-q^{2(n-1)k}\over 1-q^{2k}}
\bar{v}
= q^{2(n-1)k}B_{-k}B_k \bar{v} + k\,{1-q^{2(n-1)k}\over 1-q^{2k}}\bar{v}
\\
&&
\hskip 0.5cm 
= q^{2(n-1)k}\left(B_{-k}B_k + k\,{1-q^{-2(n-1)k}\over 1-q^{-2k}}\right)\bar{v}
= q^{2(n-1)k} B_k B_{-k} \bar{v}
= q^{2(n-1)k} B_k \bar{B_{-k} v}.
\end{eqnarray*}
The proof for $e_i$ is similar, using the commutation relation
\[
[e_i,f_j] = \delta_{ij}{q^{h_i}-q^{-h_i}\over q-q^{-1}}.
\] 
\cqfd

\noindent
Proposition~\ref{INVOLU} implies that for $|\beta| = k$,
\begin{eqnarray}
&&\bar{\V_\beta\,|\alpha)} = \V_\beta\, \bar{|\alpha)},\hskip 2cm \label{BARV}
\bar{\VV_\beta\,|\alpha)} = \VV_\beta\, \bar{|\alpha)},\\
&&\bar{\U_\beta\,|\alpha)} = q^{2(n-1)k}\,\U_\beta\, \bar{|\alpha)},\hskip 0.7cm
\bar{\UU_\beta\,|\alpha)} = q^{2(n-1)k}\,\UU_\beta\, \bar{|\alpha)}.\hskip 1cm
\end{eqnarray}
\subsection{The scalar product of $\FI$}

Define a scalar product on $\FI$ by
$\< |\alpha) \,,\,|\beta) \> = \delta_{\alpha\beta}$.
\begin{proposition}\label{SCALAR}
For $u,v \in \FI$ one has
\begin{eqnarray*}
&&\<f_i u \,,\, v \> = \<u \,,\, q^{h_i-1} e_i v\>, \quad
\<e_i u \,,\, v \> = \<u \,,\, q^{-h_i-1} f_i v\>,\\ 
&&\<\V_\alpha u \,,\, v \> = \<u \,,\, \U_\alpha v\>, \hskip 1cm
\<\VV_\alpha u \,,\, v \> = \<u \,,\, \UU_\alpha v\>.
\end{eqnarray*}
\end{proposition}
\proof This follows immediately from (\ref{ACTF}) (\ref{ACTE})
(\ref{ACTUT}) (\ref{ACTVT}) (\ref{ACTU}) (\ref{ACTV}).

\subsection{Symmetry of the bar involution}

Define a semi-linear involution $v \mapsto v'$ on $\FI$ by
$|\alpha)':=|\alpha')$, where $\alpha'$ is the partition conjugate
to $\alpha \in \PP^+$.

\begin{proposition}\label{PROPPRIM}
For $u\in\FI$ and $\beta\in\PP^+$ with $|\beta|=k$, there holds
\[
\begin{array}{ll}
(e_iu)' = q^{h_{-i}-1}e_{-i}u', &
(f_iu)' = q^{-h_{-i}-1}f_{-i}u',\\[2mm]
(\VV_\beta u)'=(-q)^{(n-1)k}\,\V_\beta u',\hskip 1cm&
(\UU_\beta u)'=(-q)^{(n-1)k}\,\U_\beta u'.
\end{array}
\]
\end{proposition}

\proof This also follows from (\ref{ACTF}) (\ref{ACTE})
(\ref{ACTUT}) (\ref{ACTVT}) (\ref{ACTU}) (\ref{ACTV}).
\cqfd

\noindent
Let $S_\beta = \sum_\alpha \kappa_{\beta,\alpha}\VV_\alpha$
be the element of $\HI$ corresponding to the Schur function
$s_\beta$. 
The third equation above implies that 
\begin{equation}\label{EQSYMS}
(S_\beta u)' = (-q)^{(n-1)k}\, S_{\beta'} u'.
\end{equation}
\begin{theorem}\label{SYMINVOL}
For $u,v \in \FI$ we have
\[
\<\bar{u}\,,\,v\> = \<u'\,,\,\bar{v'}\>\,.
\]
\end{theorem}
\proof
The proof is by induction on the degree $d$ of $u$ and $v$.
If $d=0$ this is clear. So let us assume that the
theorem is proved in degree $d < m$.
The operators $e_i$, $f_i$, $\U_k$, $\V_k$, $\UU_k$, $\VV_k$ are homogeneous
of respective degree $-1$, $+1$, $-kn$, $+kn$, $-kn$, $+kn$.
Since $\FI$ is generated by the action of the operators
$f_i$ and $\VV_k$ on the highest weight vector $|0)$, 
it is enough to prove that 
\begin{eqnarray}
&&\<\overline{(f_iu)},v\> = \<(f_iu)',\overline{v'}\>,\label{2PR1}
\\
&&\<\overline{(\VV_kw)},v\> = \<(\VV_kw)',\overline{v'}\>,\label{2PR2}
\end{eqnarray} 
for all $u,\ v,\ w$ with $\deg u = m-1$, $\deg v = m$,
$\deg w = m-kn$.

Let us prove (\ref{2PR1}).
We have
$$
\<\bar{(f_iu)},v\> = \<f_i\,\bar{u},v\> = 
\<\bar{u}, q^{h_i-1} e_i v\>
= \<u', \bar{(q^{h_i-1} e_i v)'}\>.
$$
The first equality comes from Proposition~\ref{INVOLU},
the second from Proposition~\ref{SCALAR} and the third
from the fact that $\deg u < m$.
Now, by Proposition~\ref{INVOLU},~\ref{SCALAR} and \ref{PROPPRIM}
\[
\<u', \bar{(q^{h_i-1} e_i v)'}\>
= \<u', \bar{e_{-i}v'}\>
=\<u',e_{-i}\bar{v'}\>
=\<q^{-h_{-i}-1}f_{-i}u',\bar{v'}\>
=\<(f_iu)',\bar{v'}\>,
\]
and (\ref{2PR1}) is proved.

The proof of (\ref{2PR2}) is similar.
We have
$$
\<\overline{(\VV_kw)},v\> = \<\VV_k\overline{w},v\>
=\<\overline{w},\UU_k v\> = \< w' , \overline{(\UU_kv)'}\>.
$$
The first equality comes from (\ref{BARV}),
the second from Proposition~\ref{SCALAR} and the third
from the fact that $\deg w < m$.
Then, using again  Proposition~\ref{INVOLU},~\ref{SCALAR} and \ref{PROPPRIM},
\begin{eqnarray*}
&& \< w' , \overline{(\UU_kv)'}\> = \< w' , \overline{(-q)^{(n-1)k}\,\U_k(u')}\>
= \< w' , (-q)^{(n-1)k}\,\U_k(\overline{u'})\>\\[2mm]
&&\hskip 1.8cm = \< (-q)^{(n-1)k}\V_k w' , \overline{v'}\>
= \<(\VV_kv)',\overline{w'}\>,
\end{eqnarray*}
and (\ref{2PR2}) is proved. 
\cqfd

\subsection{Canonical bases of $\FI$}

For $\beta \in \PP^+$ write 
$
\bar{|\beta)} = \sum_{\alpha\in \PP^+} b_{\alpha,\beta}(q)\,|\alpha).
$
Then, for $|\alpha|= |\beta| \le r$ it follows from (\ref{DEFINVOLU}) 
that we have 
\[
b_{\alpha,\beta}(q) = 
a_{\alpha+\rho_r^*,\,\beta+\rho_r^*}(q)
= a_{\alpha+\rho_r,\,\beta+\rho_r}(q),
\]
where the coefficients $a_{\lambda,\mu}(q)\ (\lambda,\ \mu\in P^{++}_r)$ 
have been defined in Section~\ref{SECT53}.
Hence by Corollary~\ref{COR} the matrix 
\[{\bf B}_k := [b_{\alpha,\beta}(q)], \qquad (\alpha,\beta\vdash k)\]
is unitriangular, and one can define canonical bases 
$\{\G^+_\alpha \ | \ \alpha \in \PP^+\}$,
$\{\G^-_\alpha\ | \ \alpha \in \PP^+\}$ of $\FI$
characterized by:
\begin{quote}
{\rm (i)} \quad $\overline{\G^+_\alpha}=\G^+_\alpha$, \quad
$\overline{\G^-_\alpha}=\G^-_\alpha$,

{\rm (ii)} \quad $\G^+_\alpha \equiv |\alpha) \mod q\L^+_\infty$,
\quad $\G^-_\alpha \equiv |\alpha) \mod q^{-1}\L^-_\infty$,
\end{quote}
where $\L^+_\infty$ (\resp $\L^-_\infty$) is the $\Z[q]$-submodule
(\resp $\Z[q^{-1}]$-submodule) spanned by the vectors $|\alpha)$.
Set
\[
\G^+_\beta =\sum_\alpha d_{\alpha,\beta}(q)\,|\alpha)\,,\quad\quad
\G^-_\alpha =\sum_\beta e_{\alpha,\beta}(-q^{-1})\, |\beta)\,,
\]
and 
\[
{\bf D}_k:=[d_{\alpha,\beta}(q)], \qquad 
{\bf E}_k:=[e_{\alpha,\beta}(q)], \qquad (\alpha,\beta\vdash k). 
\]
Then, for $r\ge k$ we have
\[
d_{\alpha,\beta}(q) = c_{\alpha+\rho_r,\,\beta+\rho_r}(q),\qquad
e_{\alpha,\beta}(q) = l_{\alpha+\rho_r,\,\beta+\rho_r}(q).
\]
Hence by Theorem~\ref{THVV} we get
\begin{equation}\label{EQUAE}
e_{\alpha,\beta} = P^-_{\beta+\rho_r,\alpha+\rho_r},
\end{equation}
a parabolic Kazhdan-Lusztig polynomial for $\EASG_r$
associated with the parabolic subgroup $\SG_{\nu,-n}$
which stabilizes the point $\nu \in \A_{r,-n}$ congruent to $\alpha+\rho_r$
and $\beta+\rho_r$. 
Also, putting 
$\widehat{u}_\alpha  := w(w_0(\alpha+\rho_r),-n)\,w_{0,\nu}$
and
$\widehat{u}_\beta  := w(w_0(\beta+\rho_r),-n)\,w_{0,\nu}$,
we have
\begin{equation}\label{EQUAD}
d_{\alpha,\beta} = \sum_{s\in\SG_r}(-q)^{\ell(s)}
P_{s \widehat{u}_{\alpha}\, , \, \widehat{u}_{\beta}}.
\end{equation}
Note that by Theorem~\ref{THDEO} this is also a parabolic Kazhdan-Lusztig polynomial
of negative type associated with the subgroup $\SG_r \subset \EASG_r$.
It is interesting to give another expression of $d_{\alpha,\beta}$
in terms of the action $\pi_n$ (instead of $\pi_{-n}$).
Let $\underline{P}_r = P_r/\Z(1,\ldots ,1)$ and $\lambda \mapsto \underline{\lambda}$
be the natural projection $P_r\rightarrow \underline{P}_r$.
The action $\pi_n$ of $\EASG_r$ on $P_r$ induces an action
$\underline{\pi}_n$ of $\ASG_r$ on $\underline{P}_r$
with fundamental alcove
$\underline{\A}_{r,n}:=\{\underline{\lambda}\in\underline{P}_r \ |\ 
\lambda_1\ge \cdots \ge\lambda_r,\ \lambda_1-\lambda_r \le n\}$.
Let $\xi$ be the point of $\underline{\A}_{r,n}$ congruent to 
$\underline{\alpha+\rho_r}$
and $\underline{\beta+\rho_r}$ under $\underline{\pi}_n$,
and let $w_{0,\xi}$ denote the longest element of its stabilizer.
Consider the projection $\underline{\cdot}:\ \EASG_r \rightarrow \ASG_r$
defined by $\underline{\sigma\tau^k}=\sigma \ (k\in\Z,\, \sigma\in\ASG_r)$,
and the automorphism $\sharp$ of $\EASG_r$ defined by
$s_i^\sharp = s_{-i} \ (i\in \Z/r\Z)$.
It is easy to check that, for $\lambda\in P_r^+$, 
$\underline{w(w_0\lambda,-n))}=\underline{(w(\lambda,n))}^\sharp$.
It follows that 
\begin{equation}\label{EQUAD2}
d_{\alpha,\beta} = \sum_{s\in\SG_r}(-q)^{\ell(s)}
P_{s \widehat{v}_{\alpha}\, , \, \widehat{v}_{\beta}},
\end{equation}
where $\widehat{v}_{\alpha},\,\widehat{v}_{\beta}$ are given
by $\widehat{v}_{\alpha} = \underline{w(\alpha+\rho_r,n)}w_{0,\xi},\,
\widehat{v}_{\beta} = \underline{w(\beta+\rho_r,n)}w_{0,\xi}$.

\begin{remark}{\rm Consider the $U_q(\slchap_n)$-submodule $M$ of $\FI$
generated by $|0)$.
This is an irreducible integrable representation with
highest weight $\Lambda_0$.
By Proposition~\ref{INVOLU}, the bar involution of $\FI$ induces the 
Kashiwara involution of $M$, and it follows that the subset 
$\{\G^+_\alpha \ | \ \alpha \mbox{ is $n$-regular } \}$
is the global lower crystal basis of $M$ (see \cite{LLT}).
The expressions (\ref{EQUAD}) and (\ref{EQUAD2}) of the coefficients of this
basis as Kazhdan-Lusztig polynomials have been obtained
independently by Vasserot, Varagnolo \cite{VV} and
by Goodman, Wenzl \cite{GW} respectively. \finex
}
\end{remark}
It follows from Theorem~\ref{THSTPR} that the basis $\G^-_\alpha$ satisfies
the following analogue of the Steinberg-Lusztig tensor product theorem.
Let $\alpha \in \PP^+$ of length $r$. Write $\alpha = \alpha^{(0)} + n\alpha^{(1)}$, where
$\alpha^{(0)}$ is $n$-restricted, that is,
\[
0\le \alpha_i^{(0)}-\alpha_{i+1}^{(0)} <n \qquad (1\le i \le r-1).
\]
Then
$
\G^-_{\alpha} = S_{\alpha^{(1)}}\,\G^{-}_{\alpha^{(0)}}.
$
Taking $\alpha^{(0)}  = (0)$ and writing $n\alpha$ in place of $\alpha$
we obtain that 
\begin{equation}\label{EQALPHA}
\G^-_{n\alpha} = S_{\alpha}\,|0).
\end{equation}
We can now prove the following symmetry of the basis $\{\G^-_\alpha\}$.
\begin{theorem}\label{THCONJ}
Let $\lambda,\,\mu^{(0)}, \ldots , \mu^{(n-1)}$ be partitions.
Set $k=|\lambda|$. There holds
\[
\begin{array}{ll}
\mbox{\rm (i) } &(\G^-_{n\lambda})' = (-q)^{(n-1)k}\,\G^-_{n\lambda'}\,,\\[2mm]
\mbox{\rm (ii)}\quad &c_{\mu^{(0)},\, \ldots \,,\, \mu^{(n-1)}}^\lambda(q^{-1})=q^{-(n-1)k}\,
c_{({\mu^{(n-1)})}',\, \ldots \,,\, {(\mu^{(0)})}'}^{\lambda'}(q).
\end{array}
\]
\end{theorem}
\proof By (\ref{EQALPHA}) and (\ref{EQSYMS}) we have
\[
(\G^-_{n\lambda})' = \left(S_{\lambda}\,|0)\right)'
= (-q)^{(n-1)k}\,S_{\lambda'}\,|0)
= (-q)^{(n-1)k}\,\G^-_{n\lambda'}.
\]
The second equation follows now from the fact that
if $\mu$ is the partition with $n$-quotient 
$(\mu^{(0)},\, \ldots \,,\, \mu^{(n-1)})$
then $\mu'$ has $n$-quotient 
$(({\mu^{(n-1)})}',\, \ldots \,,\, {(\mu^{(0)})}')$.
\cqfd

Let $\{\G^*_\alpha\}$ denote the basis of $\FI$
adjoint to $\{\G^+_\alpha\}$ for the above scalar product.
In other words, 
$\<\G^*_\alpha\,,\,\G^+_\beta\> = \delta_{\alpha,\beta}.$
Write 
\[
\G^*_\alpha =\sum_\beta g_{\alpha,\beta}(q)\, |\beta),
\qquad \mbox{and} \qquad
{\bf G}_k:=[g_{\alpha,\beta}(q)], \qquad (\alpha,\beta\vdash k).
\]
Since $\{|\alpha)\}$ is an orthonormal basis, we have 
${\bf G}_k = {\bf D}_k^{-1}$.
\begin{theorem}
For $\alpha \in \PP^+$ one has $(\G^*_\alpha)' = \G^-_{\alpha'}.$
\end{theorem}
\proof
We have to prove that $(\G^*_\alpha)'$ satisfies the two defining
properties of $\G^-_{\alpha'}$, namely
\[
(\G^*_\alpha)'\equiv \alpha' \mod q^{-1}\L^-_\infty, \qquad
\bar{(\G^*_\alpha)'} = (\G^*_\alpha)'.
\]
The first property is obvious. Indeed by definition
$\G^+_\alpha \equiv |\alpha) \mod q\L^+_\infty$.
Since ${\bf G}_k = {\bf D}_k^{-1}$, we deduce that
$\G^*_\alpha\equiv |\alpha) \mod q\L^+_\infty$, which implies that
$(\G^*_\alpha)'\equiv \alpha' \mod q^{-1}\L^-_\infty$.
The second property is equivalent to
\[
\<\bar{(\G^*_\alpha)'}\,,\,(\G^+_\beta)'\> = \delta_{\alpha,\beta},\qquad
(\alpha,\ \beta \vdash k ),
\]
because $\{(\G^*_\alpha)'\}$ is the basis adjoint to 
$\{(\G^+_\beta)'\}$.
Now, by Theorem~\ref{SYMINVOL},
\[
\<\bar{(\G^*_\alpha)'}\,,\,(\G^+_\beta)'\> = 
\<\G^*_\alpha\,,\,\bar{\G^+_\beta}\>
= \<\G^*_\alpha\,,\,\G^+_\beta\>
= \delta_{\alpha,\beta}.
\]
\cqfd

\begin{corollary}\label{CORINV}
Let ${\bf J}_k = [j_{\alpha,\beta}(q)]_{\alpha,\beta\vdash k}
:=[e_{\alpha',\beta'}(-q)]^{-1}_{\alpha,\beta\vdash k}$.
Then ${\bf J}_k = {\bf D}_k$.
In other words, we have
\[
\sum_{\gamma\vdash k} e_{\alpha',\gamma'}(-q)\, d_{\gamma,\beta}(q) 
= \delta_{\alpha,\beta},
\]
where $e_{\alpha',\gamma'}$ and $d_{\gamma,\beta}$ are the parabolic 
Kazhdan-Lusztig
polynomials given by {\rm(\ref{EQUAE}) (\ref{EQUAD})}.
\end{corollary}

\begin{remark}\label{REMTILT}{\rm
(i) Let $\alpha,\ \beta$ be two partitions of $k$ and take $r\ge k$.
Assuming Lusztig's conjecture (\ref{EQ28}), it follows from 
Corollary~\ref{CORINV} that 
$$d_{\alpha,\beta}(1) = j_{\alpha,\beta}(1) = [W(\alpha') : L(\beta')],$$
the multiplicity of the simple $U_\zeta(\gl_r)$-module $L(\beta')$ in  
the Weyl module $W(\alpha')$, 
as was conjectured in \cite{LT}, Conjecture~5.2.

\smallskip\noindent
(ii) For $\lambda \in P^+_r$, let $T(\lambda)$ denote the indecomposable tilting
$U_\zeta(\gl_r)$-module with highest weight $\lambda$.
By Proposition~8.2 of \cite{DPS} which states that
$$[W(\alpha') : L(\beta')] = [T(\beta) : W(\alpha)],$$
we see that $[T(\beta) : W(\alpha)] = d_{\alpha,\beta}(1)$.
Taking into account (\ref{EQUAD2}) we thus get another proof
of the character formula of Soergel \cite{So2} in type $A$.
Note that we do not need to deduce the formula for singular
weights from that for regular weights (see \cite{So1}, Remark 7.2).
In particular, we see that the formula is also valid for 
$n<r$, when all integral weights are singular.
}
\end{remark}


\section{Tables} \label{SECT8}

We illustrate our results by giving some tables of $q$-Littlewood-Richardson 
coefficients and of polynomials $d_{\alpha,\beta}(q)$.
These tables are $q$-analogues of those calculated by James
in \cite{Ja}, which were the starting point of our investigation.

\subsection{Canonical highest weight vectors of the Fock space
representation of $U_q(\widehat{\hbox{\twelvegoth sl}}_2)$}
\begin{figure}[t]
\begin{center}
\leavevmode
\epsfxsize =6cm
\epsffile{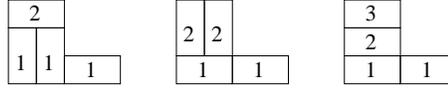}
\end{center}
\caption{\label{FIG10} The Yamanouchi domino tableaux of shape $(42^2)$}
\end{figure}
\noindent
The following tables give the coefficients $e_{2\alpha,\beta}(-q^{-1})$
of the expansion of $\G^-_{2\alpha}$ on the standard basis $\{ |\beta) \}$
for $n=2$ up to partitions of 10.
They should be read by columns, \eg
\[
\G^-_{(4)} = |4) - q^{-1}\,|3,1) + q^{-2} \,|2,2).
\]
These vectors form a basis of the subspace of primitive vectors of $\FI$.
Their coefficients are the $q$-analogues 
$c_{\mu^{(0)},\mu^{(1)}}^\lambda(-q^{-1})$ of the Littlewood-Richardson
multiplicities for all partitions $\mu^{(0)},\ \mu^{(1)}$ with
$|\mu^{(0)}|+|\mu^{(1)}|\le 5$.
They are easily computed using the combinatorial
description of \cite{CL} in terms of Yamanouchi domino tableaux. 
For example the row labelled $(42^2)$ is given by the tableaux of 
Figure~\ref{FIG10}.

{\small
$$
\begin {array}{rc}
&(2)\\
(2)&1\\
\noalign{\medskip}(1^2)&-q^{-1}
\end {array}
\hskip 1cm
\begin{array}{rcc}
&(4)&(2^2)\\
\noalign{\medskip}
(4)& 1&0\\
\noalign{\medskip}
(31)& -q^{-1} & 0\\
\noalign{\medskip}
(2^2)&q^{-2} & 1 \\
\noalign{\medskip}
(21^2)&0&-q^{-1} \\
\noalign{\medskip}
(1^4)&0&q^{-2}
\end{array}
\hskip 1cm
\begin{array}{rccc}
&(6)&(42)&(2^3)\\
\noalign{\medskip}
(6)&1&0&0\\
\noalign{\medskip}
(51)&-q^{-1}&0&0\\
\noalign{\medskip}
(42)&q^{-2}&1&0\\
\noalign{\medskip}
(41^2)&0&-q^{-1}&0\\
\noalign{\medskip}
(3^2)&-q^{-3}&-q^{-1}&0\\
\noalign{\medskip}
(31^3)&0&q^{-2}&0\\
\noalign{\medskip}
(2^3)&0&-q^{-2}&1\\
\noalign{\medskip}
(2^21^2)&0&-q^{-3}&-q^{-1}\\
\noalign{\medskip}
(21^4)&0&0&q^{-2}\\
\noalign{\medskip}
(1^6)&0&0&-q^{-3}
\end{array}
$$

$$
\begin{array}{rccccc}
&(8)&(62)&(4^2)&(42^2)&(2^4)\\
\noalign{\medskip}
(8)&1&0&0&0&0\\
\noalign{\medskip}
(71)&-q^{-1}&0&0&0&0\\
\noalign{\medskip}
(62)&q^{-2}&1&0&0&0\\
\noalign{\medskip}
(61^2)&0&-q^{-1}&0&0&0\\
\noalign{\medskip}
(53)&-q^{-3}&-q^{-1}&0&0&0\\
\noalign{\medskip}
(51^3)&0&q^{-2}&0&0&0\\
\noalign{\medskip}
(4^2)&q^{-4}&q^{-2}&1&0&0\\
\noalign{\medskip}
(431)&0&0&-q^{-1}&0&0\\
\noalign{\medskip}
(42^2)&0&q^{-2}&q^{-2}&1&0\\
\noalign{\medskip}
(421^2)&0&-q^{-3}&0&-q^{-1}&0\\
\noalign{\medskip}
(41^4)&0&0&0&q^{-2}&0\\
\noalign{\medskip}
(3^22)&0&-q^{-3}&0&-q^{-1}&0\\
\noalign{\medskip}
(3^21^2)&0&q^{-4}&q^{-2}&q^{-2}&0\\
\noalign{\medskip}
(32^21)&0&0&-q^{-3}&0&0\\
\noalign{\medskip}
(31^5)&0&0&0&-q^{-3}&0\\
\noalign{\medskip}
(2^4)&0&0&q^{-4}&q^{-2}&1\\
\noalign{\medskip}
(2^31^2)&0&0&0&-q^{-3}&-q^{-1}\\
\noalign{\medskip}
(2^21^4)&0&0&0&q^{-4}&q^{-2}\\
\noalign{\medskip}
(21^6)&0&0&0&0&-q^{-3}\\
\noalign{\medskip}
(1^8)&0&0&0&0&q^{-4}
\end{array}
$$
\small}

$$
\begin{array}{rccccccc}
&(10)&(82)&(64)&(62^2)&(4^22)&(42^3)&(2^5)\\
\noalign{\smallskip}
(10)&1&0&0&0&0&0&0\\
\noalign{\smallskip}
(91)&-q^{-1}&0&0&0&0&0&0\\
\noalign{\smallskip}
(82)&q^{-2}&1&0&0&0&0&0\\
\noalign{\smallskip}
(81^2)&0&-q^{-1}&0&0&0&0&0\\
\noalign{\smallskip}
(73)&-q^{-3}&-q^{-1}&0&0&0&0&0\\
\noalign{\smallskip}
(71^{3})&0&q^{-2}&0&0&0&0&0\\
\noalign{\smallskip}
(64)&q^{-4}&q^{-2}&1&0&0&0&0\\
\noalign{\smallskip}
(631)&0&0&-q^{-1}&0&0&0&0\\
\noalign{\smallskip}
(62^2)&0&q^{-2}&q^{-2}&1&0&0&0\\
\noalign{\smallskip}
(621^2)&0&-q^{-3}&0&-q^{-1}&0&0&0\\
\noalign{\smallskip}
(61^4)&0&0&0&q^{-2}&0&0&0\\
\noalign{\smallskip}
(5^2)&-q^{-5}&-q^{-3}&-q^{-1}&0&0&0&0\\
\noalign{\smallskip}
(532)&0&-q^{-3}&0&-q^{-1}&0&0&0\\
\noalign{\smallskip}
(531^2)&0&q^{-4}&q^{-2}&q^{-2}&0&0&0\\
\noalign{\smallskip}
(52^21)&0&0&-q^{-3}&0&0&0&0\\
\noalign{\smallskip}
(51^5)&0&0&0&-q^{-3}&0&0&0\\
\noalign{\smallskip}
(4^22)&0&q^{-4}&q^{-2}&q^{-2}&1&0&0\\
\noalign{\smallskip}
(4^21^2)&0&-q^{-5}&-q^{-3}&-q^{-3}&-q^{-1}&0&0\\
\noalign{\smallskip}
(43^2)&0&0&-q^{-3}&0&-q^{-1}&0&0\\
\noalign{\smallskip}
(431^3)&0&0&0&0&q^{-2}&0&0\\
\noalign{\smallskip}
(42^3)&0&0&q^{-4}&q^{-2}&q^{-2}&1&0\\
\noalign{\smallskip}
(42^21^2)&0&0&0&-q^{-3}&-q^{-3}&-q^{-1}&0\\
\noalign{\smallskip}
(421^4)&0&0&0&q^{-4}&0&q^{-2}&0\\
\noalign{\smallskip}
(41^6)&0&0&0&0&0&-q^{-3}&0\\
\noalign{\smallskip}
(3^31)&0&0&q^{-4}&0&q^{-2}&0&0\\
\noalign{\smallskip}
(3^22^2)&0&0&-q^{-5}&-q^{-3}&-q^{-3}&-q^{-1}&0\\
\noalign{\smallskip}
(3^221^2)&0&0&0&q^{-4}&0&q^{-2}&0\\
\noalign{\smallskip}
(3^21^4)&0&0&0&-q^{-5}&-q^{-3}&-q^{-3}&0\\
\noalign{\smallskip}
(32^21^3)&0&0&0&0&q^{-4}&0&0\\
\noalign{\smallskip}
(31^7)&0&0&0&0&0&q^{-4}&0\\
\noalign{\smallskip}
(2^5)&0&0&0&0&q^{-4}&q^{-2}&1\\
\noalign{\smallskip}
(2^41^2)&0&0&0&0&-q^{-5}&-q^{-3}&-q^{-1}\\
\noalign{\smallskip}
(2^31^4)&0&0&0&0&0&q^{-4}&q^{-2}\\
\noalign{\smallskip}
(2^21^6)&0&0&0&0&0&-q^{-5}&-q^{-3}\\
\noalign{\smallskip}
(21^8)&0&0&0&0&0&0&q^{-4}\\
\noalign{\smallskip}
(1^{10})&0&0&0&0&0&0&-q^{-5}
\end{array}
$$

\subsection{\large\bf Basis $\{\G^+_\beta\}$ of the Fock space 
representation of $U_q(\widehat{\hbox{\twelvegoth sl}}_2)$}

The following tables give the coefficients $d_{\alpha,\beta}(q)$
of the expansion of $\G^+_\beta$ on the standard basis $\{ |\alpha) \}$
for $n=2$ up to partitions of 10.
They should be read by columns, \eg
\[
\G^+_{(3,1)} = |3,1) + q\,|2,2) + q^2 \,|2,1,1).
\]
Each square matrix corresponds to a weight space of $\FI$.
(The weight space containing $|10)$ being too large, the
corresponding matrix had to be displayed on two pages.)
The 1-dimensional weight spaces corresponding to the 
partitions $(1),\ (2,1),\ (3,2,1),\ (4,3,2,1)$ have been omitted.

\vskip 2cm

$$
\begin {array}{rcc}   
(2)&1&0\\ 
\noalign{\medskip}(1^2)&q&1
\end {array}
\hskip 7 mm
\begin {array}{rcc}  
(3)&1&0\\ \noalign{\medskip}(1^3)&q&1
\end {array}
\hskip 7mm
\begin {array}{rccccc}  
(4)&1&0&0&0&0\\ 
\noalign{\medskip}(3 1)&q&1&0&0&0\\ 
\noalign{\medskip}(2^2)&0&q&1&0&0 \\ 
\noalign{\medskip}(2 1^2)&q&q^{2}&q&1&0\\ 
\noalign{\medskip}(1^4)&q^{2}&0&0&q&1\end {array}
\hskip 7mm
\begin {array}{rccccc}  
(5)&1&0&0&0&0\\ 
\noalign{\medskip}(3 2)&0&1&0&0&0\\ 
\noalign{\medskip}(3 1^2)&q&q&1&0&0\\ 
\noalign{\medskip}(2^2 1)&0&q^{2}&q&1&0\\  
\noalign{\medskip}(1^5)&q^{2}&0&q&0&1
\end{array}
$$
\vskip 2cm 

$$
\begin {array}{rcc}
(4 1)&1&0 \\ 
\noalign{\medskip}(2 1^3)&q&1
\end{array}
\hskip 20mm 
\begin {array}{rcccccccccc}  
(6)&1&0&0&0&0&0&0&0&0&0\\ 
\noalign{\medskip}(5 1)&q&1&0&0&0&0&0&0&0&0\\ 
\noalign{\medskip}(4 2)&0&q&1&0&0&0&0&0&0&0\\ 
\noalign{\medskip}(4 1^2)&q&q^{2}&q&1&0&0&0&0&0&0\\ 
\noalign{\medskip}(3^2)&0&0&q&0&1&0&0&0&0&0\\ 
\noalign{\medskip}(3 1^3)&q^{2}&q&q^{2}&q&q&1&0&0&0&0\\ 
\noalign{\medskip}(2^3)&0&0&q^{2}&q&q&0&1&0&0&0\\ 
\noalign{\medskip}(2^2 1^2)&0 &q^{2}&q^{3}&q^2&q^2&q&q&1&0&0 \\ 
\noalign{\medskip}(2 1^4)&q^{2}&q^{3}&0&q&0&q^2&0&q&1&0\\ 
\noalign{\medskip}(1^6)&q^{3}&0&0&q^2&0&0&0&0&q&1
\end {array}
$$

\vskip 10mm
$$
\begin {array}{rcccccccccc}  
(7)&1&0&0&0&0&0&0&0&0&0\\ 
\noalign{\medskip}(5 2)&0&1&0&0&0&0&0&0&0&0\\ 
\noalign{\medskip}(5 1^2 )&q&q&1&0&0&0&0&0&0&0\\ 
\noalign{\medskip}(4 2 1)&0&q^{2}&q&1&0&0&0&0&0&0\\ 
\noalign{\medskip}(3 ^2 1)&0&0&0&q&1&0&0&0&0&0 \\ 
\noalign{\medskip}(3 2^2 )&0&0&q&q^{2}&q&1&0&0&0&0\\ 
\noalign{\medskip}(3 2 1^2)&0&q&q^2&q^{3}&q^2&q&1&0&0&0\\ 
\noalign{\medskip}(3 1^4)&q^{2}& q^{2}&q&0&0&0&q&1&0&0\\ 
\noalign{\medskip}(2^2 1^3)&0&q^{3}&q^2&0&0&q&q^2&q&1&0\\ 
\noalign{\medskip}(1^7)&q^{3}&0&q^2&0&0&0&0&q&0&1
\end {array}
\hskip 10mm
\begin{array}{l}
\begin {array}{rccccc}
(6 1)&1&0&0&0&0\\
\noalign{\medskip}(4 3)&0&1&0&0&0\\
\noalign{\medskip}(4 1^3)&q&q&1&0&0\\
\noalign{\medskip}(2^3 1)&0&q^{2}&q&1&0\\
\noalign{\medskip}(2 1^5)&q^{2}&0&q&0&1
\end {array}\\[1cm]
\\
\\
\begin {array}{rcc}
(5 2 1)&1&0 \\
\noalign{\medskip}(3 2 1^3)&q&1
\end{array}
\end{array}
$$
\vskip 10mm

\small
$$
\begin {array}{rcccccccccccccccccccc} 
(8)&1&0&0&0&0&0&0&0&0&0&0&0&0&0&0&0&0&0&0&0\\
\noalign{\medskip}
(71)&q&1&0&0&0&0&0&0&0&0&0&0&0&0&0&0&0&0&0&0\\
\noalign{\medskip}
(62)&0&q&1&0&0&0&0&0&0&0&0&0&0&0&0&0&0&0&0&0\\
\noalign{\medskip}(61^2)&q&{q}^{2}&q&1&0&0&0&0&0&0&0&0&0&0&0&0&0&0&0&0\\
\noalign{\medskip}(53)&0&0&q&0&1&0&0&0&0&0&0&0&0&0&0&0&0&0&0&0\\
\noalign{\medskip}(51^3)&{q}^{2}&q&{q}^{2}&q&q&1&0&0&0&0&0&0&0&0&0&0&0&0&0&0\\
\noalign{\medskip}(4^2)&0&0&0&0&q&0&1&0&0&0&0&0&0&0&0&0&0&0&0&0\\
\noalign{\medskip}(431)&0&0&q&0&{q}^{2}&0&q&1&0&0&0&0&0&0&0&0&0&0&0&0\\
\noalign{\medskip}(42^2)&0&0&{q}^{2}&q&0&0&0&q&1&0&0&0&0&0&0&0&0&0&0&0\\
\noalign{\medskip}(421^2)&0&{q}^{2}&q{+}{q}^{3}&{q}^{2}&{q}^{2}&q&q&{q}^{2}
&q&1&0&0&0&0&0&0&0&0&0&0\\
\noalign{\medskip}(41^4)&{q}^{2}&{q}^{3}&{q}^{2}&q&{q}^{3}&{q}^{2}&{
q}^{2}&0&0&q&1&0&0&0&0&0&0&0&0&0\\
\noalign{\medskip}(3^22)&0&0&0&0&0
&0&q&{q}^{2}&q&0&0&1&0&0&0&0&0&0&0&0\\
\noalign{\medskip}(3^21^2)&0&0&{q}
^{2}&0&0&0&{q}^{2}&{q}^{3}&{q}^{2}&q&0&q&1&0&0&0&0&0&0&0
\\
\noalign{\medskip}(32^21)&0&0&{q}^{3}&{q}^{2}&{q}^{2}&q&q{+}{q}^{3}&{q}
^{4}&q{+}{q}^{3}&{q}^{2}&0&{q}^{2}&q&1&0&0&0&0&0&0
\\
\noalign{\medskip}(31^5)&{q}^{3}&{q}^{2}&{q}^{3}&{q}^{2}&0&q&0&0&0&{
q}^{2}&q&0&q&0&1&0&0&0&0&0\\
\noalign{\medskip}(2^4)&0&0&0&0&{q}^{3}
&{q}^{2}&{q}^{2}&0&0&0&0&0&0&q&0&1&0&0&0&0\\
\noalign{\medskip}(2^31^2)&0
&0&{q}^{3}&{q}^{2}&{q}^{4}&{q}^{3}&{q}^{3}&0&q&{q}^{2}&q&0&q&{q
}^{2}&0&q&1&0&0&0\\
\noalign{\medskip}(2^21^4)&0&{q}^{3}&{q}^{4}&{q}^{3}&0
&{q}^{2}&0&0&{q}^{2}&{q}^{3}&{q}^{2}&0&{q}^{2}&0&q&0&q&1&0&0
\\
\noalign{\medskip}(21^6)&{q}^{3}&{q}^{4}&0&{q}^{2}&0&{q}^{3}&0&0&0&0
&q&0&0&0&{q}^{2}&0&0&q&1&0\\
\noalign{\medskip}(1^8)&{q}^{4}&0&0&{q}^{3}&0&0&0&0&0&0&{q}^{2}
&0&0&0&0&0&0&0&q&1\end {array}
$$

\newpage \small

$$
\begin {array}{rcccccccccccccccccccc} 
(9)&1&0&0&0&0&0&0&0&0&0&0&0&0&0&0&0&0&0&0&0 \\
\noalign{\medskip}(72)&0&1&0&0&0&0&0&0&0&0&0&0&0&0&0&0&0&0&0&0\\
\noalign{\medskip}(71^2)&q&q&1&0&0&0&0&0&0&0&0&0&0&0&0&0&0&0&0&0\\
\noalign{\medskip}(621)&0&{q}^{2}&q&1&0&0&0&0&0&0&0&0&0&0&0&0&0&0&0&0\\
\noalign{\medskip}(54)&0&0&0&0&1&0&0&0&0&0&0&0&0&0&0&0&0&0&0&0\\
\noalign{\medskip}(531)&0&0&0&q&q&1&0&0&0&0&0&0&0&0&0&0&0&0&0&0\\
\noalign{\medskip}(52^2)&0&0&q&{q}^{2}&0&q&1&0&0&0&0&0&0&0&0&0&0&0&0&0\\
\noalign{\medskip}(521^2)&0&q&{q}^{2}&{q}^{3}&q&{q}^{2}&q&1
&0&0&0&0&0&0&0&0&0&0&0&0\\
\noalign{\medskip}(51^4)&{q}^{2}&{q}^{2}&q&0&{q}^{2}&0&0&q&1&0&0
&0&0&0&0&0&0&0&0&0\\
\noalign{\medskip}(4^21)&0&0&0&0&{q}^{2}&q&0&0&0&1&0&0&0&0&0&0&0&0&0&0\\
\noalign{\medskip}(421^3)&0&{q}^{3}&{q}^{2}&q&{q}^{3}&{q}^{2}&q&{
q}^{2}&q&q&1&0&0&0&0&0&0&0&0&0\\
\noalign{\medskip}(3^3)&0&0&0&0&0&{q}^{2}&q&0&0&q&0
&1&0&0&0&0&0&0&0&0\\
\noalign{\medskip}(3^21^3)&0&0&0&{q}^{2}&0&{q}^{3}&{q}^{2}&0&0&{q}^{2}
&q&q&1&0&0&0&0&0&0&0\\
\noalign{\medskip}(32^3)&0
&0&0&0&{q}^{2}&{q}^{3}&{q}^{2}&q&0&{q}^{2}&0&q&0&1&0&0&0&0&0&0\\
\noalign{\medskip}(32^21^2)&0&0&{q}^{2}&{q}^{3}&{q}^{3}&{q}^{4}
&q{+}{q}^{3}&{q}^{2}&q&{q}^{3}&{q}^{2}&{q}^{2}&q&q&
1&0&0&0&0&0\\
\noalign{\medskip}(321^4)&0&{q}^{2}&{q}^{3}&{q}^{4}
&0&0&{q}^{2}&q&{q}^{2}&0&{q}^{3}&0&{q}^{2}&0&q&1&0&0&0&0\\
\noalign{\medskip}(31^6)&{q}^{3}&{q}^{3}&{q}^{2}&0&0&0&0&{q}
^{2}&q&0&0&0&0&0&0&q&1&0&0&0\\
\noalign{\medskip}(2^41)&0
&0&0&0&{q}^{4}&0&0&{q}^{3}&{q}^{2}&0&0&0&0&{q}^{2}&
q&0&0&1&0&0\\
\noalign{\medskip}(2^21^5)&0&{q}^{4}&{q}^{3}&0&0&0&{q}^{2}&{q}^{3}
&{q}^{2}&0&0&0&0&0&q&{q}^{2}&q&0&1&0\\
\noalign{\medskip}(1^9)&{q}^{4}&0&{q}^{3}
&0&0&0&0&0&{q}^{2}&0&0&0&0&0&0&0&q&0&0&1
\end {array}
$$
\vskip 10mm

\normalsize
$$
\begin {array}{rcccccccccc} 
(81)  &1&0&0&0&0&0&0&0&0&0\\
\noalign{\medskip}(63)  &0&1&0&0&0&0&0&0&0&0\\
\noalign{\medskip}(61^3)&q&q&1&0&0&0&0&0&0&0\\
\noalign{\medskip}(432)  &0&0&0&1&0&0&0&0&0&0\\
\noalign{\medskip}(431^2)&0&q&0&q&1&0&0&0&0&0\\
\noalign{\medskip}(42^21)&0&{q}^{2}&q&{q}^{2}&q&1&0&0&0&0\\
\noalign{\medskip}(41^5)&{q}^{2}&{q}^{2}&q&0&q&0&1&0&0&0\\
\noalign{\medskip}(3^221)&0&0&0&{q}^{3}&{q}^{2}&q&0&1&0&0\\
\noalign{\medskip}(2^31^3)&0&{q}^{3}&{q}^{2}&0
&{q}^{2}&q&q&0&1&0\\
\noalign{\medskip}(21^7)&{q}^{3}&0&{q}^{2}&0&0&0&q&0&0&1
\end {array}
\hskip 10mm
\begin {array}{rccccc}
(721)&1&0&0&0&0\\
\noalign{\medskip}(541)&0&1&0&0&0\\
\noalign{\medskip}(521^3)&q&q&1&0&0\\
\noalign{\medskip}(32^31)&0&q^{2}&q&1&0\\
\noalign{\medskip}(321^5)&q^{2}&0&q&0&1
\end{array}
$$

$$
\small
\begin {array}
{rcccccccccccccccccc} 
(10)&1&0&0&0&0&0&0&0&0&0&0&0&0&0&0&0&0&0 \\
\noalign{\smallskip}(91)&q&1&0&0&0&0&0&0&0&0&0&0&0&0&0&0&0&0\\
\noalign{\smallskip}(82)        &0&q&1&0&0&0&0&0&0&0&0&0&0&0&0&0&0&0\\
\noalign{\smallskip}(81^2)&q&{q}^{2}&q&1&0&0&0&0&0&0&0&0&0&0&0&0&0&0\\
\noalign{\smallskip}(73)&0&0&q&0&1&0&0&0&0&0&0&0&0&0&0&0&0&0\\
\noalign{\smallskip}(71^3)&{q}^{2}&q&{q}^{2}&q&q&1&0&0&0&0&0
&0&0&0&0&0&0&0\\
\noalign{\smallskip}(64)&0&0&0&0&q&0&1&0&0&0&0&0&0&0&0&0&0&0\\
\noalign{\smallskip}(631)&0&0&q&0&{q}^{2}&0&q&1&0&0&0&0&0&0&0&0&0&0\\
\noalign{\smallskip}(622)&0&0&{q}^{2}&q&0&0&0&q&1&0&0&0&0&0&0&0&0&0\\
\noalign{\smallskip}(621^2)&0&{q}^{2}&q{+}{q}^{3}&{q}^{2}&{q}^{2}&q&q&{q}
^{2}&q&1&0&0&0&0&0&0&0&0\\
\noalign{\smallskip}(61^4)&{q}^{2}&{q}^{3}&{q}^{2}&q&{q}^{3}&{q}^{2}
&{q}^{2}&0&0&q&1&0&0&0&0&0&0&0\\
\noalign{\smallskip}(5^2)&0&0&0&0&0&0&q&0&0&0&0&1&0&0&0&0&0&0\\
\noalign{\smallskip}(532)&0&0&0&0&0&0&q&{q}^{2}&q&0&0&0&1&0&0&0&0&0\\
\noalign{\smallskip}(531^2)&0&0&{q}^{2}&0&q&0&2\,{q}^{2
}&{q}^{3}&{q}^{2}&q&0&q&q&1&0&0&0&0\\
\noalign{\smallskip}(52^21)&0&0&{q}^{3}&{q}^{2}&{q}^{2}
&q&{q}^{3}&{q}^{4}&q{+}{q}^{3}&{q}^{2}&0&0&{q}^{2}&q&1&0&0&0\\
\noalign{\smallskip}(51^5)&{q}^{3}&{q}^{2}&{q}^{3}&{q}^{2}&{q}
^{2}&q&{q}^{3}&0&0&{q}^{2}&q&{q}^{2}&0&q&0&1&0&0\\
\noalign{\smallskip}(4^22)&0&0&0&0&0&0&{q}^{2}&0&0&0&0&q&q&0&0&0&1&0\\
\noalign{\smallskip}(4^21^2)&0&0&0&0&{q}^{2}&0&{q}^{3}&0&0&0&0&{
q}^{2}&{q}^{2}&q&0&0&q&1\\
\noalign{\smallskip}(43^2)&0&0&0&0&0&0&0&0&q&0&0&0&{q}^{2}&0&0
&0&q&0\\
\noalign{\smallskip}(431^3)&0&0&{q}^{2}&0&{q}^{3}&0&{q}^{2}&q
&{q}^{2}&q&0&q&{q}^{3}&{q}^{2}&0&0
&{q}^{2}&q \\
\noalign{\smallskip}(42^3)&0&0&0&0&{q}^{3}&{q}^{2}&{q}^{2}&0&{q}^{2}&q
&0&q&{q}^{3}&{q}^{2}&q&0&{q}^{2}&q\\
\noalign{\smallskip}(42^21^2)&0&0&{q}^{3}&{q}^{2}&{q}^{4}&{
q}^{3}&2\,{q}^{3}&{q}^{2}&q{+}{q}^{3}&2\,{q}^{2}&q&2\,{q}^{2}&{q}
^{4}&{q}^{3}&{q}^{2}&0&{q}^{3}&{q}^{2}\\
\noalign{\smallskip}(421^4)&0&{q}^{3}&{q}^{4}{+}{q}^{
2}&{q}^{3}&{q}^{3}&{q}^{2}&{q}^{4}&{q}^{3}&{q}^{2}&q{+}{q}^{3}&{q
}^{2}&{q}^{3}&0&{q}^{2}&0&q&0&q\\
\noalign{\smallskip}(41^6)&{q}^{3}&{q}^{4}&{q}^{3}&{q}
^{2}&{q}^{4}&{q}^{3}&0&0&0&{q}^{2}&q&0&0&{q}^{3}&0&{q}^{2}&0
&{q}^{2} \\
\noalign{\smallskip}(3^31)&0&0&0&0&0&0&0&0&{q}^{2}&0&0&q&{q}^{3}&{q}
^{2}&q&0&{q}^{2}&q \\
\noalign{\smallskip}(3^22^2)&0&0&0&0&0&0&{q}^{3}&0&{q}^{3}&{q}^{2}&0&2\,
{q}^{2}&{q}^{4}&{q}^{3}&{q}^{2}&0&{q}^{3}&{q}^{2}\\
\noalign{\smallskip}(3^221^2)&0&0&0&0
&0&0&{q}^{4}&{q}^{3}&{q}^{4}{+}{q}^{2}&{q}^{3}&0&2\,{q}^{3}&{q}^{
5}&{q}^{4}&{q}^{3}&0&{q}^{4}&q{+}{q}^{3}\\
\noalign{\smallskip}(3^21^4)&0&0
&{q}^{3}&0&0&0&0&{q}^{4}&{q}^{3}&{q}^{2}&0&0&0&0&0&0&0&{q}^{2}\\
\noalign{\smallskip}(32^21^3)&0
&0&{q}^{4}&{q}^{3}&{q}^{3}&{q}^{2}&{q}^{4}&{q}^{5}&{q}^{4}{+}{q}^
{2}&2\,{q}^{3}&{q}^{2}&{q}^{3}&0&{q}^{2}&q&q&0&q{+}{q}^{3}\\
\noalign{\smallskip}(31^7)&{q}^{4}&{q}^{3}&{q}^{4}&{q}^{3}&0&{q}^{2}&0
&0&0&{q}^{3}&{q}^{2}&0&0&0&0&q&0&0\\
\noalign{\smallskip}(2^5)&0&0&0&0&0&0&{q}^{4}&0&0
&{q}^{3}&{q}^{2}&{q}^{3}&0&0&0&0&0&0\\
\noalign{\smallskip}(2^41^2)&0&0&0&0&{q}^{4}&{q}^{
3}&{q}^{5}&0&0&{q}^{4}&{q}^{3}&{q}^{4}&0&{q}^{3}&{q}^{2}&{q}^
{2}&0&{q}^{2}\\
\noalign{\smallskip}(2^31^4)&0&0&{q}^{4}&{q}^{3}&{q}^{5}&{q}^{4}
&0&0&{q}^{2}&{q}^{3}&{q}^{2}&0&0&{q}^{4}&{q}^{3}&{q}^{3}&0&{q}^{3}
\\
\noalign{\smallskip}(2^21^6)&0&{q}^{4}&{q}^{5}&{q}^{4}&0&{q}^{3}&0&0&{q}
^{3}&{q}^{4}&{q}^{3}&0&0&0&0&{q}^{2}&0&0\\
\noalign{\smallskip}(21^8)&
{q}^{4}&{q}^{5}&0&{q}^{3}&0&{q}^{4}&0&0&0&0&{q}^{2}&0&0&0&0
&{q}^{3}&0&0 \\
\noalign{\smallskip}(1^{10})&{q}^{5}&0&0&{q}^{4}&0&0&0&0&0&0&{q}^{3}&0&0
&0&0&0&0&0
\end {array}
$$
$$
\begin {array}
{rcccccccccccccccccc}
(43^2)&1&0&0&0&0&0&0&0&0&0&0&0&0&0&0&0&0&0 \\
\noalign{\smallskip}(431^2) &q&1&0&0&0&0&0&0&0&0&0&0&0&0&0&0&0&0\\
\noalign{\smallskip}(42^3)  &q&0&1&0&0&0&0&0&0&0&0&0&0&0&0&0&0&0\\
\noalign{\smallskip}(42^21^2) &{q}^{2}&q&q&1&0&0&0&0&0&0&0&0&0&0&0&0&0&0\\
\noalign{\smallskip}(421^4) &0&{q}^{2}&0&q&1&0&0&0&0&0&0&0&0&0&0&0&0&0\\
\noalign{\smallskip}(41^6) &0&0&0&0&q&1&0&0&0&0&0&0&0&0&0&0&0&0 \\
\noalign{\smallskip}(3^31) &q&0&0&0&0&0&1&0&0&0&0&0&0&0&0&0&0&0\\
\noalign{\smallskip}(3^22^2) &{q}^{2}&0&q&0&0&0&q&1&0&0&0&0&0&0&0&0&0&0\\
\noalign{\smallskip}(3^221^2) &q{+}{q}^{3}&{q}^{2}&{q}^{2}&q&0&0&{q}^{2}&q&1
&0&0&0&0&0&0&0&0&0\\
\noalign{\smallskip}(3^21^4) &{q}^{2}&{q}^{3}&0&{q}^{2}&q&0&0&0&q&1
&0&0&0&0&0&0&0&0 \\
\noalign{\smallskip}(32^21^3)
&{q}^{3}&{q}^{4}&{q}^{2}&q{+}{q}^{3}&{q}^{2}&0&0&q&{q}^{2}&q&1&0&0&0&0&0&0&0\\
\noalign{\smallskip}(31^7) &0&0&0&0&{q}^{2}&q&0&0&0&q&0&1&0&0&0&0&0&0\\
\noalign{\smallskip}(2^5) &0&0&{q}^{2}&q&0&0&0&q&0&0&0&0&1&0&0&0&0&0\\
\noalign{\smallskip}(2^41^2)&0&0&{q}^{3}&{q}^{2}&0&0&0&{q}^{2}&0&0&q&0&q&1
&0&0&0&0\\
\noalign{\smallskip}(2^31^4) &0&0&0&q&{q}^{2}&q&0&0&0&q&{q}^{2}&0&0&q&1&0&0&0
\\
\noalign{\smallskip}(2^21^6) &0&0&0&{q}^{2}&{q}^{3}&{q}^{2}
&0&0&0&{q}^{2}&0&q&0&0&q&1&0&0\\
\noalign{\smallskip}(21^8)&0&0&0&0&0&q&0&0&0&0&0&{q}^{2}&0&0&0&q&1&0 \\
\noalign{\smallskip}(1^{10}) &0&0&0&0&0&{q}^{2}&0&0&0&0&0&0&0&0&0&0&q&1
\end {array}
$$

%

\normalsize
\vskip 2cm

\subsection*{Acknowledgements}
We want to thank B. Feigin and W. Soergel for helpful and stimulating 
discussions.
This work was done during our stay at R.I.M.S. Kyoto University
for the R.I.M.S. Project 1998 
``Combinatorial methods in representation theory''.
We would like to thank the organizers 
K. Koike, M. Kashiwara, S. Okada,  H.-F. Yamada, I. Terada
for the invitation, and the R.I.M.S. for its warm hospitality. 

\bigskip\bigskip


\bigskip\bigskip
\begin{thebibliography}{ABC} \scriptsize
%
\bibitem{BV} {\sc D. Barbasch, D. Vogan}, {\it Primitive ideals and orbital 
integrals in complex classical groups}, Math. Ann. 259, (1982) 153-199.
%
\bibitem{CL} {\sc C. Carr\'e, B. Leclerc}, {\it Splitting the square of a 
Schur function into its symmetric and antisymmetric parts}, 
J. Algebraic Combinatorics, {\bf 4} (1995), 201-231.
%
\bibitem{CP} {\sc V. Chari, A. Pressley}, {\it A guide to quantum groups},
Cambridge Univ. Press. 1995.
%
\bibitem{De1} {\sc V. V. Deodhar}, {\it On some geometric aspects of
Bruhat orderings II. The parabolic analogue of Kazhdan-Lusztig
polynomials}, J. Algebra, {\bf 111} (1987), 483-506.
%
\bibitem{De2} {\sc V. V. Deodhar}, {\it Duality in parabolic set up for
questions in Kazhdan-Lusztig theory}, 
J. Algebra, {\bf 142} (1991), 201-209.
%
\bibitem{Du} {\sc J. Du}, {\it IC bases and quantum linear groups},
Proc. Symp. Pure Math. {\bf 56} (1994), part 2, 135-148.
%
\bibitem{DPS} {\sc J. Du, B. Parshall, L. Scott}, 
{\it Quantum Weyl reciprocity and tilting modules},
Commun. Math. Phys. {\bf 195} (1998), 321-352.
%
\bibitem{Fu} {\sc W. Fulton}, {\it Eigenvalues of sums of 
Hermitian matrices [after A. Klyachko]}, S\'eminaire Bourbaki, Juin 1998. 
%
\bibitem{Ga} {\sc  D. Garfinkle}, {\it On the classification of primitive 
ideals for complex classical Lie algebras, I}, Compositio Mathematica, 75 
(1990) 2, 135-169.
%
\bibitem{GW} {\sc F. Goodman, H. Wenzl}, {\it Crystal bases of 
quantum affine algebras and affine Kazhdan-Lusztig polynomials},
math.QA/9807014.
%
\bibitem{IM} {\sc N. Iwahori, H. Matsumoto}, {\it On some Bruhat 
decomposition and the structure of the Hecke ring of $p$-adic
Chevalley groups}, Publ. IHES, {\bf 25} (1965), 5-48.
%
\bibitem{Ja} {\sc G. James}, {\it The decomposition matrices of $GL_n(q)$
for $n\le 10$}, Proc. London Math. Soc., {\bf 60} (1990), 225-265.
%
\bibitem{JK} {\sc G. James, A. Kerber}, {\it The representation theory of the
symmetric group}, Addison Wesley, 1981.
%
\bibitem{Jim} {\sc M. Jimbo}, {\it A $q$-analogue of $U(gl(N+1))$,
Hecke algebra and the Yang-Baxter equation}, Lett. Math. Phys.
{\bf 11} (1986), 247-252.
%
\bibitem{KMS}{\sc M. Kashiwara, T. Miwa, E. Stern}, {\it Decomposition
of $q$-deformed Fock spaces},
Selecta Math. {\bf 1}  (1996) 787.
%
\bibitem{KLLT} {\sc A. N. Kirillov, A. Lascoux, B. Leclerc, J.-Y. Thibon}, 
{\it S\'eries g\'en\'eratrices pour les tableaux de dominos}, 
C. R. Acad. Sci. Paris, t. 318, S\'erie I, (1994) 395-400.
%
\bibitem{KSS} {\sc A. N. Kirillov, A. Schilling,  M. Shimozono},
{\it On a bijection from Littlewood-Richardson tableaux to rigged
configurations}, (in preparation).
%
\bibitem{KS} {\sc A.N. Kirillov, M. Shimozono}, {\it A generalization
of Kostka-Foulkes polynomials}, math.QA/9803062.
%
\bibitem{K} {\sc A.A. Klyachko}, {\it Stable vector bundles and 
hermitian operators}, IGM, Universit\'e de Marne-la-Vall\'ee,
Preprint 1994.
%
\bibitem{LLT} {\sc A. Lascoux, B. Leclerc, J.-Y. Thibon}, 
{\it Hecke algebras at roots of unity 
and crystal bases of quantum affine
algebras}, Commun. Math. Phys. {\bf 181} (1996), 205-263.
%
\bibitem{LLT2} {\sc A. Lascoux, B. Leclerc, J.-Y. Thibon},
{\it Ribbon tableaux, Hall-Littlewood functions, quantum
affine algebras, and unipotent varieties},
J. Math. Phys. {\bf 38} (1997), 1041-1068.
%
\bibitem{LT} {\sc B. Leclerc, J.-Y. Thibon}, {\it Canonical bases
of $q$-deformed Fock spaces}, Int. Math. Res. Notices, {\bf 9} (1996),
447-456.
%
\bibitem{vL} {\sc M. van Leeuwen}, {\em A Robinson-Schensted
algorithm in the geometry of flags for Classical Groups}, Thesis, 1989.
%
\bibitem{vL2} {\sc M. van Leeuwen}, {\em Some bijective correspondences
involving domino tableaux}, Preprint CWI MAS-R9708, 1997.
%
\bibitem{Lit} {\sc D. E. Littlewood}, 
{\it Modular representations of symmetric groups}, 
Proc. Roy. Soc. {\bf 209} (1951), 333-353.
%
\bibitem{Lu0} {\sc G. Lusztig}, {\it Some problems in the
representation theory of a finite Chevalley group}, 
Proc. Symp. Pure Math. AMS, {\bf 37} (1980), 313-317. 
%
\bibitem{Lu1} {\sc G. Lusztig}, {\it Green polynomials and singularities
of unipotent classes}, Advances in Math. {\bf 42} (1981), 169-178.
%
\bibitem{Lu2} {\sc G. Lusztig}, {\it Singularities, character formulas,
and a $q$-analog of weight multiplicities}, Analyse et topologie sur les 
espaces singuliers (II-III), Ast\'erisque {\bf 101}-{\bf 102}
(1983), 208-227.
%
\bibitem{LuAMS} {\sc G. Lusztig}, {\it Some examples of square integrable 
representations of semisimple $p$-adic groups},
Trans. AMS {\bf 227} (1983), 623-653.
%
\bibitem{Lu4} {\sc G. Lusztig}, {\it Modular representations and
quantum groups}, Contemp. Math. {\bf 82} (1989), 58-77.
%
\bibitem{Lu5} {\sc G. Lusztig}, {\it On quantum groups},
J. Algebra, {\bf 131} (1990), 466-475.
%
\bibitem{Lu6}{\sc G. Lusztig}, {\it Canonical bases arising from quantized
enveloping algebras}, J. Amer. Math. Soc. {\bf 3} (1990), 447-498.
%
\bibitem{Mcd} {\sc I. G. Macdonald}, {\it Symmetric functions and 
Hall polynomials}, 2nd edition, Oxford 1995.
%
%
\bibitem{SchWa} {\sc A. Schilling, O. Warnaar}, {\it Inhomogeneous
lattice paths, generalized Kostka-Foulkes polynomials, and 
$A_{n-1}$-supernomials}, math.QA/9802111.
%
\bibitem{Sc} {\sc M.P. Sch\"utzenberger}, {\it Propri\'et\'es nouvelles 
des tableaux de Young}, S\'eminaire Delange-Pisot-Poitou, 19{\`eme ann\'ee}, 
{\bf 26}, 1977/78.
%
\bibitem{ShW} {\sc M. Shimozono, J. Weyman}, {\it Characters of 
modules supported in the closure of a nilpotent conjugacy class},
math.QA/9804036.
%
\bibitem{Sh1} {\sc M. Shimozono}, {\it A cyclage poset structure
for Littlewood-Richardson tableaux}, math.QA/9804037.
%
\bibitem{Sh2} {\sc M. Shimozono}, {\it Affine type A crystal structure 
on tensor products of rectangles, Demazure characters, and
nilpotent varieties}, math.QA/9804039.
%
\bibitem{So1} {\sc W. Soergel}, {\it Kazhdan-Lusztig-Polynome und 
eine Kombinatorik f\"ur 
Kipp-Moduln}, Represent. Theory 1 (1997), 37-68 (english 83-114).
%
\bibitem{So2} {\sc W. Soergel},  {\it Charakterformeln f\"ur 
Kipp-Moduln \"uber Kac-Moody-Algebren},
Represent. Theory 1 (1997), 115-132.
%
\bibitem{SW} {\sc D. Stanton,  D. White}, {\em A Schensted algorithm for 
rim-hook tableaux}, J. Comb. Theory A 40, (1985), 211-247.
%
\bibitem{VV} {\sc M. Varagnolo, E. Vasserot}, {\it Canonical bases
and the Lusztig conjecture for quantized $sl(n)$ at roots of
unity}, math.QA/9803023.
%
\bibitem{Z} {\sc A. Zelevinsky}, {\it Littlewood-Richardson semigroups},
MSRI Berkeley, Preprint 1997.
%
\end{thebibliography}
\end{document}